\documentclass[11pt]{article}
\usepackage{fancyhdr}
\usepackage[english]{babel}
\usepackage[utf8]{inputenc}
\usepackage{amsmath}
\usepackage{amssymb}
\usepackage{amsthm}
\usepackage{bbm}
\usepackage{tikz}
\usetikzlibrary{patterns}
\usepackage{enumerate}
\usepackage{graphicx}
\graphicspath{ {./Figures/} }
\usepackage{enumerate}
\numberwithin{equation}{section}
\usepackage[backend = biber, sorting = nyt, defernumbers=true, style = numeric, url = false, isbn = false, doi =true]{biblatex}
\usepackage{url}
\usepackage[margin = 1.5cm, a4paper]{geometry}
\usepackage{csquotes}
\usepackage{mathtools}

\usepackage{pgfplots}
\usepackage{tkz-euclide}
\usepackage{subcaption}
\usepackage{orcidlink}
\usepackage{longtable}
\RequirePackage{doi}

\usepackage{hyperref}
\usetikzlibrary{decorations.pathreplacing}
\usetikzlibrary{decorations.pathmorphing}

\DeclareMathOperator{\Var}{Var}

\newcommand{\R}{\mathbb{R}}
\newcommand{\N}{\mathbb{N}}
\newcommand{\Z}{\mathbb{Z}}

\newcommand{\PlusUps}{+}
\newcommand{\MinusUps}{-}

\makeatletter
\newcommand*\wt[1]{\mathpalette\wthelper{#1}}
\newcommand*\wthelper[2]{%
        \hbox{\dimen@\accentfontxheight#1%
                \accentfontxheight#11.075\dimen@
                $\m@th#1\widetilde{#2}$%
                \accentfontxheight#1\dimen@
        }%
}

\newcommand*\accentfontxheight[1]{%
        \fontdimen5\ifx#1\displaystyle
                \textfont
        \else\ifx#1\textstyle
                \textfont
        \else\ifx#1\scriptstyle
                \scriptfont
        \else
                \scriptscriptfont
        \fi\fi\fi3
}
\makeatother

\DeclareMathOperator{\es}{es}
\DeclareMathOperator*{\esssup}{ess\,sup}

\theoremstyle{plain}
\newtheorem{Sa}{Theorem}[section]
\newtheorem{Def}[Sa]{Definition}
\newtheorem{Prop}[Sa]{Proposition}
\newtheorem{Kor}[Sa]{Corollary}
\newtheorem{Lem}[Sa]{Lemma}

\theoremstyle{definition}

\newtheorem{Bem}[Sa]{Remark}

\newcommand{\adjustedaccent}[1]{%
  \mathchoice{}{}
    {\mbox{\raisebox{-.5ex}[0pt][0pt]{$\scriptstyle#1$}}}
    {\mbox{\raisebox{-.35ex}[0pt][0pt]{$\scriptscriptstyle#1$}}}
}
\newcommand\smileacc[1]{\overset{\adjustedaccent{\smallsmile}}{#1}}
\newcommand\frownacc[1]{\overset{\adjustedaccent{\smallfrown}}{#1}}

\allowdisplaybreaks

\makeatletter
\g@addto@macro\bfseries{\boldmath}
\makeatother

\makeatletter
\let\blx@rerun@biber\relax
\makeatother

\hypersetup{
     pdfauthor={Xaver Kriechbaum},
   pdftitle = {BRWRE},
   colorlinks=true,
   linkcolor=blue,
    linktocpage
}

\title{Tightness for branching random walk in a\\ space-inhomogeneous random environment}
\author{Xaver Kriechbaum\footnote{Weizmann Institute of Science, Israel. \url{xaver.kriechbaum@weizmann.ac.il}}\ \orcidlink{0009-0005-4470-299X}}
\date{}

\newcommand{\barrierprobUB}{\smileacc{p}_n^{x}(y)}
\newcommand{\barrierprobLB}{\frownacc{p}_n^{y_0}(y)}

\bibliography{LitSpaceInh}
\pgfplotsset{compat=1.16}
\usepackage{microtype}
\begin{document}
\maketitle
\begin{abstract}
We consider the maximum $M_t$ of branching random walk in a space-inhomogeneous random environment on $\Z$. In this model the branching rate while at some location $x\in\Z$ is randomized in an i.i.d.\@ manner. We prove that there is a centering $\wt{m}_t$ depending only on the environment such that $(M_t-\wt{m}_t)_{t\ge 0}$ is tight in an annealed sense.
\end{abstract}
\setcounter{tocdepth}{2}
\tableofcontents
\section{Introduction}
\subsection{Model and main results}\label{Sec:Definitions and main result}
We study the maximum of a branching random walk in space-inhomogeneous random environment. We start with an informal description of the model we use. We note that the same model has been studied in \cite{CernyDrewitz} and is similar to the model from \cite{CometsPopov}. Our particles live in $\Z$ and move like a time-continuous simple random walk with jump rate 1. To each location $x\in\Z$ we associate a random branching rate $\xi(x)$. While a particle is at location $x$ it splits into two particles with rate $\xi(x)$. If the environment $\xi$ is fixed, we call the law associated with this random process the \emph{quenched} law, if we take $(\xi(x))_{x\in\Z}$ i.i.d.\@ under some measure $\mathbb{P}$ we call the law of the process averaged over $\mathbb{P}$ its \emph{annealed} law. We denote by $\mathbf{H}_n$ the first time that any of the particles hit $n\in\N$ and by $M_t$, $t\ge 0$, the maximal position of the particles alive at time $t$. In the case that $\xi(x) \equiv \xi>0$ is constant the behavior of $(M_t)_{t\ge 0}$ is well understood, which we discuss in slightly more detail in Section~\ref{Sec:Literature}. The main result of this paper is that there is a centering $(m_n)_{n\in\N}$ depending on the environment $(\xi(x))_{x\in\Z}$ such that, under the annealed law, the sequence $(\mathbf{H}_n-m_n)_{n\in\N}$ is tight. As a corollary, we establish that there is a centering $\wt{m}_t$ such that the sequence $(M_t-\wt{m}_t)_{t\ge 0}$ is tight, which answers the third open question in \cite[Section 8]{CernyDrewitz}.

Next, we recall the model studied in \cite{CernyDrewitz} in more detail. Let $(\xi(x))_{x\in\Z}$ be an i.i.d. collection of random variables on a probability space $(\Omega,\mathcal{F},\mathbb{P})$ with 
\begin{equation}\label{eq:EssinfEsssup}
0<\mathrm{ei}:= \mathrm{ess\,inf}\,\xi(0)< \mathrm{ess\,sup}\,\xi(0) =:\mathrm{es}<\infty.
\end{equation} We use $\mathbb{E}$ to denote the expected value corresponding to $\mathbb{P}$. Given a realization of $\xi$ and an initial condition $x\in\Z$, start with one particle at site $x$. All particles move independently according to a continuous-time simple random walk with jump rate 1. While at site $x$ a particle splits into two at rate $\xi(x)$ independently of everything else. These particles then evolve independently according to the same mechanism. We write $P_{x}^\xi$ and $E_{x}^\xi$ (the \emph{quenched} law and expectation respectively) for the law of the process conditioned on starting with a single particle at $x$ and environment $\xi$. The \emph{annealed} law of the process is defined as $\mathbb{P}\otimes P_{x}^\xi$, when starting with a single particle at location $x$, respectively $\mathbb{P}\otimes P^\xi$ if the starting particle is located at zero. We also use the shortened notation $P_x$ or $P$ to denote the annealed law of the process.

Let $N(t)$ denote the set of particles alive at time $t$. For $Y\in N(t)$, we denote by $(Y_s)_{s\in[0,t]}$ the path of the particle and its ancestors up to time $t$. We are interested in the first time any particle hits $n$, i.e.\@ $\mathbf{H}_n := \inf\{t\ge 0: \exists_{Y\in N(t)} Y_t = n\}$ as well as the maximal particle position $M_t := \max_{Y\in N(t)} Y_t$. 

Next, we introduce the quantities necessary to define the centering $m_n$ around which $\mathbf{H}_n$ is tight. This setup has already been used by \cite{CernyDrewitz} and is copied from their Section 4.1 without significant changes. The aim is to introduce a tilted measure under which the underlying random walk has a drift to the right such that the behavior of the extremal particle is typical for the tilted random walk.  For this, let $(X_t)_{t\ge0}$ be a continuous time SRW on $\Z$ with jump rate $1$, starting at $x$ under $P_x^\xi$. Let $H_k := \inf\{s\ge 0 : X_s = k\}$ and $\tau_k := H_k-H_{k-1}$. We write $\zeta(x) = \xi(x)-\es\le 0$. For $n\ge 1$, $A\in\sigma(X_{s\wedge H_n}, s\in [0,\infty))$ and $\eta\le 0$ define
\begin{equation} 
P_{(n)}^{\zeta,\eta}[A] := (Z_{(n)}^{\zeta,\eta})^{-1}E_0^\xi\bigg[\exp\bigg(\int_0^{H_n}(\zeta(X_s)+\eta)\;\mathrm{d}s\bigg)\cdot \mathbf{1}_{A}\bigg], \label{eq:tilt}
\end{equation}
where $ Z_{(n)}^{\zeta,\eta} :=E_0^\xi\big[\exp\big(\int_0^{H_n}(\zeta(X_s)+\eta)\;\mathrm{d}s\big)\big]$. One can extend $P_{(n)}^{\zeta,\eta}$ to a measure on $\sigma(X_s, s\ge 0)$ using Kolmogorov's extension theorem, where the necessary consistency follows from the strong Markov property of simple random walk. We denote the extension by $P^{\zeta,\eta}$. Furthermore, define associated logarithmic moment generating functions
\begin{align}
L_k^\zeta(\eta) &:= \log\bigg( E_{k-1}^\xi\bigg[\exp\bigg(\int_0^{H_k}(\zeta(X_s)+\eta)\;\mathrm{d}s\bigg)\bigg] \bigg),
\quad
\overline{L}_n^\zeta(\eta) &:= \frac{1}{n}\sum_{k=1}^nL_i^\zeta(\eta),\quad L(\eta) := \mathbb{E}[L_1^\zeta(\eta)]. \label{eq:LogMomMean}
\end{align}

With this definition one has $Z_{(n)}^{\zeta,\eta} = \exp\big(n\overline{L}_n^\zeta(\eta)\big)$. Several useful properties of these logarithmic moment generating functions can be found in \cite[Appendix~A.1]{CernyDrewitz}.

Next, we want to define the asymptotic (linear) velocity of $M_t$. In \cite[p.\@ 5]{CernyDrewitz} it is explained that there is a deterministic function $\lambda:\R\to\R$ such that for all $v\in\R$
\begin{equation}\label{eq:Lyap}
\lambda(v) = \lim\limits_{t\to\infty} \frac{1}{t}\log\big(E^\xi\big[ |\{Y\in N(t) : Y_t\ge \lfloor t\cdot v\rfloor\}|\big] \big)
\end{equation}
and that there is a $v_c\in (0,\infty)$ such that $\lambda$ is linear on $(0,v_c)$ and strictly concave on $(v_c,\infty)$. Furthermore, they explain that there is a unique $v_0\in (0,\infty)$ such that $\lambda(v_0) = 0$. The asymptotic velocity of $M_t$ is equal to $v_0$. As they did, we assume throughout the paper that
\begin{equation}
v_0>v_c. \label{eq:Vel}
\end{equation}
From \cite[Lemma 4.2]{CernyDrewitz} we see that \eqref{eq:Vel} ensures that there exists a tilting parameter $\overline{\eta} \le 0$ for which we have $\mathbb{E}[E^{\zeta,\overline{\eta}}[H_k]] = k/v_0$, i.e.\@ that we can tilt the random walk to have the same first order behavior as the maximal particle. By the same lemma $\overline{\eta}$ can be defined via
\begin{equation}
\overline{\eta} = \mathrm{argmax}_{\eta\le 0}\left( \eta/v_0-L(\eta)\right). \label{eq:TiltPara}
\end{equation}
 Since we only care about $\eta = \overline{\eta}$ in this paper, we set $P^{\zeta} := P^{\zeta,\overline{\eta}}$. Furthermore, we set $ \vartheta^\ast := \mathrm{es}-\overline{\eta}$. We are now ready to define the ``linear part'' of the centering $m_n$. For this, we define
\begin{align}
K_n^\zeta &:= -\smash[t]{\sum_{k=1}^n} L_k^\zeta(\overline{\eta}) = -n\overline{L}_n^{\zeta}(\overline{\eta}), \label{Def:Knzeta}\\
m_n &:=(\vartheta^\ast)^{-1}(K_n^\zeta-\log(p_n)), \label{Def:mn}
\end{align}
where $p_n$ is a barrier probability that will be defined in \eqref{Def:pn}.

For $Y\in N(t)$ we define $H_k(Y) := \inf\{s\le t : Y_s = k\}$. Furthermore, we recall $\mathbf{H}_n := \inf\{t\ge 0: \exists_{Y\in N(t)} Y_t = n\}$. We aim to prove
\begin{Sa}\label{Theo:Main}
The sequence
$
(\mathbf{H}_n-m_n)_{n\in\N}
$
is tight under the annealed measure $P$.
\end{Sa}

\begin{Bem}
Similarly to the calculations done in \cite[Section 6]{CernyDrewitz} we can move to a point where we need to compare certain barrier probabilities with different starting heights, which for the Gaussian case has been handled in  \cite{TimeInhKriech}. An added difficulty in the present setting is that for a fixed particle $Y$ the sequence $(H_k(Y))_{k\in\N}$ is neither Gaussian nor a random walk with i.i.d.\@ increments under the quenched measure. Thus, we approximate it by a Gaussian random walk to perform barrier computations. This is the main technical difficulty of the paper and will be handled in Section~\ref{Sec:BarrComps}.

We also note that roughly speaking $m_n \approx n/v_0+Z_n n^{1/2}+O(\log(n))$, where $Z_n$ is a function of the environment which is typically of order 1, see Remark~\ref{Rem:OrdKn} and Lemma~\ref{Lem:LBpn} for an explanation.

\end{Bem}

As a corollary of Theorem~\ref{Theo:Main} we establish tightness of $M_t$ around the centering
\[
\wt{m}_t := \inf\{k\in\N : m_{k+1}\ge t\}.
\]
We note that, by the previous remark, roughly speaking $\wt{m}_t \approx v_0\cdot t+Z_t't^{1/2}-O(\log(t))$, where $Z_t'$ again is a function of the environment that is typically of order 1. This is consistent with \cite[Theorem 2.1]{CernyDrewitz}.

\begin{Kor}\label{Kor:Max}
The sequence
$
(M_t-\wt{m}_t)_{t\ge 0}
$
is tight under the annealed measure $P$.
\end{Kor}

Two natural follow-up questions to our results are:
\begin{enumerate}[(i)]
\item Do $(\mathbf{H}_n-m_n)_{n\in\N}$ and $(M_t-\wt{m}_t)_{t\ge0}$ converge in distribution as $n\to\infty$, $t\to\infty$?
\item Does an analogous result hold in higher dimension, i.e.\@ if the particles live in $\Z^d$?
\end{enumerate}
Regarding (i) it is worth referring to the simulations pictured in \cite[Figure 1]{vcerny2022tightness}, which seem to suggest that the variance of $M_t-\wt{m}_t$ fluctuates in $t$ in a way depending on the environment. However, this does not contradict the possibility that there is some function $\varphi(\xi,t)$ such that $\varphi(\xi,t)(M_t-\wt{m}_t)$ converges in distribution and as far as we know this question is open.

Regarding (ii) we want to emphasize that our methods are only viable for the one-dimensional case. The reason for this is that we consider the times $H_k(X)$ at which a (tilted) random walk $(X_t)_{t\ge 0}$ hits $k$, and use that the sequence $(H_{k+1}(X)-H_k(X))_{k\in\N}$ is quenched independent with the distribution of $H_{k+1}(X)-H_k(X)$ a function of the environment. Meanwhile in higher dimensions if one considers the hitting times $\wt{H}_k(X)$ of hyperplanes, it is not deterministic where on such a hyperplane $X$ is at time $\wt{H}_k(X)$ and thus the distribution of $\wt{H}_{k+1}(X)-\wt{H}_k(X)$ depends on $X_{H_k(X)}$.
\subsection{Literature}\label{Sec:Literature}
In the homogeneous situation, i.e.\@ $\xi(x) \equiv \xi\in\R$ for all $x\in\Z$, the behavior of $M_t$ is well understood, the citations we give will sometimes be time-discrete instead of time-continuous. In \cite{Biggins,Hammersley,Kingman} it has been proved that $\lim_{n\to\infty} M_t/t = v\in \R$ exists. In \cite{AddarioBerryReed,BramsonZeitouni,HuShi} it has been shown that there is a centering $m_t$ such that $(M_t-m_t)_{t\ge 0}$ is tight, and that the correct centering is given by $m_t = vt-c\log(t)$ for an explicit $c\in\R$. In \cite{Aidekon,BramsonDingZeitouni} this has been extended to $(M_t-m_t)_{t\ge0}$ converging in distribution. Finally, \cite{Madaule,MalleinExtProc} have studied the extremal process, i.e.\@ the position of the particles at time $t$ shifted by the position of the maximum. They show that the extremal process converges to a decorated Poisson process.

Next, we give an overview of the literature regarding the model we use and closely related ones. Our main motivation is \cite{CernyDrewitz} in which for the same model we consider an invariance principle for the maximum has been proved. More precisely they show that under suitable conditions there is an explicit $\overline{\sigma}_{v_0}\in (0,\infty)$, such that in annealed distribution
\[
\bigg(\frac{M_{nt}-v_0nt}{\overline{\sigma}_{v_0}\sqrt{n}}\bigg)_{t\ge 0}
\]
converges to standard Brownian motion as $n\to\infty$. 

Before them \cite{CometsPopov} have considered multi-dimensional time-discrete BRW in an i.i.d.\@ environment. They allow not only the number of offspring but also the increments of the offspring to depend on space. Among other things, they established a shape theorem for the branching random walk, which implies a law of large numbers for the position of the maximum.

In \cite{vcerny2022tightness,vcerny2023bounded,drewitz2022invariance,10.1214/23-EJP956} the space-continuous version of the model from \cite{CernyDrewitz} has been studied, in which the environment is replaced by a suitably stationary and mixing random process $(\xi(x))_{x\in\R}$, and particle movement is replaced by Brownian motion. They also allow particles to have $k$ instead of 2 children, $k\in\N$.  In this situation \cite{10.1214/23-EJP956} establishes a law of large numbers and an (annealed) invariance principle for the maximum. In \cite{vcerny2022tightness} it is proved that for almost every realization of the environment $\xi$ the maximum is tight around its quenched median under the quenched law of the process. While this result is stronger than our tightness result in Theorem~\ref{Theo:Main} for the discrete situation, the proof in \cite{vcerny2022tightness} uses the randomized F-KPP equation with Heaviside initial condition $w_0^y = \mathbf{1}_{[y,\infty)}$,
\begin{equation} \label{eq:RandFKPP}
\partial_t w(t,x) = \frac{1}{2}\partial_x^2 w(t,x)+\xi(x)F(w(t,x)),\quad w(0,x) = w_0(x),\quad  t>0, x\in\R,\\
\end{equation}
where $F$ is a polynomial which is determined by the branching probabilities $(p_k)_{k\in\N}$, that a particle splits into $k$ children. Concretely they use that $w^y(t,x) = P^{\xi}_x[M_t\ge y]$ and a ``Sturmian principle'' for a related PDE.   In \cite{drewitz2022invariance} the PDE \eqref{eq:RandFKPP} is studied, they show an invariance principle for the front of \eqref{eq:RandFKPP}, i.e.\@ for the median of $w(t,\cdot)$. Additionally, they show an invariance principle for the front of the linearized version of \eqref{eq:RandFKPP}, the parabolic Anderson model.\footnote{Here the solution takes values outside of $[0,1]$ and the concept of a quantile is generalized in the appropriate manner.} They also show that the two fronts have at most distance $\log(t)$. In \cite{vcerny2023bounded} it is shown that for the parabolic Anderson model the transition front of the solution is bounded, which means that for $\varepsilon<M<\infty$, the difference between the $\varepsilon$-quantile and the $M$-quantile is a bounded function of $\varepsilon$, $M$, which does not depend on $t$. In contrast for \eqref{eq:RandFKPP} there exist environments $\xi$ such that the length of the transition front is logarithmic in $t$, at least along an increasing subsequence $t_n$, $n\in\N$.\footnote{We refer to \cite[Section 2.1]{vcerny2022tightness} for a detailed discussion why this does not contradict the tightness of $M_t$ around its quenched median. Roughly speaking it is related to the fact, that the space coordinate of the PDE corresponds to the starting point of $M_t$, so that the environment in the beginning may vary quite a lot. In contrast, for tightness of $M_t$ around its median the starting point is fixed to be 0.} 

After this paper appeared, an approach similar to \cite{vcerny2022tightness} based on studying discrete space linear PDEs has been used in \cite{vcerny2026tightness} to show 
\emph{quenched} tightness for the model considered in this paper. Additionally, they also don't need the assumption \eqref{eq:Vel}.

In \cite{10.1214/21-ECP386} it has been established, that there are deterministic increasing subsequences $(t_n)_{n\in\N}$ of $(0,\infty)$ such that $(M_{t_n}-E_0^\xi[M_{t_n}])_{n\in\N}$ is tight. The argument used is similar to the classical Dekking-Host argument from \cite{DH91}. 

Next, we cite a non-exhaustive assortment of other papers which study space- or time-inhomogeneous branching random walk models, without going into details.

In \cite{fang2012tightness} sufficient conditions for $\mathbb{P}$-a.s.\@ quenched tightness of $M_n$ around its median are given using a recursive structure of $F_{M_n}$. They allow time dependence and local dependence between siblings.

In \cite{LubetzkyThornettZeitouni,10.1214/23-EJP960} branching Brownian motion in a periodic space-inhomogeneous environment has been studied, the corresponding PDE has been studied in \cite{hamel2016logarithmic}.

In \cite{10.1214/aop/1176991677,10.2307/2244200} branching Brownian motion in a deterministic continuous space-inhomogeneous environment has been studied.

Branching random walks in random time-inhomogeneous environments have been studied in \cite{bhattacharya2021extreme,BiggKyp,HuangLiu,huang2022central,TimeInhKriech,MalleinInterfaces, MM_timeinh,zhiqiang2014central}.

Branching random walk in deterministic time-inhomogeneous environments have been studied in \cite{BovHart,FangZei,MaillZei,Mall,doi:10.1080/03605302.2014.972744,FOu}.

\subsection{Structure of the paper}
In Section~\ref{Sec:MT1} we do some preliminary work, by defining barrier events and stating the many-to-one and many-to-two formulas we use.

In Section~\ref{Sec:Structure} we give an outline of the proof of Theorem~\ref{Theo:Main}, in which we omit a lot of the technical details.

In Section~\ref{Sec:Mainbody} we prove Theorem~\ref{Theo:Main} assuming certain estimates for barrier probabilities. We note that the content of this section are mostly standard arguments in the BRW literature, see for example \cite{OZLN}. There are some minor complications because we consider the hitting times $\mathbf{H}_k$ instead of the maximum, a similar argument can be found in \cite{LubetzkyThornettZeitouni}.

In Section~\ref{Sec:BarrComps} we prove the estimates for barrier probabilities we used in Section~\ref{Sec:Mainbody}. This section is the meat of our paper and also the technically most involved part.

\subsection{Notation}\label{Sec:Notation}
Given a family of random variables $(X_{k,n})_{k,n\in\N}$ and two functions $g :\N\to \R$, $f: \N\times \N\to \R$ we write $X_{n,k}\le \wt{O}_{\text{unif},g}^n(f(k,n))$ to denote that there is a tight sequence of random variables $C_n$ on $(\Omega,\mathcal{F},\mathbb{P})$ such that $ \sup_{k\le g(n)} \frac{X_{k,n}}{f(k,n)}\le C_n$. Sometimes we will write $\wt{O}_{\text{unif}, k\le g(n)}^{n}(f(k,n))$ or $\wt{O}_{\text{unif}, k\le g}^{n}(f(k,n))$ instead. For a family of random variables $(X_n)_{n\in\N}$ we write $X_n\le \wt{O}_n(f(n))$ to denote that there is a tight sequence of random variables $C_n$ on $(\Omega,\mathcal{F},\mathbb{P})$ so that $\frac{X_n}{f(n)}\le C_n$. We also use the deterministic Landau symbols $f\le O_n(g(n))$ respectively $f\le o_n(g(n))$ to denote that $\limsup_{n\to\infty} \frac{|f(n)|}{g(n)}<\infty$ respectively $\lim_{n\to\infty} \frac{|f(n)|}{g(n)} = 0$, where the index $n$ will be omitted if it is clear from context.

We use the notation $E[X;A] := E[X\cdot\mathbf{1}_A]$.

Throughout the paper, we use $y_0\ge e+1$ to denote a fixed integer constant, whereas $c,C\ge 0$ will denote constants which may change from line to line. Similarly, $c_\varepsilon$, $C_\varepsilon$ will often denote positive constants, which depend on some $\varepsilon$ and also are allowed to change from line to line.

For $n\in\N$ we set $[n] := \{0,\dots,n\}$, for $x\in \R$ we set $J_x := [x-1,x]$. For a family $(X_s)_{s\in I}$  and $t_1,t_2\in I$ we set $X_{t_2;t_1} := X_{t_2}-X_{t_1}$.

For $\mu\in\R$, $\sigma^2\in (0,\infty)$ we use $g_{\mu,\sigma^2}(x) := \frac{1}{\sqrt{2\pi \sigma^2}}\smash{e^{-\frac{(x-\mu)^2}{2\sigma^2}}}$ to denote the Gaussian density with mean $\mu$ and variance $\sigma^2$.
We set $\mathbb{N}_{\ge x} := \N\cap[x,\infty)$, $x\in\R$.

There is an index of symbols at the end of this document.

\paragraph{Acknowledgements} This project was supported by Israel Science Foundation grant number 421/20.

Thanks to Ofer Zeitouni for suggesting the problem and for many useful discussions.

Thanks to Haotian Gu for pointing out multiple mistakes in earlier versions of this document. 
\section{Many-to-one and barrier events} \label{Sec:MT1}
We first give a rather general definition of barrier events. This is mainly useful to shorten notation.
\begin{Def}[Barrier Events] \label{Def:BarrierEvents}
For $I\subseteq \R$ closed, $t := \max I$, $(Z_s)_{s\in I}$ a real-valued process, $f:I\to \R$, $y\in\R$ and $J\subseteq \R$ an interval, define
\begin{equation}\label{eq:GenBarrierEvent}
\mathcal{B}_{I,f}^{y,J}(Z_{\cdot}) := \{\forall_{s\in I}\, y+Z_s+f(s)\ge 0,\, y+Z_t+f(t)\in J\}.
\end{equation}
If $f = 0$, $y=0$ or $J =\R$ we omit them from notation.
\end{Def}
We emphasize that our barrier events deal with the path of a random process staying \emph{above} a barrier. Next, we define the logarithmic part of $m_n$.
\begin{Def}\label{Def:LogCorr}
For $k\in\N$ set $\chi_k^2(\zeta) := \Var_{P^{\zeta}}[\tau_k]$ and $\sigma_k^2(\zeta) = \sum_{j=1}^k \chi_j^2(\zeta)$, where the dependency on the environment will be omitted most times. Under $P^{\zeta}$ let $(B_s^\zeta)_{s\ge 0}$ be a continuous Gaussian process with $(B_s^\zeta)_{s\ge0} \stackrel{d} = \big(\int_0^s \chi_{\lceil r\rceil}^2 \;\mathrm{d}B_r\big)_{s\ge0}$, where $B_r$ is a standard Brownian motion under $P^{\zeta}$. We note that under $P^{\zeta}$ we have $B_k^\zeta \sim \mathcal{N}(0,\sigma_k^2)$.

Next, we take
\begin{equation}
 W_n := (\vartheta^\ast)^{-1}K_n^\zeta-\smash[t]{\sum_{k=1}^n} (L_k^\zeta)'(\overline{\eta}). \label{Def:Wnzeta}
\end{equation}
We recall that $y_0\ge e+1$ is an integer fixed throughout the paper and set
\begin{equation}
p_n := P^{\zeta}\left[\mathcal{B}_{[0,n]}^{y_0,J_{y_0}}(B_{\cdot}^\zeta-W_{\cdot}^\zeta)\right], \label{Def:pn}
\end{equation}
where we extend $(W_k)_{k\in\N}$ to $(W_s)_{s\ge 0}$ via linear interpolation.
\end{Def}

\begin{Bem}\label{Rem:OrdKn}
First, we note that using \cite[Lemma 4.2,  eq.\@ (A.9)]{CernyDrewitz} one can see that $W_n$ corresponds to $R_n'$ from \cite[p.\@ 42]{CernyDrewitz}, we call it $W_n$ for consistency with \cite{TimeInhKriech}. In particular \cite[Lem. 6.3]{CernyDrewitz} states that $n^{-1/2}W_{nt}$ converges to a Brownian motion in $\mathbb{P}$-distribution, such that we approximately can think of $p_n$ as the probability that a Brownian motion stays above another (quenched) Brownian motion.

Next, we observe that similar to \cite[Lemma A.1]{CernyDrewitz} there are $c_{\chi}, C_{\chi}\in (0,\infty)$ such that $\mathbb{P}$-a.s.\@ for all $k\in\N$ we have $\chi_k^2 \in [c_{\chi},C_{\chi}]$ and there are $c_{L^\zeta}, C_{L^\zeta} \in (0,\infty)$ such that $\mathbb{P}$-a.s.\@ for all $k\in\N$ we have that $L_k^{\zeta}(\overline{\eta})\in [-C_{L^\zeta},-c_{L^\zeta}]$.

Finally, we comment on the size of $K_n^{\zeta}$. Since as we explained above $W_n\approx \varphi B_n$ for $\varphi>0$, $(B_k)_{k\ge0}$ a standard Brownian motion, we get that $K_n^\zeta\approx \sum_{k=1}^n (L_k^\zeta)'(\overline{\eta})+\varphi B_n$. Our choice of $\overline{\eta}$ ensures that
$
\mathbb{E}[(L_k^\zeta)'(\overline{\eta})] = \mathbb{E}[E^{\zeta}[\tau_k]] = 1/v_0,
$
such that
\[
K_n^\zeta \approx n/v_0+\sum_{k=1}^n ((L_k^\zeta)'(\overline{\eta})-\mathbb{E}[(L_k^\zeta)'(\overline{\eta})])+\varphi B_n.
\]
In \cite[Lem.~6.3]{CernyDrewitz} it is proven that the sum in the last display also converges to a Brownian motion in $\mathbb{P}$-distribution such that overall we get that
\[
K_n^\zeta \approx n/v_0+\varphi'B_n,
\]
i.e.\@ $K_n^\zeta$ is like $n/v_0$ plus fluctuations that are typically of order $n^{1/2}$. 
\end{Bem}

The quantity $p_n$ decays at most polynomially in $n$. 
\begin{Lem} \label{Lem:LBpn}
There is a $\gamma \in (0,\infty)$ such that $\mathbb{P}$-a.s.
$
\limsup_{n\to \infty} |\log(p_n)|/\log(n) \le \gamma.
$
\end{Lem}%
\begin{Bem}\label{Rem:VarStruct}
Since $p_n$ is defined as barrier probability of the Brownian motion (with time-dependent variance) $B_s^\zeta$ this can be proven similarly to \cite[Lemma 2.6]{TimeInhKriech}.\footnote{We recall that in \cite{TimeInhKriech} we considered the probability that a Brownian motion stays below a (random) barrier, while in this paper we deal with the probability that it says above a (random) barrier. This sign change is inconsequential for the bound we get but still should be kept in mind.} The proof from \cite{TimeInhKriech} cannot be used verbatim since in the setup of this paper the Brownian motion we use has a time-dependent variance. This can be fixed by changing the linear interpolation in the definition of $g_t$ in \cite[Def.~10.1, pp. 44, 45]{TimeInhKriech} to be ``adapted to (the local variance of) $B_s^{\zeta}$''. Since a similar but more involved calculation is done in Lemma~\ref{Lem:UnifCrudeLB} we don't expand on this change here and omit a full proof of Lemma~\ref{Lem:LBpn}.
\end{Bem}

We frequently need to reduce barrier probabilities involving the paths of all particles $(Y_s)_{Y\in N(t), s\le t}$ to barrier probabilities involving the path of a single particle, i.e.\@ a time-continuous simple random walk. This can be done using the next lemma, which follows from the results in \cite[Chapter~2,~Section~3.3]{maillardMT1}.  To keep the displays short we locally define a variant to the barrier events from Definition \ref{Def:BarrierEvents} via 
$
\mathcal{B}_{I,c^{\pm}}(Z_\cdot) := \{\forall_{k\in I}\, c_k^- \le Z_k\le c_k^+\}.
$
\begin{Lem}\label{Lem:MT1}
Fix $n\in\N$, $t\ge 0$ and two sequences $(c_k^-)_{k\le n}$, $(c_k^+)_{k\le n}$ with values in $[0,t]$, such that for all $k$ we have $c_k^-\le c_k^+$. We have $\mathbb{P}$-a.s.
\begin{equation}
E_0^\xi\big[|\{Y\in N(t) : \mathcal{B}_{[n], c^{\pm}}(H_\cdot(Y)) \big] = E_0^\xi\bigg[\exp\bigg(\int_0^{t} \xi(X_r)\;\mathrm{d}r\bigg);  \mathcal{B}_{[n],c^{\pm}}(H_{\cdot}(X))\bigg], \label{eq:MT1}
\end{equation}
Now, we additionally assume that $(c_k^-)_{k\le n}$ is monotonically increasing, $c_k^+ \equiv t$ and set $k_s := \min\{k\in \N : c_k^- \ge s\}$. Then we have $\mathbb{P}$-a.s.
\begin{multline}\label{eq:MT2}
E_0^\xi\big[|\{Y\in N(t) : \mathcal{B}_{[n],c^{\pm}}(H_{\cdot}(Y))\}|^2\big] =E_0^\xi\left[\exp\left(\int_0^{t} \xi(X_r)\;\mathrm{d}r\right) : \forall_{k\le n}\, c_k^-\le H_k(X)\le t\right]\\
\begin{aligned}+2\int_0^{t} &E_0^\xi\Bigg[\exp\left(\int_0^s \xi(X_r)\;\mathrm{d}r\right)\xi(X_s)\mathbf{1}_{\{\forall_{k< k_s}\, c_k^-\le H_k(X), H_{k_s}(X)\ge s\}}\\
&\times \left(E^{\xi}_{X_s}\left[\exp\left(\int_0^{t-s}\xi(X_r)\,\mathrm{d}r\right);\forall_{k\ge k_s}\,c_k^-\le H_k(X)+s\le t\right)\right)^2\Bigg]\;\mathrm{d}s
\end{aligned}
\end{multline}
\end{Lem}

Next, we formulate a way to lower bound barrier probabilities with index set $I$ by barrier probabilities for two subsets $I_1, I_2\subseteq I$, the proof is immediate via the Markov property.
\begin{Lem}\label{Lem:BarrSplit}
Let $I\subseteq\R$ be closed, $t:= \max I$, $(Z_s)_{s\in I}$ be a real-valued and translation invariant Markov process on some probability space $(\wt{\Omega},\wt{\mathcal{F}},\wt{\mathbb{P}})$, $f: I \to \R$ a function and $z_1,z_2\in \R$. Take $t_0\in I$ and set $I_1 := \{s\in I : s\le t_0\}$, $I_2 := \{s\in I : s\ge t_0\}$. We have that for $x_0\in \R$,
\begin{align}
\wt{\mathbb{P}}\big[\mathcal{B}_{I,f}^{z_1,J_{z_2}}(Z_{\cdot})\big] &\ge \wt{\mathbb{P}}\big[\mathcal{B}_{I_1,f}^{z_1-1,J_{x_0}+f(t_0)}(Z_{\cdot})\big] \inf_{x\in J_{{x_0}+1}}\wt{\mathbb{P}}\big[\mathcal{B}_{I_2,f}^{x,J_{z_2}}(Z_{\cdot}-Z_{t_0})\big]. \label{eq:BarrSplit1}
\end{align}%
\end{Lem}

\section{Strategy of the proof}\label{Sec:Structure}
We follow the classical strategy in branching random walks of doing a first and second moment calculation to get bounds on the tails of $\mathbf{H}_n-m_n$. We use the first moment calculation for an upper bound on the \emph{left} tail of $\mathbf{H}_n-m_n$ (Section~\ref{Sec:FirstMom}). The second moment calculation gives a lower bound on the left tail of $\mathbf{H}_n$ (Section~\ref{Sec:SecondMom}), which is then utilized together with the tree structure to get an upper bound on the right tail of $\mathbf{H}_n-m_n$ (Section~\ref{Sec:UseTree}). A version of this argument in the setting with constant environment can be found in \cite[Section 2.5]{OZLN} and the references therein.

We use a first moment calculation to get an upper bound on $P_0^\xi\left[\mathbf{H}_n-m_n\le -y\right]$, $y>0$. In Lemma~\ref{Lem:UBLD} the region $y\ge \log(n)^2$ is handled by applying Lemma~\ref{Lem:MT1} without further complications. Thus, we are left with getting an upper bound for $P_0^\xi[\mathbf{H}_n-m_n \in[-y-1,-y]]$ for $y\in [0,\log(n)^2]$. For this purpose, we introduce the barrier 
\[
t_{n;y}(k) := (\vartheta^\ast)^{-1}K_k^\zeta-\sigma_k^2\cdot (\vartheta^\ast \sigma_n^2)^{-1}\cdot \log(p_n)-y-1,\quad k\le n.
\] and note that $t_{n;y} := t_{n;y}(n) = m_n-y-1$. The $\sigma_k^2/(\vartheta^\ast\sigma_n^2)\cdot \log(p_n)$ should be thought of as linear interpolation between $0$ and $\log(p_n)$, which is adapted to the variance structure $(H_k)_{k\le n}$ has under $P^{\zeta}$ from \eqref{eq:tilt}. We then can apply standard arguments and Lemma~\ref{Lem:MT1} to arrive at 
\[
P_0^\xi[\mathbf{H}_n-m_n \in[-y-1,-y]]\le C e^{-\vartheta^\ast y} p_n^{-1}P^{\zeta}\big[\mathcal{B}_{[n], -t_{n;1}(\cdot)}^{[0,1]}(H_{\cdot}) \big]
\]  
All that is left is controlling the quotient in the last display. This is done in Lemma~\ref{Lem:FMomSimplBarrComp} by showing that there is a $c_\varepsilon>0$ such that
\begin{equation}
\liminf_{n\to\infty} \mathbb{P}\left[ p_n^{-1} P^{\zeta}\big[\mathcal{B}_{[n], -t_{n;0}(\cdot)}^{y,[0,1]}(H_{\cdot}(Y)) \big] \le (y+2)^{c_\varepsilon} \right] \ge 1-\varepsilon. \label{eq:BarrCompFM}
\end{equation}
See Figure~\ref{fig:FMFigH} for an illustration of an event similar to the one in the quenched probability in \eqref{eq:BarrCompFM}; the difference is that in \eqref{eq:BarrCompFM} $(H_k(Y))_{k\le n}$ needs to stay above the dashed line instead of above the thick line of Figure~\ref{fig:FMFigH}.

 At this point let us only note that there are two differences between the two factors in \eqref{eq:BarrCompFM}: We have replaced $H_k$ by a Brownian motion with the same mean and variance structure, and we have moved the starting point from $y$ to $y_0$. We explain more about how to do these barrier computations later in this section.

The next step in the proof of Theorem~\ref{Theo:Main} is getting a lower bound on $P_0^\xi[\mathbf{H}_n-m_n\le -y]$, $y>0$. For this we let $Z_n(y)$ be the number of particles for which $(H_k(Y))_{k\le n}$ stays above $t_{n;y}(k)-1$ and $H_n(Y)\in [t_{n;y}-1,t_{n;y}]$. Using Cauchy-Schwarz and $t_{n;y} = m_n-y-1$ yields 
\[
P_0^\xi[\mathbf{H}_n-m_n\le -y] \ge P_0^\xi[Z_n(y)\ge 1]\ge E_0^\xi[Z_n(y)]^2/E_0^\xi[Z_n(y)^2].
\]
Getting a lower bound on $E_0^\xi[Z_n(y)]$ is similar to the first moment calculation sketched above, only that instead of the upper bound \eqref{eq:BarrCompFM}, we need to establish that for all $\varepsilon>0$ there is a $\wt{c}_{\varepsilon}>0$ such that
\begin{equation}
\liminf_{n\to\infty} \mathbb{P}\left[ p_n^{-1}P^{\zeta}\big[\mathcal{B}_{[n], -t_{n;1}(\cdot)}^{y,[0,1]}(H_{\cdot}(Y))\big] \ge (y+2)^{-\wt{c}_\varepsilon} \right] \ge 1-\varepsilon. \label{eq:BarrCompSM}
\end{equation}
See Figure~\ref{fig:SMFigH} for an illustration of an event similar to the one in the quenched probability of \eqref{eq:BarrCompSM}; the difference is that in \eqref{eq:BarrCompSM} $(H_k(Y))_{k\le n}$ needs to stay above the dashed line instead of the thick line of Figure~\ref{fig:SMFigH}.

Next, we need an upper bound on $E_0^\xi[Z_n(y)^2]$. This upper bound can be established using \eqref{eq:MT2} and standard arguments. The details of this computation are in Section~\ref{Sec:Zny2UB}.

Leveraging the lower bound on $P_0^\xi[\mathbf{H}_n-m_n\le -y]$ into an upper bound on $P_0^\xi[\mathbf{H}_n-m_n\ge y]$ can be done as in the homogeneous situation. We note that at time $s\ll m_n$ we have exponentially in $s$ many particles at location $0$. On the one hand the behavior of the children of these particles is independent, on the other hand the lower bound on the left tail of $\mathbf{H}_n-m_n$ gives a lower bound on the probability that a child of a fixed particle at the origin at time $s$ reaches $n$ by time $m_n$. These two fact together give an upper bound on $P_0^\xi\left[\mathbf{H}_n-m_n \ge y\right]$. The details for this argument are in Section~\ref{Sec:UseTree}.

Lastly, we sketch how to establish the barrier estimates \eqref{eq:BarrCompFM}, \eqref{eq:BarrCompSM}. We emphasize that this is both the main technical work of the paper and the part which is new. For this let $\bar{H}_k := H_k-E^{\zeta}[H_k] = H_k-\sum_{i=1}^k (L_i^\zeta)'(\overline{\eta})$, where the second equality follows from \cite[Lemma~A.1]{CernyDrewitz}. We recall \eqref{Def:Knzeta} and Definition~\ref{Def:LogCorr}. We need to get bounds on 
\[
P^{\zeta}\big[ \mathcal{B}_{\{0,\dots, n\},-W_\cdot}^{y, J_1}(\bar{H}_{\cdot})\big]/P^{\zeta}\big[\mathcal{B}_{\{0,\dots,n\}, -W_\cdot}^{y_0, J_{y_0}}(B_\cdot^\zeta)\big].
\]
In \cite[Sections 8, 9]{TimeInhKriech} one can find\footnote{Strictly speaking \cite{TimeInhKriech} does not allow the variance of the Brownian motion $B_s^\zeta$ to be time-dependent. However, since (a) the local variance is piecewise constant and (b) we take care to adapt the barriers we use to the local variance of $B_s^\zeta$ as in Remark~\ref{Rem:VarStruct} this does not produce trouble. } bounds for fractions of this type, when in the numerator we replace $\bar{H}_k$ by $B_k^\zeta$. Thus, we want to couple the random walk $\bar{H}_k$ to a sum of Gaussian random variables with the same variances. A classical tool for this is the KMT-coupling for sums of independent random variables from \cite{Sakhanenko}, which roughly speaking allows us to couple $(\bar{H}_k, B_k^\zeta)_{k\le n}$ so that 
\[
P^{\zeta}\big[\sup_{k\le n} |\bar{H}_k-B_k^\zeta| \ge \log(n)^{1.5} \big] \le Ce^{-c\log(n)^{1.5}}.
\] The problem is, that if we use this on all of $[0,n]$, we introduce an error of size $\log(n)^{1.5}$ into the starting point of the Brownian motion and the interval we end up in, i.e.\@ we'd get
\begin{equation} \label{eq:BadKMT}
P^{\zeta}\left[\mathcal{B}_{[n],-W_\cdot}^{y,J_0}(\bar{H}_\cdot)\right]\lessapprox P^{\zeta}\left[\mathcal{B}_{[n], -W_\cdot}^{\log(n)^{1.5}, [-1,2\log(n)^{1.5}]}(B_\cdot^\zeta)\right].
\end{equation}
From the homogeneous situation, we expect that 
\[
P^{\zeta}\left[\mathcal{B}_{[n], -W_\cdot}^{\log(n)^{1.5}, [-1,2\log(n)^{1.5}]}(B_\cdot^\zeta)\right]/P^{\zeta}\left[\mathcal{B}_{\{0,\dots,n\}, -W_\cdot}^{y_0, J_{y_0}}(B_\cdot^\zeta)\right] \approx \log(n)^{C},
\]
for some $C>0$ so that \eqref{eq:BadKMT} is not good enough. The problem of the ending interval being blown up can be fixed by employing the classical argument (cf.~\cite[Proof of (3.4), (3.5), pp. 11--13]{MM_timeinh} for an use of the technique) of cutting $[n]$ into three equal length pieces and using the Markov-property. This reduces the problem to approximating $P^{\zeta}\big[\mathcal{B}_{\{0,\dots, n/3\},-W_\cdot}^{y,[\sqrt{n},2\sqrt{n}]}(\bar{H}_{\cdot})\big]$ by the corresponding Gaussian probability. We still plan to use KMT, however since now the interval we end up in is of length $\sqrt{n}$, it doesn't matter that the application of KMT blows it up by $\log(n)^{1.5}$. 

Thus, it remains to deal with the problem that the application of KMT shifts the starting point by $\log(n)^{1.5}$. For this, we use a dyadic induction scheme and couple the random walk on time-scale $[2^k,2^{k+1}]$ with a Brownian motion one by one, starting with $[2^{k_0},n/3]$, where $k_0$ is the maximal $k$ with $2^k\le (n/3)^{1/6}$, and then continuing backwards inductively. For this argument we set $c_k := 2^k$ and recall the notation $W_{t_2;t_1} := W_{t_2}-W_{t_1}$, $\bar{H}_{t_2;t_1} = \bar{H}_{t_2}-\bar{H}_{t_1}$. An important step, which we don't detail here, is to show that at times $c_k$, $k\le k_0$, we have $y+H_{c_k}-W_{c_k}\approx y+\sqrt{c_k}$. Assuming this, we can use the Markov property at time $c_{k-1}$ to get 
\begin{multline*}
P^{\zeta}\big[ \mathcal{B}_{\{0,\dots, c_k\},-W_\cdot}^{y,y+[\sqrt{c_k},2\sqrt{c_k}]}(\bar{H}_{\cdot})\big]\\\lessapprox P^{\zeta}\big[\mathcal{B}_{\{0,\dots, c_{k-1}\}, -W_{\cdot}}^{y,y+[\sqrt{c_{k-1}},2\sqrt{c_{k-1}}]}(\bar{H}_{\cdot})\big]\cdot \max_{z\in y+[\sqrt{c_{k-1}},2\sqrt{c_{k-1}}]} P^{\zeta}\big[\mathcal{B}_{\{c_{k-1}+1,\dots, c_k\},-W_{\cdot;c_{k-1}}}^{z, y+[\sqrt{c_k},2\sqrt{c_k}]}(\bar{H}_{\cdot;c_{k-1}})\big]
\end{multline*}
Next, we apply KMT for the second probability in the last display. Since we use KMT for a time-interval of length $c_{k-1}$ this will shift the starting point by $\log(c_{k-1})^{1.5}$ and blow up the final interval by $2\log(c_{k-1})^{1.5}$, such that we get
\begin{multline*}
P^{\zeta}\big[ \mathcal{B}_{\{0,\dots, c_k\},-W_\cdot}^{y,y+[\sqrt{c_k},2\sqrt{c_k}]}(\bar{H}_{\cdot})\big] \lessapprox  P^{\zeta}\big[\mathcal{B}_{\{0,\dots, c_{k-1}\}, -W_{\cdot}}^{y,[y+\sqrt{c_{k-1}},y+2\sqrt{c_{k-1}}]}(\bar{H}_{\cdot})\big]\\
\times\max_{z\in [y+\sqrt{c_{k-1}},y+2\sqrt{c_{k-1}}]} P^{\zeta}\left[\mathcal{B}_{\{c_{k-1}+1,\dots, c_k\},-W_{\cdot;c_{k-1}}}^{z+\log(c_k-c_{k-1})^{1.5}, [y+\sqrt{c_k},y+2\sqrt{c_k}+2\log(c_{k-1})^{1.5}]}(B_{\cdot;c_{k-1}}^\zeta)\right].
\end{multline*} For the final factor we can use Gaussian calculations to show that we can just take $z = y+\sqrt{c_{k-1}}$ instead of the maximum over $y+[\sqrt{c_{k-1}},\sqrt{c_k}]$ and ignore the blow-up of the final interval. Setting $I_k(y) := y+[\sqrt{c_k},2\sqrt{c_k}+2\log(c_{k-1})^{1.5}]$ and iterating the argument this yields
\begin{equation*}
P^{\zeta}\big[\mathcal{B}_{[n/3],-W_{\cdot}}^{y,[\sqrt{n},2\sqrt{n}]}(\bar{H}_{\cdot}) \big]\lessapprox P^{\zeta}\big[\mathcal{B}_{[c_0],-W_{\cdot}}^{y, I_0(y)}(B_{\cdot}^{\zeta})\big]\prod_{k=1}^{k_0+1} P^{\zeta}\big[\mathcal{B}_{\{c_{k-1},\dots, c_k\},-W_{\cdot;c_{k-1}}}^{y+\sqrt{c_{k-1}}+\log(c_{k-1})^{1.5}, I_k(y)}(B_{\cdot;c_{k-1}}^\zeta)\big],
\end{equation*}
where we have set $c_{k_0+1} = n$ in the last display only. Using the Markov-property at times $c_k$ we can continue on from the last display\footnote{Here, we once more assume that we know that $B_{c_k}^\zeta-W_{c_k} \approx y+\sqrt{c_k}$ on the event that $B_\cdot^\zeta$ stays above $W_{\cdot}$} to get that
\begin{equation}
P^{\zeta}\left[\mathcal{B}_{\{0,\dots, n/3\},-W_{\cdot}}^{y,[\sqrt{n},2\sqrt{n}]}(\bar{H}_{\cdot}) \right] \lessapprox P^{\zeta}\left[\mathcal{B}_{\{0,\dots, n/3\}, -W_{\cdot}+C\log(\cdot)^{1.5}}^{y, [\sqrt{n},2\sqrt{n}]}(B_{\cdot}^{\zeta})\right]. \label{eq:SmartKMT}
\end{equation}
So by using this inductive argument, we do not move the starting point by $\log(n)^{1.5}$ by applying the KMT, but instead add an extra logarithmic shift to the barrier, however from \cite[Section 9]{TimeInhKriech} we know that shifts of that magnitude do not impact the overall probability by more than a multiplicative constant so that the standard argument of cutting $\{0,\dots, n\}$ into three parts, \eqref{eq:SmartKMT} and \cite[Section 8,9]{TimeInhKriech} together yield \eqref{eq:BarrCompFM}. Equation \eqref{eq:BarrCompSM} can be established similarly. The full technical details of the inductive scheme approximating the random walk $\bar{H}_k$ by $B_k^\zeta$ can be found in Section~\ref{Sec:RWtoBMHard}.
\section{Proof of Theorem~\ref{Theo:Main} assuming certain barrier computations} \label{Sec:Mainbody}
\subsection[\texorpdfstring{Upper bound on the left tail of $\mathbf{H}_n-m_n$ -- first moment calculation}{Upper bound on the left tail -- first moment calculation}]{Upper bound on the left tail of $\mathbf{H}_n-m_n$ -- first moment calculation}\label{Sec:FirstMom}
In this section, we prove the following theorem.
\begin{Sa}\label{Theo:UBLT}
For all $\varepsilon>0$, there exists a $y_\varepsilon>0$ such that
\[
\liminf_{n\to \infty} \mathbb{P}\left[P_0^\xi[\mathbf{H}_n-m_n\le -y_\varepsilon] \le \varepsilon\right] \ge 1-\varepsilon.
\]
\end{Sa}
The proof of Theorem~\ref{Theo:UBLT} is split into dealing with $P_0^\xi[\mathbf{H}_n-m_n \le -\log(n)^2]$, which can be handled using the next lemma, and $P_0^\xi[\mathbf{H}_n-m_n \in [\log(n)^2, -y_\varepsilon]]$, which is handled in Proposition~\ref{Prop:UBLT}.
\begin{Lem}\label{Lem:UBLD}
We have that for all $n\in\N$ and $y\ge0$ $\mathbb{P}$-a.s.\@
\begin{align*}
P_0^\xi\big[\mathbf{H}_n+y\le (\vartheta^\ast)^{-1}K_n^\zeta\big]\le e^{-\vartheta^\ast y}.
\end{align*}
\end{Lem}

\begin{proof}
We recall that $\zeta(x) = \xi(x)-\mathrm{es} \le 0$ and $\vartheta^\ast = \mathrm{es}-\overline{\eta}$. Using the Markov inequality and \eqref{eq:MT1} yields that $\mathbb{P}$-a.s. 
\begin{align*}
P_0^\xi\left[\mathbf{H}_n+y\le (\vartheta^\ast)^{-1} K_n^\zeta\right] &\le E_0^\xi\left[ \exp\left(\int_0^{(\vartheta^\ast)^{-1}K_n^\zeta -y}\xi(X_s)\;\mathrm{d}s\right); H_n(X)+y\le (\vartheta^\ast)^{-1} K_n^\zeta\right]\\
&= e^{K_n^\zeta-\vartheta^\ast y} E_0^\xi\left[\exp\left(\int_0^{(\vartheta^\ast)^{-1}K_n^{\zeta}-y} \zeta(X_s)+\overline{\eta}\;\mathrm{d}s\right);H_n(X)+y\le (\vartheta^\ast)^{-1}K_n^\zeta\right]\\
&\le  e^{K_n^\zeta-\vartheta^\ast y} E_0^\xi\left[ \exp\left(\int_0^{H_n}\zeta(X_s)+\overline{\eta}\;\mathrm{d}s \right); H_n(X)+y\le (\vartheta^\ast)^{-1}K_n^\zeta\right],
\end{align*}
where the last step uses that $\zeta(x)+\overline{\eta}\le 0$ for all $x\in\Z$. Next, we use the definition of $P^{\zeta}$ and $K_n^{\zeta}$ in \eqref{eq:tilt}, \eqref{Def:Knzeta} to get from the last display
\begin{align*}
 P_0^{\xi}\left[\mathbf{H}_n+y\le (\vartheta^\ast)^{-1}K_n^\zeta\right] &\le e^{K_n^\zeta-\vartheta^\ast y} e^{n\overline{L}_n^{\zeta}(\overline{\eta})} P^{\zeta}\left[H_n(X)+y\le (\vartheta^\ast)^{-1}K_n^\zeta\right]\le e^{-\vartheta^\ast y}. \qedhere
\end{align*}
\end{proof}

\begin{Prop}\label{Prop:UBLT}
For all $\varepsilon>0$ there is a $c_\varepsilon>0$ such that
\[
\liminf_{n\to\infty} \mathbb{P}\bigg[ \bigcap_{y\in[y_0,\log(n)^2]\cap\Z} \left\{ P_0^\xi\big[\mathbf{H}_n-m_n\in[y_0-y-1,y_0-y]\big]\le e^{-\vartheta^\ast y} y^{c_\varepsilon}\right\} \bigg] \ge 1-\varepsilon.
\]
\end{Prop}

\begin{proof}[Proof of Theorem~\ref{Theo:UBLT} assuming Proposition~\ref{Prop:UBLT}]
Take $y_\varepsilon' \in\N_{\ge y_0}$ such that for $c_\varepsilon$ as in Proposition~\ref{Prop:UBLT} we have $
\sum_{y= y_\varepsilon'}^\infty  y^{c_\varepsilon}e^{-\vartheta^\ast y} \le \varepsilon.$
Then, by Proposition~\ref{Prop:UBLT}, 
\begin{equation}
\liminf_{n\to\infty} \mathbb{P}\big[ P_0^\xi[ \mathbf{H}_n-m_n \in[y_0-\log(n)^2,y_0-y_\varepsilon'] ]\le \varepsilon \big] \ge 1-\varepsilon. \label{eq:SubUBLT1}
\end{equation}
On the other hand, Lemma~\ref{Lem:UBLD} implies that $\mathbb{P}$-a.s.\@
\begin{equation}
P_0^\xi[\mathbf{H}_n-m_n\le y_0-\log(n)^2] \le e^{-\vartheta^\ast \log(n)^2+\vartheta^\ast y_0-\log(p_n)}.\label{eq:SubUBLT2}
\end{equation}
Combining \eqref{eq:SubUBLT2} with  Lemma~\ref{Lem:LBpn} implies that in $\mathbb{P}$-probability
\begin{equation}
\lim\limits_{n\to\infty} P_0^\xi[\mathbf{H}_n-m_n\le y_0-\log(n)^2] =0. \label{eq:SubUBLT3}
\end{equation}
Combining \eqref{eq:SubUBLT1} with \eqref{eq:SubUBLT3} finishes the proof.
\end{proof}

We need a few additional definitions before we can proceed with proving Proposition~\ref{Prop:UBLT}. The basic technique used in the proof of Proposition~\ref{Prop:UBLT} is a first moment argument for which we bound certain barrier probabilities. In this process as usual we add a curved banana $h_n$, which is much smaller than $(k\cdot (n-k))^{1/2}$, to the barrier.

 However, in our current setup the underlying process $H_k$ for which we estimate barrier probabilities does not have identically distributed increments, because of this we scale the increments of the banana by the local variance, see \eqref{eq:IncrPartBan}. We also will take the middle part of the banana between $n/3$ and $2n/3$ to be an interpolation adapted to the variance of $B_s^\zeta$ in the sense of \eqref{Rem:VarStruct}. This is done to ensure that the banana starts and ends at 0.  Recall the definitions of $\chi_k^2$, $\sigma_k^2$, $\sigma_n^2$ from Definition~\ref{Def:LogCorr} as well as the definition of $\vartheta^\ast$ immediately before \eqref{Def:Knzeta}.
 \begin{figure}
 \centering
 \begin{subfigure}[t]{0.45\linewidth}
 \includegraphics[scale = 0.57]{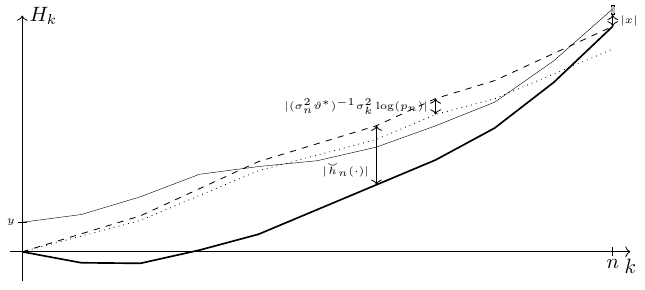}
 \caption{Visualization of the barrier event in \protect{\eqref{Def:bprobUB}} in terms of $(y+H_k)_{k\le n}$. The dashed line is $((\vartheta^\ast)^{-1}K_k^\zeta)_{k\le n}$, the dotted line is $((\vartheta^\ast)^{-1}K_k^\zeta+(\vartheta^\ast \sigma_n^2)^{-1}\sigma_k^2 \log(p_n))_{k\le n}$ and the thick line is $((\vartheta^\ast)^{-1}K_k^\zeta + \smileacc{m}_{n,\smileacc{h}}(k))_{k\le n}$. Finally, the thin black line is an interpolation of a sample path of $(y+H_k)_{k\le n}$ realizing the event in \protect{\eqref{Def:bprobUB}}.}
 \label{fig:FMFigH}
 \end{subfigure}\hskip1cm
 \begin{subfigure}[t]{0.45\linewidth}
 \includegraphics[scale = 0.57]{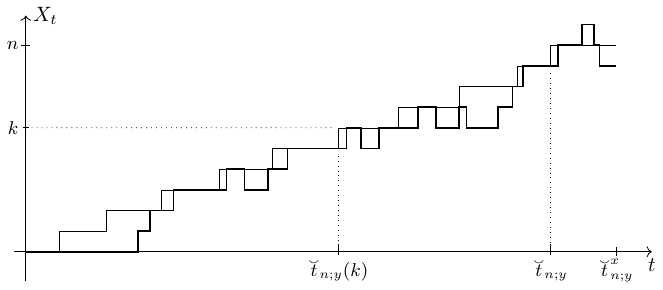}
 \caption{Visualization of the event in \protect{\eqref{Def:bprobUB}} in terms of of $(X_t)_{t\le \smileacc{t}_{n;y}^x}$. The thin black line is the barrier and the thick black line is a sample of $(X_t)_{t\le \smileacc{t}_{n;y}^x}$ realizing the event in \protect{\eqref{Def:bprobUB}}.}
 \label{fig:FMFigX}
 \end{subfigure}
\caption{Visualization of some quantities in Definition~\protect{\ref{Def:HelpstuffFirstMoment}}.}
\label{fig:FMFig}
 \end{figure}
\begin{Def} \label{Def:HelpstuffFirstMoment}
For $n\in\N$, $k\le n$ set $\smileacc{g}_n(k) := -\delta((1+k)^{1/6}\wedge(1+n-k)^{1/6}-1)$, where $\delta>0$ will be chosen later. Furthermore, we set $\Delta_{k,k-1}(\smileacc{g}_n) := \smileacc{g}_n(k)-\smileacc{g}_n(k-1)$. Now, we define $\smileacc{h}_n(k)$, $k\in [n]$, recursively. Set $\smileacc{h}_n(0) = \smileacc{h}_n(n) = 0$. For $k< \lfloor n/3\rfloor$ we define $\smileacc{h}_n(k+1)$ as a function of $\smileacc{h}_n(k)$ via
\begin{equation} \label{eq:IncrPartBan}
\smileacc{h}_n(k+1) = \smileacc{h}_n(k)+\chi_{k+1}^2\Delta_{k+1,k}(\smileacc{g}_n).
\end{equation}
For $k>n-\lfloor n/3\rfloor$ we define $\smileacc{h}_n(k-1)$ as a function of $\smileacc{h}_n(k)$ via
\begin{equation}\label{eq:IncrPartEnd}
\smileacc{h}_n(k-1) = \smileacc{h}_n(k)-\chi_k^2\Delta_{k,k-1}(\smileacc{g}_n).
\end{equation}

 For $k\in \{\lfloor n/3\rfloor+1,\dots, n-\lfloor n/3\rfloor-1\}$ we set $\smileacc{h}_n$ to be the (weighted) linear interpolation of $\smileacc{h}_n(\lfloor n/3\rfloor)$ and $\smileacc{h}_n(n-\lfloor n/3\rfloor)$, i.e.\@ we define
\begin{equation}\label{eq:IncrPartMid}
\smileacc{h}_n(k) = \smileacc{h}_n(\lfloor n/3\rfloor)+(\smileacc{h}_n(n-\lfloor n/3\rfloor)-\smileacc{h}_n(\lfloor n/3\rfloor))\frac{\sum_{j= \lfloor n/3\rfloor+1}^k \chi_j^2}{\sum_{j=\lfloor n/3\rfloor+1}^{n-\lfloor n/3\rfloor} \chi_j^2}.
\end{equation}

Furthermore, for $n\in\N$, $y,x\in\R$, $k\in[n]$ we set
\begin{align}
\smileacc{m}_{n,\smileacc{h}}(k) &:= \smileacc{h}_n(k)-\frac{\sigma_k^2}{\vartheta^\ast \sigma_n^2}\log(p_n),\quad  \barrierprobUB := P^{\zeta}\left[ \mathcal{B}_{[n],-\smileacc{m}_{n,\smileacc{h}}}^{y, J_x,}(H_\cdot-(\vartheta^\ast)^{-1}K_{\cdot}^\zeta ) \right] \label{Def:bprobUB}\\
\smileacc{t}_{n;y}^x(k) &:= (\vartheta^\ast)^{-1}K_k^\zeta+x-y+\smileacc{m}_{n,\smileacc{h}}(k),\ \smileacc{t}_{n;y}(k) := \smileacc{t}_{n;y}^0(k),\ \smileacc{t}_{n;y}^x := \smileacc{t}_{n;y}^x(n). \label{Def:Deftnyx}
\end{align}
See Figure~\ref{fig:FMFig} for an illustration\footnote{In the figures we do not account for the fact that the increments of $\smileacc{h}_n$ are scaled by the local variance of $H_\cdot$.} of the event in \eqref{Def:bprobUB}. We now fix $\delta>0$ so that $\smileacc{t}_{n;y}^x(\cdot)$ is monotonically increasing and note that this is possible since $(K_{k+1}^{\zeta}-K_k^{\zeta})_{k\in\N}$ is uniformly lower bounded as we have explained above \eqref{Def:Wnzeta}. Furthermore, we set
\begin{equation}
\mathbf{\tau}_n(y) := \inf\{ k\in [n] : \exists_{Y\in N(\smileacc{t}_{n;y}(k))} H_k(Y)-(\vartheta^\ast)^{-1}K_k^\zeta-\smileacc{m}_{n,\smileacc{h}}(k)+y+1\le 0\}\wedge n \label{eq:tauny}
\end{equation}
and
\begin{align} \label{eq:Anyx}
A_n(y,x)&:= \smash{\bigcup_{Y\in N(\smileacc{t}_{n;y}^x)} \mathcal{B}_{[n],-\smileacc{m}_{n,\smileacc{h}}(k)}^{y,J_x}(H_\cdot(Y)-(\vartheta^\ast)^{-1} K_\cdot^\zeta).}
\end{align}
Finally, define
\begin{equation} \label{eq:qkend}
q_{k,\text{end}}(x) := P^{\zeta}\left[\mathcal{B}_{\{n-k-1,\dots, n\}, -(\smileacc{m}_{n,\smileacc{h}}(\cdot)-\smileacc{m}_{n,\smileacc{h}}(n-k-1))}^{x, J_{y_0}}(H_{\cdot;n-k-1}-(\vartheta^\ast)^{-1}K_{\cdot;n-k-1}^\zeta) \right].
\end{equation}
In words, $q_{k,\text{end}}(x)$ is the probability that with $H_{n-k-1}$ being  $x$ above $\smileacc{t}_{n;y}(n-k-1)$, the sequence $(H_{n-k},\dots, H_n)$ stays above the barrier $(\smileacc{t}_{n;y}(n-k),\dots, \smileacc{t}_{n;y}(n))$ and $H_n$ is in $[t_{n;y}^{y_0}-1,t_{n;y}^{y_0}]$. It roughly corresponds to restricting the $x$-axis in Figure~\ref{fig:FMFigH} to $[n-k-1,n]$ and fixing the starting height at time $n-k-1$ to be~$x$. 
\end{Def}

The following two barrier estimates are pivotal for the proof of Proposition~\ref{Prop:UBLT} and are proved in Section~\ref{Sec:BarrComps}.
\begin{Lem}\label{Lem:FMomSimplBarrComp}
For all $\varepsilon>0$, there is a $c_\varepsilon>0$ such that 
\[
\liminf_{n\to \infty} \mathbb{P}\bigg[\bigcap_{\genfrac{}{}{0pt}{}{x \in[y_0,\log(n)]\cap\N}{y\in[y_0,\log(n)^2]\cap\N}} \left\{ p_n^{-1}\barrierprobUB \le y^{c_\varepsilon} x^{c_\varepsilon}\right\} \bigg] \ge 1-\varepsilon.
\]
\end{Lem}
\begin{Lem}\label{Lem:FMomHardBarrComp}
For all $\varepsilon>0$, $C>0$ there is a $C_\varepsilon>0$ such that
\[
\liminf_{n\to\infty} \mathbb{P}\left[\sum_{k=1}^{\lfloor \log(n)^7\rfloor} e^{\vartheta^\ast \smileacc{h}_n(n-k)} \left(\inf_{x\in [1,C+2]} q_{k,\text{end}}(x)\right)^{-1}\le C_\varepsilon\right]\ge 1-\varepsilon.
\]
\end{Lem}

In further preparation towards the proof of Proposition~\ref{Prop:UBLT} note that
\begin{equation} \label{eq:FMDec}
\begin{aligned}
P_0^\xi\left[ \mathbf{H}_n-m_n \in [y_0-y-1,y_0-y]\right] &= P_0^\xi\left[ \mathbf{H}_n-m_n \in [y_0-y-1,y_0-y], \tau_n(y)=n\right] \\
&\quad+P_0^\xi\left[ \mathbf{H}_n-m_n \in [y_0-y-1,y_0-y], \tau_n(y) <n\right]\\
&\le P_0^\xi\left[A_n(y,y_0)\right]+P_0^\xi\left[\tau_n(y)<n\right]. 
\end{aligned}
\end{equation}
Here the first summand corresponds to all particles hitting $k$ later than at time $\smileacc{t}_{n;y}^{-1}(k)$, while at least one particle hits $n$ between the times $m_n+y_0-y-1$ and $m_n+y_0-y$; this event can be controlled by applying  Lemmata~\ref{Lem:MT1} and~\ref{Lem:FMomSimplBarrComp}, which we do in the following Lemma~\ref{Lem:Anyxi}. The second summand in \eqref{eq:FMDec} corresponds to a particle hitting some location $k$ before time $\smileacc{t}_{n;y}^{-1}(k)$; we prove an upper bound for this probability in Lemmata~\ref{Lem:taueasy},~\ref{Lem:tauhard} below. We note that \eqref{eq:FMDec} is the place in which we have introduced a barrier, since its right-hand side is dominated by $\tau_n(y) \approx n$.

\begin{Lem}\label{Lem:Anyxi}
For all $\varepsilon>0$, there exists a $c_\varepsilon>0$ such that
\[
\liminf_{n\to\infty} \mathbb{P}\bigg[\bigcap_{y\in[y_0,\log(n)^2]\cap \Z} \left\{P_0^\xi[A_n(y,y_0)]\le e^{-\vartheta^\ast y} y^{c_\varepsilon} \right\} \bigg] \ge 1-\varepsilon.
\]
\end{Lem}
\begin{proof}[Proof of Lemma~\ref{Lem:Anyxi} assuming Lemma~\ref{Lem:FMomSimplBarrComp}]
We recall that $\vartheta^\ast = \mathrm{es}-\overline{\eta}$, $\zeta(x) = \xi(x)-\mathrm{es}$. Using the Markov inequality and applying Lemma~\ref{Lem:MT1} as well as $\zeta(X_s)+\overline{\eta}\le 0$ yields that
\begin{align*}
P_0^\xi\left[A_n(y,y_0) \right] & \le  e^{\vartheta^\ast \smileacc{t}_{n;y}^{y_0}}E_0^\xi\left[\exp\left(\int_0^{\smileacc{t}_{n;y}^{y_0}} \zeta(X_s)+\overline{\eta}\;\mathrm{d}s\right); \mathcal{B}_{[n],-\smileacc{m}_{n,\smileacc{h}}(k)}^{y,J_{y_0}}\left(H_\cdot - (\vartheta^\ast)^{-1} K_\cdot^\zeta\right)\right]\\
&\le e^{\vartheta^\ast \smileacc{t}_{n;y}^{y_0}}E_0^\xi\left[\exp\left(\int_0^{H_n} \zeta(X_s)+\overline{\eta}\;\mathrm{d}s\right);  \mathcal{B}_{[n],-\smileacc{m}_{n,\smileacc{h}}(k)}^{y,J_{y_0}}\left(H_\cdot - (\vartheta^\ast)^{-1} K_\cdot^\zeta\right) \right]\\
&\stackrel{\eqref{eq:tilt}, \eqref{Def:bprobUB}}{=} e^{\vartheta^\ast \smileacc{t}_{n;y}^{y_0}+n\overline{L}_n^\zeta(\overline{\eta})} \smileacc{p}_n^{y_0}(y) \stackrel{\eqref{Def:Deftnyx}}{=} e^{\vartheta^\ast x-\vartheta^\ast y}\frac{\smileacc{p}_n^{y_0}(y)}{p_n}.
\end{align*}
The claim now follows by applying Lemma~\ref{Lem:FMomSimplBarrComp}.
\end{proof}

The two forthcoming lemmata allow us to control the second summand in \eqref{eq:FMDec}. The cases $\tau_n(y)\not \approx n$ and $\tau_n(y)\approx n$ are handled separately. 
\begin{Lem}\label{Lem:taueasy}
There is a $C>0$ such that
\[
\liminf_{n\to \infty} \mathbb{P}\bigg[ \bigcap_{y\in [y_0,\log(n)^2]\cap\N} \big\{P_0^\xi[\tau_n(y)\le n-\log(n)^7 ] \le C e^{-\vartheta^\ast y} \big\} \bigg] =1.
\]
\end{Lem}

\begin{Lem}\label{Lem:tauhard}
For all $\varepsilon>0$ there is a $c_\varepsilon>0$ such that 
\[
\liminf_{n\to \infty} \mathbb{P}\bigg[ \bigcap_{y\in [y_0,\log(n)^2]\cap\N} \big\{P_0^\xi\big[ n-\log(n)^7< \tau_n(y)<n \big] \le e^{-\vartheta^\ast y} y^{c_\varepsilon} \big\} \bigg] \ge 1-\varepsilon.
\]
\end{Lem}

\begin{proof}[Proof of Lemma~\ref{Lem:taueasy}]
Similar to the computation in the proof of Lemma~\ref{Lem:Anyxi} we get that for all $y>0$ and $k\le n-\log(n)^7$ we have that $\mathbb{P}$-a.s.
\begin{equation}\label{eq:taunEasyMain}
\begin{aligned}
P_0^\xi\left[\tau_n(y) = k \right] &\stackrel{L.\ref{Lem:MT1}}{\le} E_0^\xi\left[ \exp\left(\int_0^{\smileacc{t}_{n;y}^{-1}(k)} \xi(X_s)\;\mathrm{d}s\right);  H_k-(\vartheta^\ast)^{-1}K_k^\zeta-\smileacc{m}_{n,\smileacc{h}}(k)+y+1\le 0  \right]\\
&\le e^{-\vartheta^\ast y} p_n^{-\frac{\sigma_k^2}{\sigma_n^2}}e^{\vartheta^\ast \smileacc{h}_n(k)} \le e^{-\vartheta^\ast y} p_n^{-C\frac{k}{n}}e^{\vartheta^\ast \smileacc{h}_n(k)},
\end{aligned}
\end{equation}
where the last step uses that by \cite[Lemma~A.1]{CernyDrewitz} we have $\chi_k^2 \in [c_\chi,C_\chi]$ for some deterministic constants $c_\chi, C_\chi \in (0,\infty)$.  By Lemma~\ref{Lem:LBpn} we know that there is a $\gamma\in [0,\infty)$ such that $\mathbb{P}$-a.s. 
\begin{equation}\label{eq:pnPoly}
\limsup_{n\to\infty} \frac{|\log(p_n)|}{\log(n)}\le \gamma.
\end{equation} We combine \eqref{eq:taunEasyMain}, \eqref{eq:pnPoly} and $\sup_{n\in\N}\sum_{k=1}^{\lfloor n-\log(n)^7\rfloor} c n^{C\frac{\gamma k}{n}} e^{\vartheta^\ast \smileacc{h}_n(k)} <\infty$ to finish the proof.
\end{proof}

\begin{proof}[Proof of Lemma~\ref{Lem:tauhard} assuming Lemmata~\ref{Lem:FMomSimplBarrComp},~\ref{Lem:FMomHardBarrComp}]
Repeating the calculation in the proof of Lemma~\ref{Lem:taueasy} yields that for $k\in [0,\dots,\lfloor \log(n)^7\rfloor]$
\begin{align}
&P_0^\xi\left[\tau_n(y) = n-k \right]  \label{eq:taunkhardzwisch}\\
&\le e^{-\vartheta^\ast y+\vartheta^\ast \smileacc{h}_n(n-k)} p_n^{-\frac{\sigma_{n-k}^2}{\sigma_n^2}}\notag\\
&\quad\times\underbrace{P^{\zeta}\left[\forall_{j\le n-k-1}\, H_j-(\vartheta^\ast)^{-1}K_j^{\zeta}-\smileacc{m}_{n,\smileacc{h}}(j)+y+1\ge 0, H_{n-k}-(\vartheta^\ast)^{-1}K_{n-k}^\zeta-\smileacc{m}_{n,\smileacc{h}}(n-k)+y+1\le 0\right]}_{=: q_k(y)}\notag
\end{align}

We want to compare $q_k(y)$ to  $p_n^{-\sigma_{n-k}^2/\sigma_n^2}\le p_n^{-1}$ by using Lemma~\ref{Lem:FMomSimplBarrComp}.  To do this we partition $q_k(y)$ with regard to the value $H_{n-k-1}$ takes. One peculiarity of the model we use is, that $(H_k)_{k\in\N}$ is monotonically increasing. Furthermore, by \cite[Lemma A.1]{CernyDrewitz} $|K_{k+1}^\zeta-K_k^\zeta| \le C <\infty$  is $\mathbb{P}$-a.s.\@ bounded. Finally, we have the $\mathbb{P}$-a.s.\@ uniform bound $|\smileacc{m}_{n,\smileacc{h}}(n-k)-\smileacc{m}_{n,\smileacc{h}}(n-k-1)|\le C<\infty$. where we use Lemma~\ref{Lem:LBpn} to control the $\log(p_n)$ part of $\smileacc{m}_{n,\smileacc{h}}$. Thus we get that for some $C_1>0$ independent of $n$ and $y$%
\[
q_k(y)\le C_1\cdot \max_{l\in [0,C_1]} P^{\zeta}\left[\mathcal{B}_{[n-k-1],-\smileacc{m}_{n,\smileacc{h}}}^{y+1, J_{l+1}}(H_{\cdot}-(\vartheta^\ast)^{-1} K_{\cdot}^\zeta) \right].
\] %
Using Lemma~\ref{Lem:BarrSplit} we can bound 
\begin{align*}
&P^{\zeta}\left[\mathcal{B}_{[n-k-1],-\smileacc{m}_{n,\smileacc{h}}}^{y-1, J_{l+1}}(H_{\cdot}-(\vartheta^\ast)^{-1} K_{\cdot}^\zeta) \right] \le  \smileacc{p}_n^{y_0}(y)\\
&\quad\quad\quad\times\left(\inf_{x \in J_{l+2}}P^{\zeta}\left[\mathcal{B}_{\{n-k-1,\dots, n\}, -(\smileacc{m}_{n,\smileacc{h}}(\cdot)-\smileacc{m}_{n,\smileacc{h}}(n-k-1))}^{x, J_{y_0}}(H_{\cdot}-H_{n-k-1}-(\vartheta^\ast)^{-1}(K_\cdot^\zeta-K_{n-k-1}) \right] \right)^{-1}\\
&\le \smileacc{p}_n^{y_0}(y)\cdot  \left(\inf_{x\in [1,C_1+2]} q_{k,\text{end}}(x)\right)^{-1}
\end{align*}
Plugging this into \eqref{eq:taunkhardzwisch} yields 
\[
P_0^\xi[\tau_n(y) = n-k] \le e^{-\vartheta^\ast y+\smileacc{h}_n(n-k)}p_n^{-1} \smileacc{p}_n^{y_0}(y)\left(\inf_{x\in [1,C_1+2]} q_{k,\text{end}}(x)\right)^{-1}.
\]
Combining the last display with Lemmata~\ref{Lem:FMomSimplBarrComp},~\ref{Lem:FMomHardBarrComp} implies the statement of Lemma~\ref{Lem:tauhard}.
\end{proof}
We are finally ready to prove Proposition~\ref{Prop:UBLT}.
\begin{proof}[Proof of Proposition~\ref{Prop:UBLT}]
The proposition follows from combining \eqref{eq:FMDec} with Lemmatas~\ref{Lem:Anyxi}--\ref{Lem:tauhard}.
\end{proof}

\subsection[\texorpdfstring{Lower bound on the left tail of $\mathbf{H}_n-m_n$ -- second moment calculation}{Lower bound on the left tail -- second moment calculation}]{Lower bound on the left tail of $\mathbf{H}_n-m_n$ -- second moment calculation}\label{Sec:SecondMom}
In this section we prove the following theorem.
\begin{Sa}\label{Theo:LBLT}
For all $\varepsilon>0$ there is a $c_\varepsilon>0$ such that for all $y\ge y_0$
\[
\liminf_{n\to\infty} \mathbb{P}\big[ P_0^\xi[\mathbf{H}_n-m_n \le y_0-y ]\ge e^{-\vartheta^\ast y} y^{-c_\varepsilon} \big]\ge 1-\varepsilon.
\]
\end{Sa}

We prove Theorem~\ref{Theo:LBLT} using a second moment approach. For this, we need the following definition, which is mostly analogous to Definition~\ref{Def:HelpstuffFirstMoment} in the previous section, i.e.\@ it sets up the banana we need to add to the barrier in order to carry out the second moment calculation, and defines the relevant barrier and barrier probability. We recall the definitions of $\chi_k^2$, $\sigma_k^2$, $K_k^\zeta$ , $p_n$ and $\mathcal{B}$ from Section~\ref{Sec:Definitions and main result} and Definition~\ref{Def:BarrierEvents}.

\begin{figure}
 \centering
 \includegraphics[scale = 0.8]{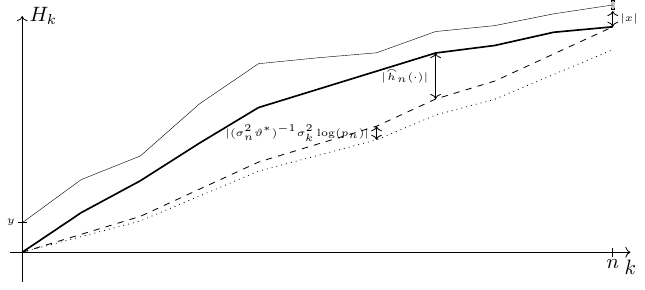}
 \caption{Visualization of the barrier event in \protect{\eqref{Def:bprobLB}} in terms of $(H_k)_{k\le n}$. The dashed line is $((\vartheta^\ast)^{-1}K_k^\zeta)_{k\le n}$, the dotted line is $((\vartheta^\ast)^{-1}K_k^\zeta+(\vartheta^\ast \sigma_n^2)^{-1}\sigma_k^2 \log(p_n))_{k\le n}$ and the thick line is $((\vartheta^\ast)^{-1}K_k^\zeta + \frownacc{m}_{n,\frownacc{h}}(k))_{k\le n}$. Finally, the thin black line is the linear interpolation of a sample path of $(H_k)_{k\le n}$ realizing the event in \protect{\eqref{Def:bprobLB}}.}
 \label{fig:SMFigH}
 \end{figure}

\begin{Def} \label{Def:HelpstuffSecMoment}
For $n\in\N$, $k\le n$ set $\frownacc{g}_n(k) := \delta((1+k)^{1/6}\wedge (1+n-k)^{1/6}-1)$, where $\delta>0$ will be chosen later on. Furthermore, set $\Delta_{k,k-1}(\frownacc{g}_n) := \frownacc{g}_n(k)-\frownacc{g}_n(k-1)$. Now define $\frownacc{h}_n(k)$, $k\in [n]$, recursively. Set $\frownacc{h}_n(0) = \frownacc{h}_n(n) = 0$. For $k< \lfloor n/3\rfloor$ we define $\frownacc{h}_n(k+1)$ as a function of $\frownacc{h}_n(k)$ via
\begin{equation}\label{eq:SecMomhn1}
\frownacc{h}_n(k+1) = \frownacc{h}_n(k)+\chi_{k+1}^2\Delta_{k+1,k}(\frownacc{g}_n).
\end{equation}
For $k>n-\lfloor n/3\rfloor$ we define $\frownacc{h}_n(k-1)$ as a function of $\frownacc{h}_n(k)$ via 
\begin{equation}\label{eq:SecMomhn2}
\frownacc{h}_n(k-1) = \frownacc{h}_n(k)-\chi_{k}^2\Delta_{k,k-1}(\frownacc{g}_n).
\end{equation}

For $k\in\{\lfloor n/3\rfloor+1,\dots, n-\lfloor n/3\rfloor-1\}$ we let $\frownacc{h}_n$ be the (weighted) linear interpolation of $\frownacc{h}_n(\lfloor n/3\rfloor)$, $\frownacc{h}_n(n-\lfloor n/3\rfloor)$, i.e. we set\@
\begin{equation}\label{eq:SecMomhn3}
\frownacc{h}_n(k) = \frownacc{h}_n(\lfloor n/3\rfloor)+(\frownacc{h}_n(n-\lfloor n/3\rfloor)-\frownacc{h}_n(\lfloor n/3\rfloor))\frac{\sum_{j=\lfloor n/3\rfloor+1}^k \chi_j^2}{\sigma_{n-\lfloor n/3\rfloor}^2-\sigma_{\lfloor n/3\rfloor}^2}.
\end{equation}

Furthermore, for $n\in\N$, $y,x\in\R$ we set
\begin{align}
\frownacc{m}_{n,\frownacc{h}}(k)&:= \frownacc{h}_n(k)-\frac{\sigma_k^2}{\vartheta^\ast \sigma_n^2}\log(p_n),\hskip3pt  \barrierprobLB := P^{\zeta}\big[\mathcal{B}_{[n],-\frownacc{m}_{n,\frownacc{h}}(k)}^{y,J_{y_0}}(H_\cdot - (\vartheta^\ast)^{-1} K_\cdot^\zeta)\big], \label{Def:bprobLB}\\
\frownacc{t}_{n;y}^x(k) &:= (\vartheta^\ast)^{-1}K_k^\zeta+x-y+\frownacc{m}_{n,\frownacc{h}}(k),\ \frownacc{t}_{n;y}(k) := \frownacc{t}_{n;y}^0(k),\ \frownacc{t}_{n;y}^x := \frownacc{t}_{n;y}^x(n), \label{Def:tnyxk}\\
Z_n(y) &:= \sum_{Y\in N(\frownacc{t}_{n;y}^{y_0})}\mathbf{1}_{\mathcal{B}_{[n],-\frownacc{m}_{n,\frownacc{h}}}^{y,J_{y_0}}(H_\cdot(Y)-(\vartheta^\ast)^{-1}K_\cdot^\zeta)}. \label{Def:Zny}
\end{align}
See Figure~\ref{fig:SMFigH} for an illustration\footnote{In the figure we do not account for the fact that the increments of $\frownacc{h}_n$ are scaled by the local variance of $H_{\cdot}$.} of the event in \eqref{Def:bprobLB}. We fix $\delta$ so that for all $x$ the sequence $\frownacc{t}_{n;y}^x(k)$ is increasing in $k$, which is possible since $(K_{k+1}^{\zeta}-K_k^{\zeta})_{k\in\N}$ is uniformly lower bounded by the bounds discussed in Remark~\ref{Rem:OrdKn}.  Finally, we set
\begin{equation}\label{eq:qkstart}
q_{k,\text{start}} := \inf_{x\in [1,2]} P^{\zeta}\big[\mathcal{B}_{\{0,\dots,k\},-\frownacc{m}_{n,\frownacc{h}}}^{y_0-1, J_x}(H_\cdot-(\vartheta^\ast)^{-1}K_\cdot^\zeta) \big].
\end{equation}
In words, this corresponds to the probability of $H_l$ staying above the barrier $\frownacc{t}_{n;y_0}(l)$ for $l\le k$ and $H_k$ being in $[\frownacc{t}_{n;y}^x(k)-1,\frownacc{t}_{n;y}^x(k)]$.
\end{Def}

Using Cauchy-Schwarz and the notations introduced in Definition~\ref{Def:HelpstuffSecMoment} yields that
\begin{equation}
P_0^\xi[\mathbf{H}_n-m_n \le y_0-y ] \ge P_0^\xi[Z_n(y)\ge1] \ge E_0^\xi[Z_n(y)]^2/E_0^\xi[Z_n(y)^2]. \label{eq:CSchwartz}
\end{equation}
Thus, we have proved Theorem~\ref{Theo:LBLT} if we can establish the following two results.
\begin{Prop}\label{Prop:EZn}
For all $\varepsilon>0$, there exists $c_\varepsilon>0$ such that for all $y\ge y_0$
\[
\liminf_{n\to \infty} \mathbb{P}\big[E_0^\xi[Z_n(y)]\ge e^{-\vartheta^\ast y} y^{-c_\varepsilon} \big] \ge 1-\varepsilon.
\]
\end{Prop}
\begin{Prop}\label{Prop:EZn2}
For all $\varepsilon>0$, there exists $c_\varepsilon>0$ such that for all $y\ge y_0$
\[
\liminf_{n\to \infty} \mathbb{P}\big[E_0^\xi[Z_n(y)^2]\le  e^{-\vartheta^\ast y} y^{c_\varepsilon}+E_0^\xi[Z_n(y)]\big]\ge 1-\varepsilon. 
\]
\end{Prop}
We emphasize that in the above two propositions and Lemmata~\ref{Lem:SecMomBarrCompLB},~\ref{Lem:SecMomBarrCompUB} the $c_{\varepsilon}$ are not equal, when applying several of these statements at once, we sometimes use $c_\varepsilon$ and $c_{\varepsilon}'$ to make the distinction clearer.
\begin{proof}[Proof of Theorem~\ref{Theo:LBLT} assuming Propositions~\ref{Prop:EZn},~\ref{Prop:EZn2}]
Fix $\varepsilon>0$. Combining \eqref{eq:CSchwartz} with Proposition~\ref{Prop:EZn2} yields that there exists $c_\varepsilon>0$ such that for $y\ge y_0$.
\[
1-\varepsilon/2\le\liminf_{n\to\infty}\mathbb{P}\bigg[P^{\zeta}[\mathbf{H}_n-m_n\le y_0-y] \ge \frac{E_0^\xi[Z_n(y)]}{1+e^{-\vartheta^\ast y}y^{c_\varepsilon}/E_0^\xi[Z_n(y)]}\bigg].
\]
By Proposition~\ref{Prop:EZn} there exists $c_{\varepsilon}'>0$ such that 
\begin{align*}
1-\varepsilon \le \liminf_{n\to\infty} \mathbb{P}\bigg[P^{\zeta}\left[\mathbf{H}_n-m_n\le y_0-y\right] \ge e^{-\vartheta^\ast y} y^{-c_{\varepsilon'}}\frac{1}{1+e^{-2\vartheta^\ast y} y^{c_\varepsilon+c_\varepsilon'}} \bigg].
\end{align*}
The fact that $\inf_{y\ge y_0} \frac{1}{1+e^{-2\vartheta^\ast y} y^{c_\varepsilon+c_\varepsilon'}} >0$ finishes the proof.
\end{proof}

To prove Propositions~\ref{Prop:EZn} and~\ref{Prop:EZn2} we will need several barrier computations, which are stated and assumed in this section but will be proved later in Section~\ref{Sec:BarrComps}.

\begin{Lem}\label{Lem:SecMomBarrCompLB}
For all $\varepsilon>0$ there exists a $c_\varepsilon>0$ such that for all $y\ge y_0$
\[
\liminf_{n\to \infty} \mathbb{P}\big[ p_n^{-1} \barrierprobLB\ge y^{-c_\varepsilon} \big] \ge 1-\varepsilon.
\]
\end{Lem}

\begin{Lem}\label{Lem:SecMomBarrCompUB}
For all $\varepsilon>0$ there exists a $c_\varepsilon>0$ such that
\[
\liminf_{n\to \infty} \mathbb{P}\bigg[\bigcap_{y\in[y_0,\log(n)^2]\cap\Z} \{ p_n^{-1} \barrierprobLB \le y^{c_\varepsilon}\} \bigg] \ge 1-\varepsilon.
\]
\end{Lem}

\begin{Lem}\label{Lem:SecMomBarrCompComplicated}
For all $\varepsilon>0$, there is a $C_\varepsilon>0$ such that
\[
\liminf_{n\to \infty}\mathbb{P}\left[\sum_{k=0}^{\lfloor\log(n)^7\rfloor } e^{-\vartheta^\ast \frownacc{h}_n(k)} (q_{k,\text{start}})^{-1}\le C_\varepsilon \right] \ge 1-\varepsilon.
\]
\end{Lem}

We are ready to prove Proposition~\ref{Prop:EZn} and note that the proof very closely resembles the proof of Lemma~\ref{Lem:Anyxi}.

\begin{proof}[Proof of Proposition~\ref{Prop:EZn} assuming Lemma~\ref{Lem:SecMomBarrCompLB}]
We recall that $\zeta(X_r) = \xi(X_r)-\mathrm{es}$, $\vartheta^\ast = \mathrm{es}-\overline{\eta}$. Using Lemma~\ref{Lem:MT1} we get that
\begin{align*}
E_0^\xi[Z_n(y)] &= E_0^\xi\left[\exp\left(\int_0^{\frownacc{t}_{n;y}^{y_0}} \xi(X_r)\;\mathrm{d}r\right);\mathcal{B}_{[n],-\frownacc{m}_{n,\frownacc{h}}}^{y,J_{y_0}}(H_\cdot-(\vartheta^\ast)^{-1}K_\cdot^\zeta) \right]\\
&= e^{\vartheta^\ast \frownacc{t}_{n;y}^{y_0}}E_0^\xi\left[\exp\left(\int_0^{\frownacc{t}_{n;y}^{y_0}} \zeta(X_r)+\overline{\eta}\right);\mathcal{B}_{[n],-\frownacc{m}_{n,\frownacc{h}}}^{y,J_{y_0}}(H_\cdot-(\vartheta^\ast)^{-1}K_\cdot^\zeta) \right].
\end{align*}
On $\mathcal{B}_{[n],-\frownacc{m}_{n,\frownacc{h}}}^{y,J_{y_0}}(H_\cdot-(\vartheta^\ast)^{-1}K_\cdot^\zeta)$ we have that $\frownacc{t}_{n;y}^{y_0}-H_n\le 1$. We use this and the definition of $\frownacc{t}_{n;y}^{y_0}$ from \eqref{Def:tnyxk} to get 
\begin{align*}
E_0^{\xi}[Z_n(y)] &\ge e^{-\vartheta^\ast } e^{K_n^{\zeta}+\vartheta^\ast(y_0-y)} p_n^{-1} E_0^\xi\left[\exp\left(\int_0^{H_n} \zeta(X_r)+\overline{\eta}\right);\mathcal{B}_{[n],-\frownacc{m}_{n,\frownacc{h}}}^{y,J_{y_0}}(H_\cdot-(\vartheta^\ast)^{-1}K_\cdot^\zeta) \right]\\
& \stackrel{\eqref{eq:tilt}}{=}  e^{\vartheta^\ast(y_0-1)} e^{-\vartheta^\ast y} p_n^{-1}\barrierprobLB
\end{align*}
and applying Lemma~\ref{Lem:SecMomBarrCompLB} finishes the proof.
\end{proof}

\begin{figure}
\centering
\includegraphics[scale = 0.9]{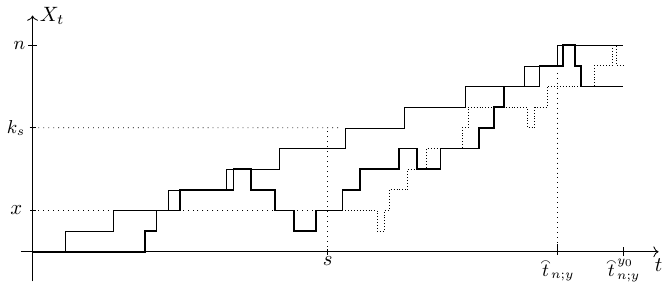}
\caption{Visualization of some quantities in Definition~\protect{\ref{Def:TwoPartOnePartDefs}} in terms of $(X_t)_{t\le \frownacc{t}_{n;y}^{y_0}}$. The thin black line is the barrier used in \eqref{Def:bprobLB}. The thick black line until $s$ is the path of a particle realizing the event $A_{n}^y(s)$ from \eqref{eq:Anys}. The thick black line after $s$ is the path of a particle realizing the event in the definition of $q_n(s,x)$ from \protect{\eqref{eq:qnsx}}. The combination of the thick black line and the densely dotted line are the paths of a pair of particles counted in $Z_n(y)$.}
\label{fig:SecMom}
\end{figure}

\subsubsection[\texorpdfstring{Upper bound on $E_0^\xi[Z_n(y)^2]$}{Upper bound on the second moment}]{Upper bound on $E_0^\xi[Z_n(y)^2]$} \label{Sec:Zny2UB}
The random variable $Z_n(y)^2$ counts the number of pairs of particles, which both follow a certain path. Heuristically, we decompose this according to the time at which these particles split. In continuous time the formal tool for this is the many-to-two lemma \eqref{eq:MT2}. We introduce a few additional definitions to shorten the displays in this calculation.
\begin{Def}\label{Def:TwoPartOnePartDefs}
Let $n\in\N$, $y>0$. For $s\in [0,\frownacc{t}_{n;y}^{y_0}]$ we define $k_s := \inf \{ k\in [n] : \frownacc{t}_{n;y}(k)>s\}\wedge n$ and set 
\begin{equation}\label{eq:Anys}
A_n^{y}(s) := \{\forall_{k< k_s}\, y+H_k-(\vartheta^\ast)^{-1}K_k^\zeta-\frownacc{m}_{n,\frownacc{h}}(k)\ge 0, H_{k_s} > s \}.
\end{equation}
Additionally, for $k\ge 0$ we set $\zeta^{(k)}(\cdot) := \zeta(k+\cdot)$ and for $s\in[0,\frownacc{t}_{n;y}^{y_0}]$, $x\in\Z$ we set
\begin{equation} \label{eq:qnsx}
q_n(s,x) := P^{\zeta^{(x)}}\big[\mathcal{B}_{\{k_s,\dots, n\},-\frownacc{m}_{n,\frownacc{h}}(\cdot)}^{y+s,J_{y_0}}(H_{\cdot-x}-(\vartheta^\ast)^{-1}K_{\cdot}^\zeta)\big].
\end{equation}
We set $I_{k,n} := [\frownacc{t}_{n;y}(k),\frownacc{t}_{n;y}(k+1)]$ for $k\le n-2$ and $I_{n-1,n} := [\frownacc{t}_{n;y}(n-1), \frownacc{t}_{n;y}^{y_0}(n)]$. We note that by our choice of $\delta$ in Definition~\ref{Def:HelpstuffSecMoment} all of these intervals are non-empty and their union covers $[0,\frownacc{t}_{n;y}^{y_0}(n)]$. Finally, we define
\begin{equation} 
q_n^\ast(k) := \sup_{s\in I_{k,n},  x\le k} e^{-\mathrm{sgn}(x)\sum_{j= (x\wedge 0)+1}^{0\vee x} L_j^{\zeta}(\overline{\eta})+\sum_{j=1}^{k} L_j^\zeta(\overline{\eta})} q_n(s,x).\label{eq:qnastk}
\end{equation}
 See Figure~\ref{fig:SecMom} for a visualization of some of the defined quantities.
\end{Def}

The next Lemma splits $\barrierprobLB$ into what happens on $[0,s]$ and what happens on $[s,n]$ using the Markov property at time $s$. 
\begin{Lem}\label{Lem:PathInterr}
For all $s\ge 0$ we have that 
\begin{align*}
\barrierprobLB =E_0^\xi\bigg[\exp\bigg(\int_0^s \xi(X_r)\;\mathrm{d}r\bigg) e^{-\mathrm{sgn}(X_s)\sum_{j= X_s\wedge 0+1}^{X_s\vee 0}L_j^{\zeta}(\overline{\eta})-\vartheta^\ast s}q_n(s,X_s); A_n^{y}(s)\bigg].
\end{align*}
\end{Lem}
\begin{proof}

Let $(\wt{X}_r)_{r\ge0}$ be a time-continuous simple random walk, which is independent of $(X_s)_{s\ge0}$ and set $\wt{H}_k := H_k(\wt{X}_r)$ to be the hitting time of $k$ for that time-continuous simple random walk. 
By using  \eqref{eq:tilt}, \eqref{Def:Knzeta} and the Markov property for $(X_r)_{r\ge0}$ at time $s$ we have that
\begin{align*}
\barrierprobLB &= e^{K_n^\zeta}E_0^\xi\left[\exp\left(\int_0^{H_n}(\zeta(X_s)+\overline{\eta}\;\mathrm{d}s\right);\mathcal{B}_{[n],-\frownacc{m}_{n,\frownacc{h}}(k)}^{y,J_{y_0}}\left(H_{\cdot}-(\vartheta^\ast)^{-1}K_{\cdot}^{\zeta}\right)\right]\\
&=e^{K_n^\zeta}E_0^\xi\Bigg[\exp\left(\int_0^s (\zeta(X_r)+\overline{\eta})\;\mathrm{d}r\right)\mathbf{1}_{A_n^{y,y_0}(s)}\\
&\hskip3cm\times E_{X_s}^\xi\left[\exp\left(\int_0^{H_n-s} \zeta(\wt{X}_r)\;\mathrm{d}r )\right);\mathcal{B}_{\{k_s,\dots, n\}, -\frownacc{m}_{n,\frownacc{h}}}^{y,J_{y_0}}\left(H_{\cdot}-(\vartheta^\ast)^{-1}K_{\cdot}^{\zeta}\right)\right]\Bigg]\\
&=e^{K_n^\zeta}E_0^\xi\Bigg[\exp\left(\int_0^s \xi(X_r)\;\mathrm{d}r\right)e^{-\vartheta^\ast s}\mathbf{1}_{A_n^{y,y_0}(s)}\\
&\hskip3cm\times E_0^\xi\left[\exp\left(\int_0^{\wt{H}_{n-X_s}} \zeta^{(X_s)}(\wt{X}_r)\;\mathrm{d}r )\right);\mathcal{B}_{\{k_s,\dots, n\}, -\frownacc{m}_{n,\frownacc{h}}}^{y+s,J_{y_0}}\left(\wt{H}_{\cdot-X_s}-(\vartheta^\ast)^{-1}K_{\cdot}^{\zeta}\right)\right]\Bigg]\\
&= e^{K_n^\zeta}E_0^\xi\left[\exp\left(\int_0^s \xi(X_r)\;\mathrm{d}r\right) e^{-\vartheta^\ast s} q(s,X_s) e^{(n-X_s)\overline{L}_{n-X_s}^{\zeta^{(X_s)}}(\overline{\eta}) }; A_n^{y}(s)\right]\\
&=E_0^\xi\left[\exp\left(\int_0^s \xi(X_r)\;\mathrm{d}r\right) e^{-\mathrm{sgn}(X_s)\sum_{j= X_s\wedge 0+1}^{X_s\vee 0}L_j^{\zeta}(\overline{\eta})-\vartheta^\ast s}q_n(s,X_s); A_n^{y}(s)\right]. \qedhere
\end{align*}
\end{proof} 
The next Lemma reduces controlling $E_0^\xi[Z_n(y)^2]$ to barrier events for a single particle.
\begin{Lem}\label{Lem:SecondMomMain}
For all $n\in\N$, $y>0$ we have that 
\begin{align*}
E_0^\xi[Z_n(y)^2]-E_0^\xi[Z_n(y)] \le C\sum_{k=0}^{n-1} p_n^{\sigma_k^2/\sigma_n^2-2}e^{-\vartheta^\ast y-\vartheta^\ast \frownacc{h}_n(k)} q_n^\ast(k)\barrierprobLB.
\end{align*}
\end{Lem}
\begin{proof}
Applying the many-to-two formula\footnote{We note that for this it is important that we have chosen the sequence $((\frownacc{t}_{n;y}(k))_{k\le n-1}, \frownacc{t}_{n;y}^{y_0})$ to be increasing.}  \eqref{eq:MT2} while recalling the definition \eqref{Def:Zny} and using that $\xi(\cdot)$ is uniformly bounded yields that
\begin{align}
&E_0^\xi[Z_n(y)^2]  \notag\\
&\le E_0^\xi\left[\exp\left(\int_0^{\frownacc{t}_{n;y}^{y_0}}\xi(X_r)\;\mathrm{d}r\right);\mathcal{B}_{[n],-\frownacc{m}_{n,\frownacc{h}}}^{y,J_{y_0}}(H_\cdot-(\vartheta^\ast)^{-1}K_\cdot^\zeta)\right]+C\int_0^{\frownacc{t}_{n;y}^{y_0}} E_0^\xi\Bigg[ \exp\left(\int_0^s \xi(X_r)\;\mathrm{d}r\right)\mathbf{1}_{A_n^{y,y_0}(s)} \notag\\
&\qquad\times E_{0}^\xi\left[\exp\left(\int_0^{\frownacc{t}_{n;y}^{y_0}-s}\xi(\wt{X}_r+X_s)\;\mathrm{d}r\right);  \mathcal{B}_{\{k_s,\dots, n\}, -\frownacc{m}_{n,\frownacc{h}}(\cdot)}^{y+s,J_{y_0}}(\wt{H}_{\cdot-X_s}-(\vartheta^\ast)^{-1}K_{\cdot}^\zeta) | X_s\right]^2\Bigg]\;\mathrm{d}s, \label{eq:Zny21}
\end{align}
where $(\wt{X}_s)_{s\ge0}$ is a time-continuous simple random walk independent of $(X_s)_{s\ge0}$, and $\wt{H}_k := H_k((\wt{X}_s)_{s\ge0})$ denotes the time this random walk hits $k$. The integrand of the second summand in \eqref{eq:Zny21} corresponds to a pair of particles, which both are counted in $Z_n(y)$, splitting at time $s$. See Figure~\ref{fig:SecMom} for the illustration of the paths of such a pair of particles. The rough plan going forward is to do a worst-case bound on one of the paths after time $s$, thus relating it to $q_n^\ast(k)$, and to combine the other path with the path until time $s$ to get a factor $\barrierprobLB$ by Lemma~\ref{Lem:PathInterr}. 

 We note that by applying the many-to-one formula \eqref{eq:MT1} the first summand in the last display equals $E_0^\xi[Z_n(y)]$. We recall that $\zeta^{(k)}(\cdot) = \zeta(k+\cdot)$ denotes environment shifted in space. We can upper bound the probability which is conditioned on $X_s$ in the last display via
\begin{align}
 &E_{0}\left[\exp\left(\int_0^{\frownacc{t}_{n;y}^{y_0}-s}\xi(\wt{X}_r+X_s)\;\mathrm{d}r\right); \mathcal{B}_{\{k_s,\dots, n\}, -\frownacc{m}_{n,\frownacc{h}}(\cdot)}^{y+s,J_{y_0}}(\wt{H}_{\cdot-X_s}-(\vartheta^\ast)^{-1}K_{\cdot}^\zeta) | X_s\right] \notag\\
 &\le e^{\vartheta^\ast(\frownacc{t}_{n;y}^{y_0}-s)}E_0^\xi\left[\exp\left(\int_0^{\wt{H}_{n-X_s}} \zeta(\wt{X}_r+X_s)+\overline{\eta}\;\mathrm{d}r \right) ; \mathcal{B}_{\{k_s,\dots, n\}, -\frownacc{m}_{n,\frownacc{h}}(\cdot)}^{y+s,J_{y_0}}(\wt{H}_{\cdot-X_s}-(\vartheta^\ast)^{-1}K_{\cdot}^\zeta) | X_s \right]\notag \\
 &= e^{\vartheta^\ast(\frownacc{t}_{n;y}^{y_0}-s)} e^{(n-X_s) \overline{L}_{n-X_s}^{\zeta^{(X_s)}}(\overline{\eta})}P^{\zeta^{(X_s)},\overline{\eta}}\left[ \mathcal{B}_{\{k_s,\dots, n\}, -\frownacc{m}_{n,\frownacc{h}}(\cdot)}^{y+s,J_{y_0}}(\wt{H}_{\cdot-X_s}-(\vartheta^\ast)^{-1}K_{\cdot}^\zeta) | X_s\right] \notag\\
 &\le p_n^{-1}e^{\vartheta^\ast (y_0-y-s)} e^{-\mathrm{sgn}(X_s)\sum_{j=  X_s\wedge 0+1 }^{0\vee X_s} L_j^{\zeta}(\overline{\eta})} q_n(s,X_s), \label{eq:Zny22}
\end{align}
Plugging \eqref{eq:Zny22} into \eqref{eq:Zny21} yields that 
\begin{align}
&E_0^\xi[Z_n(y)^2]-E_0^\xi[Z_n(y)] \notag\\
&\le p_n^{-2}e^{2\vartheta^\ast(y_0-y)}c\int_0^{\frownacc{t}_{n;y}^{y_0}} E_0^\xi\left[\exp\left(\int_0^s \xi(X_r)\;\mathrm{d}r\right) e^{-2\mathrm{sgn}(X_s)\sum_{j=X_s\wedge0+1}^{0\vee X_s}L_j^{\zeta}(\overline{\eta})-2\vartheta^\ast s}q_n(s,X_s)^2; A_n^{y}(s) \right]\;\mathrm{d}s. \label{eq:Zny2-Zny}
\end{align}
Next, we upper bound the integral term. Motivated by Lemma~\ref{Lem:PathInterr} we want to pull one copy of $e^{-\mathrm{sgn}(X_s)\cdots} q_n(s,X_s)$ out of the expectation by doing a worst case upper bound. For this, we recall the definition $I_{k,n} = [\frownacc{t}_{n;y}(k),\frownacc{t}_{n;y}(k+1)]$, $k\le n-2$, $I_{n-1,n}= [\frownacc{t}_{n;y}(n-1), \frownacc{t}_{n;y}^{y_0}]$. For $s\in I_{k,n}$ we can upper bound $e^{-\vartheta^\ast s} \le e^{-\frownacc{t}_{n;y}(k)}$. Furthermore, we have $k_s = k+1$, such that on $A_n^{y}(s)$ we have $H_{k+1}>s$ and thus $X_s\le k$. Overall, we get
\begin{align}
\int_0^{\frownacc{t}_{n;y}^{y_0}} E_0^\xi&\left[\exp\left(\int_0^s \xi(X_r)\;\mathrm{d}r\right) e^{-2\mathrm{sgn}(X_s)\sum_{j= X_s\wedge 0+1}^{X_s\vee 0}L_j^{\zeta}(\overline{\eta})-2\vartheta^\ast s}q_n(s,X_s)^2; A_{n}^{y,y_0}(s) \right]\;\mathrm{d}s \notag \\
&= \sum_{k=0}^{n-1} \int_{I_{k,n}}E_0^\xi\left[\exp\left(\int_0^s \xi(X_r)\;\mathrm{d}r\right) e^{-2\mathrm{sgn}(X_s)\sum_{j= X_s\wedge 0+1}^{X_s\vee 0}L_j^{\zeta}(\overline{\eta})-2\vartheta^\ast s}q_n(s,X_s)^2;  A_{n}^{y,y_0}(s) \right]\;\mathrm{d}s \notag\\
&\le \sum_{k=0}^{n-1} \left(\sup_{s\in I_{k,n}, x\le k} e^{-\mathrm{sgn}(x)\sum_{j= (x\wedge 0)+1}^{0\vee x} L_j^{\zeta}(\overline{\eta})+\sum_{j=1}^{k} L_j^\zeta(\overline{\eta})-\vartheta^\ast y-\vartheta^\ast \frownacc{h}_n(k)+\frac{\sigma_k^2}{\sigma_n^2}\log(p_n)} q_n(s,x)\right) \notag \\
&\qquad\times\int_{I_{k,n}}E_0^\xi\left[\exp\left(\int_0^s \xi(X_r)\;\mathrm{d}r\right) e^{-\mathrm{sgn}(X_s)\sum_{j= X_s\wedge 0+1}^{X_s\vee 0}L_j^{\zeta}(\overline{\eta})-\vartheta^\ast s}q_n(s,X_s); A_n^{y}(s) \right], \label{eq:SplitIntoTnSum}
\end{align}
where we have also used the definition of $\frownacc{t}_{n;y}(k)$, see \eqref{Def:tnyxk}.
Plugging \eqref{eq:SplitIntoTnSum} into \eqref{eq:Zny2-Zny} and applying Lemma~\ref{Lem:PathInterr} yields that
\begin{align}\label{eq:Zny2-ZnySecLast}
E_0^\xi[Z_n(y)^2]-E_0^\xi[Z_n(y)] &\le \sum_{k=0}^{n-1} p_n^{\sigma_k^2/\sigma_n^2-2}e^{-\vartheta^\ast y-\vartheta^\ast \frownacc{h}_n(k)} q_n^\ast(k)\barrierprobLB \cdot |I_{k,n}|.
\end{align}
Now we use the fact that $\mathbb{P}$-a.s.\@ for all $k\in\N$ we have $\chi_k^2 \in [c_\chi,C_\chi]$ and Lemma~\ref{Lem:LBpn} to see that 
\[
|I_{k,n}| \le |L_{k+1}^\zeta(\overline{\eta})|+y_0+|\frownacc{h}_n(k+1)-\frownacc{h}_n(k)|-\frac{\chi_{k+1}^2}{\vartheta^\ast \sigma_n^2}\log(p_n)\le C<\infty
\]
is bounded uniformly in $n\in\N$, $k\le n$, $y>0$. Combining the last display with \eqref{eq:Zny2-ZnySecLast} finishes the proof.
\end{proof}
Armed with Lemmata~\ref{Lem:SecMomBarrCompUB}  and~\ref{Lem:SecondMomMain} the only step missing in the proof of Proposition~\ref{Prop:EZn2} is showing that for $k\le \log(n)^7$ we have that there is some function $g$ such that $q_n^\ast(k) \le g(k)\cdot \barrierprobLB$ with $\sup_n\sum_{k=0}^{\lfloor \log(n)^7\rfloor } g(k)e^{-\vartheta^\ast \frownacc{h}_n(k)}<\infty$, which is done in the next lemma.
\begin{Lem} \label{Lem:EZn2Hilf1}
For all $\varepsilon>0$, there is a $C_\varepsilon>0$ such that
\[
\liminf_{n\to\infty} \mathbb{P}\bigg[\bigcap_{k\le \lfloor \log(n)^7\rfloor} \\{p_n^{-1}q_n^\ast(k) \le C_\varepsilon (q_{k,\text{start}})^{-1}\\}\bigg]\ge 1-\varepsilon.
\]
\end{Lem}

We postpone the proof of Lemma~\ref{Lem:EZn2Hilf1} to finish the proof of Proposition~\ref{Prop:EZn2} first.

\begin{proof}[Proof of Proposition~\ref{Prop:EZn2} assuming Lemmata~\ref{Lem:SecMomBarrCompUB},~\ref{Lem:SecMomBarrCompComplicated},~\ref{Lem:EZn2Hilf1}]
Fix $y\ge y_0$ and $\varepsilon>0$. Applying Lemma~\ref{Lem:SecondMomMain} yields that 
\begin{align}
&E_0^\xi[Z_n(y)^2]-E_0^\xi[Z_n(y)] \notag \\
&\le C\sum_{k=0}^{n-1} p_n^{\sigma_k^2/\sigma_n^2-2}e^{-\vartheta^\ast y-\vartheta^\ast \frownacc{h}_n(k)} q_n^\ast(k)\barrierprobLB\notag \\
&\le C\sum_{k=0}^{\lfloor\log(n)^7\rfloor} e^{-\vartheta^\ast \frownacc{h}_n(k)} p_n^{-2}e^{-\vartheta^\ast y} q_n^\ast(k) \barrierprobLB+C\sum_{k= \lceil\log(n)^7\rceil}^n e^{-\vartheta^\ast \frownacc{h}_n(k)} p_n^{\sigma_k^2/\sigma_n^2-2} e^{-\vartheta^\ast y} \barrierprobLB. \label{eq:ProofEZn2Zwischen1}
\end{align}
For the first summand we combine Lemmata~\ref{Lem:SecMomBarrCompUB},~\ref{Lem:SecMomBarrCompComplicated} and~\ref{Lem:EZn2Hilf1} to get that there exists a $ c_\varepsilon>0$ such that
\begin{equation}
1-\varepsilon/2 \le \liminf_{n\to \infty} \mathbb{P}\left[C\sum_{k=0}^{\lfloor \log(n)^7\rfloor} e^{-\vartheta^\ast \frownacc{h}_n(k)} p_n^{-2}e^{-\vartheta^\ast y} q_n^\ast(k) \barrierprobLB \le  e^{-\vartheta^\ast y}y^{c_\varepsilon}\right], \label{eq:ProofEZn2Sum1}
\end{equation}
where we have absorbed multiplicative constants into the $c_\varepsilon$, which is possible since $y\ge y_0>1$.

For the second summand we apply Lemma~\ref{Lem:SecMomBarrCompUB} and use that by Lemma~\ref{Lem:LBpn} there is a constant $\gamma>0$ such that $\mathbb{P}$-a.s.\@ $\limsup_{n\to\infty} \frac{|\log(p_n)|}{\log(n)} \le \gamma$ as well as $\sup_n\sum_{k=\lceil \log(n)^7\rceil}^n e^{-\vartheta^\ast \frownacc{h}_n(k)} n^{\gamma\cdot (1-\sigma_k^2/\sigma_n^2)} <\infty$ to deduce that there is a $c_\varepsilon'>0$ such that
\begin{equation}
1-\varepsilon/2\le\liminf_{n\to\infty} \mathbb{P}\left[C\sum_{k=\lceil \log(n)^7\rceil}^n e^{-\vartheta^\ast \frownacc{h}_n(k)}p_n^{\sigma_k^2/\sigma_n^2-2}e^{-\vartheta^\ast y}\barrierprobLB \le  e^{-\vartheta^\ast y} y^{c_\varepsilon'} \right]. \label{eq:ProofEZn2Sum2}
\end{equation}
Combining \eqref{eq:ProofEZn2Zwischen1} with \eqref{eq:ProofEZn2Sum1} and \eqref{eq:ProofEZn2Sum2} yields that for $c_{\varepsilon}'' = \max\{c_\varepsilon, c_{\varepsilon}'\}$
\[
1-\varepsilon \le\liminf_{n\to \infty} \mathbb{E}\left[ E_0^\xi[Z_n(y)^2]-E_0^\xi[Z_n(y)] \le  2e^{-\vartheta^\ast y} y^{c_\varepsilon''}\right],
\]
finishing the proof of Proposition~\ref{Prop:EZn2}.
\end{proof}
Thus to finish the proof of Theorem~\ref{Theo:LBLT}, barring the proof of the barrier computations in Lemmata~\ref{Lem:SecMomBarrCompLB}--\ref{Lem:SecMomBarrCompComplicated}, we only need to prove Lemma~\ref{Lem:EZn2Hilf1}.
\begin{proof}[Proof of Lemma~\ref{Lem:EZn2Hilf1} assuming Lemmata~\ref{Lem:SecMomBarrCompUB},~\ref{Lem:SecMomBarrCompComplicated}]
Fix $k\le \log(n)^7$, $s\in [\frownacc{t}_{n;y}(k),\frownacc{t}_{n;y}(k+1)]$. First, we note that $k_s = k+1$, such that 
\[
q_n(s,x) = P^{\zeta^{(x)}}\left[\mathcal{B}_{\{k+1,\dots, n\},-\frownacc{m}_{n,\frownacc{h}}(\cdot)}^{y+s,J_{y_0}}(H_{\cdot-x}-(\vartheta^\ast)^{-1}K_{\cdot}^{\zeta})\right].
\]
Now, we set $H_{x,k}$ to be the hitting time of $k$ when starting at $x$. By partitioning with regard to the value of $H_{x,k}$ we get 
\begin{equation}\label{eq:qnsxb}
\begin{aligned}
q_n(s,x) &\le \sum_{l=0}^\infty P^{\zeta}[H_{x,k} \in (l,l+1]] \sup_{z\in(l,l+1]} P^{\zeta^{(k)},\overline{\eta}}\left[ \mathcal{B}_{\{k+1,\dots, n\},-\frownacc{m}_{n,\frownacc{h}}(\cdot)}^{y+s+z,J_{y_0}}(H_{\cdot}-(\vartheta^\ast)^{-1}K_{\cdot}^{\zeta})\right]\\
&=\sum_{l=0}^\infty P^{\zeta^{(x)}}[H_{k-x} \in (l,l+1]] \sup_{z\in(l,l+1]} P^{\zeta^{(k)},\overline{\eta}}\left[ \mathcal{B}_{\{k+1,\dots, n\},-\frownacc{m}_{n,\frownacc{h}}(\cdot)}^{y+s+z,J_{y_0}}(H_{\cdot}-(\vartheta^\ast)^{-1}K_{\cdot}^{\zeta})\right],
\end{aligned}
\end{equation}
For $z\ge 0$ set $\{z\} := \lceil z\rceil-z$ to be the fractional part of $z$. We have
\begin{equation}\label{eq:qnsxBound}
\begin{aligned}
\frownacc{p}_n^{y_0}(y_0+\lceil z\rceil-2) &\ge P^{\zeta}\left[\mathcal{B}_{[n], -\frownacc{m}_{n,\frownacc{h}}}^{y_0+\lceil z\rceil-2,J_{y_0}}(H_\cdot-(\vartheta^\ast)^{-1}K_\cdot^\zeta), X_{y-y_0+s+2-\{z\}}=  k\right]\\
&\ge P^{\zeta}\left[ \mathcal{B}_{[k],-\frownacc{m}_{n,\frownacc{h}}}^{y_0+\lceil z\rceil-2,\R}(H_\cdot-(\vartheta^\ast)^{-1}K_{\cdot}^\zeta), X_{y-y_0+s+2-\{z\}} = k\right]\cdot \\
&\qquad\cdot P^{\zeta^{(k)},\overline{\eta}}\left[ \mathcal{B}_{\{k+1,\dots, n\}, -\frownacc{m}_{n,\frownacc{h}}(\cdot)}^{y+s+z,J_{y_0}}(H_{\cdot-k}-(\vartheta^\ast)^{-1}K_{\cdot}^{\zeta}) \right].
\end{aligned}
\end{equation}
We note that the event $\mathcal{B}_{[k],-\frownacc{m}_{n,\frownacc{h}}}^{y_0+\lceil z\rceil-2,\R}(H_\cdot-(\vartheta^\ast)^{-1}K_{\cdot}^\zeta)$ is monotonically increasing in $\lceil z\rceil$, such that plugging \eqref{eq:qnsxBound} into \eqref{eq:qnsxb} yields that
\begin{multline*}
q_n(s,x) \\
\le \sup_{z\in[1,2]} \left(P^{\zeta}\left[\mathcal{B}_{[k],-\frownacc{m}_{n,\frownacc{h}}}^{y_0-1,\R}(H_{\cdot}-(\vartheta^\ast)^{-1}K_{\cdot}^{\zeta}), X_{y-y_0+s+z} = k\right]\right)^{-1}\sum_{l=0}^\infty \frownacc{p}_n^{y_0}(y_0+l)P^{\zeta^{(x)}}[H_{k-x}\in [l,l+1]] .
\end{multline*}
We note that using $s\in [\frownacc{t}_{n;y}(k),\frownacc{t}_{n;y}(k+1)]$ and \eqref{Def:tnyxk} we get for $z\in [1,2]$
\[
y-y_0+s+z  \in [(\vartheta^\ast)^{-1}K_k^{\zeta}+\frownacc{m}_{n,\frownacc{h}}(k)-y_0+1, (\vartheta^\ast)^{-1}K_{k+1}^{\zeta}+\frownacc{m}_{n,\frownacc{h}}(k+1)-y_0+2] \subseteq (\vartheta^\ast)^{-1}K_k^\zeta+\frownacc{m}_{n,\frownacc{h}}(k)-y_0+[1,2].
\]
 Thus for $z\in [0,1]$
\begin{align*}
&P^{\zeta}\left[\mathcal{B}_{[k],-\frownacc{m}_{n,\frownacc{h}}}^{y_0-1,\R}(H_{\cdot}-(\vartheta^\ast)^{-1}K_\cdot^\zeta), X_{y-y_0+s+z} = k\right]\\
&\ge \inf_{x\in [1,2]} P^{\zeta}\left[\mathcal{B}_{[k],-\frownacc{m}_{n,\frownacc{h}}}^{y_0-1,\R}(H_{\cdot}-(\vartheta^\ast)^{-1}K_\cdot^\zeta), X_{(\vartheta^\ast)^{-1}K_k^{\zeta}+\frownacc{m}_{n,\frownacc{h}}(k)-y_0+x} = k\right]\\
&\ge \inf_{x\in [1,2]} P^{\zeta}\left[\mathcal{B}_{[k],-\frownacc{m}_{n,\frownacc{h}}}^{y_0-1,J_x}(H_{\cdot}-(\vartheta^\ast)^{-1}K_{\cdot}^{\zeta})\right]\cdot \inf_{r\in [0,2]} P^{\zeta^{(k)},\overline{\eta}}\left[ \forall_{u\le r}\, X_u = k| X_0= k\right]\\
&\ge c q_{k,\text{start}},
\end{align*}
for some $c\in (0,1]$. Here the last step uses that $\zeta^{(k)}(x)$, $x\in\Z$, is uniformly lower bounded, i.e.\@ even under $P^{\zeta^{(k)},\overline{\eta}}$ the jump rate of $(X_s)_{s\ge0}$ is uniformly bounded away from infinity, such that the probability for the particle not to jump for two units otf time is uniformly lower bounded.

This yields that 
\[
q_n(s,x) \le C (q_{k,\text{start}})^{-1} \sum_{l=0}^\infty \frownacc{p}_n^{y_0}(y_0+l)P^{\zeta^{(x)}}\left[H_{k-x}\in[l,l+1]\right]
\]
and it remains to show that for all $\varepsilon>0$ there is a $C_\varepsilon>0$ such that 
\begin{equation}
1-\varepsilon\le \liminf_{n\to\infty} \mathbb{P}\left[ \sup_{x\le k}  e^{-\mathrm{sgn}(x) \sum_{j = (x\wedge 0)+1}^{x\vee 0} L_j^{\zeta}(\overline{\eta})+\sum_{j=1}^k L_j^{\zeta}(\overline{\eta})}\sum_{l=0}^{\infty} p_n^{-1}\frownacc{p}_n^{y_0}(y_0+l)P^{\zeta^{(x)}}\left[H_{k-x}\in[l,l+1]\right] \le C_\varepsilon \right]. \label{eq:AimSecMomSecPart}
\end{equation}
Heuristically, we know from Lemma~\ref{Lem:SecMomBarrCompUB} that $\frac{\frownacc{p}_n^{y_0}(y_0+l)}{p_n} \le  (l+1)^{c_\varepsilon}$ with high probability, such that 
\[
\sum_{l=0}^\infty p_n^{-1}\frownacc{p}_n^{y_0}(y_0+l)P^{\zeta^{(x)}}\left[H_{k-x}\in [l,l+1]\right] \le C_{\varepsilon}E^{\zeta^{(x)}}[H_{k-x}^{c_\varepsilon}].
\] Under $P^{\zeta^{(x)}}$ the random walk $(X_s)_{s\ge0}$ has drift towards $+\infty$, such that $E^{\zeta^{(x)}}[H_{k-x}^{c_\varepsilon}]$ is polynomial in $(k-x)$, which allows us to conclude \eqref{eq:AimSecMomSecPart}. We make the above outline rigorous next.

Fix $\varepsilon>0$. We split the sum over $l$ into $l\le \log(n)^2-y_0$ and $l>\log(n)^2-y_0$, apply Lemma~\ref{Lem:SecMomBarrCompUB} in the first case and bound $\frownacc{p}_n^{(y_0)}(y_0+l)\le 1$ in the latter case to get that there exists a $c_{\varepsilon}>0$ such that
\begin{equation} \label{eq:MainThingyHkx}
\begin{aligned}
1-\varepsilon\le \liminf_{n\to\infty}\mathbb{P}\Bigg[
 \sum_{l=0}^\infty p_n^{-1}\frownacc{p}_{n}^{y_0}(y_0+l)& P^{\zeta^{(x)}}[H_{k-x}\in[l,l+1]] \\
 &\quad \le  C_{\varepsilon}E^{\zeta^{(x)}}[H_{k-x}^{c_\varepsilon}]+p_n^{-1} P^{\zeta^{(x)}}[H_{k-x} \ge \log(n)^2-y_0]
\Bigg].
\end{aligned}
\end{equation}
We handle the two summands on the right-hand side separately. For the first summand, we let $\tau_k$ be the hitting time of $k$, when starting at $k-1$. Then
\begin{equation} \label{eq:RewriteH}
H_{k-x}-E^{\zeta^{(x)}}[H_{k-x}] = \sum_{l=1}^{k-x} \tau_l-\sum_{l=1}^{k-x} E^{\zeta^{(x)}}[\tau_l].
\end{equation}
From \cite[Eq. (4.9), A.1.]{CernyDrewitz} we know that the $\tau_k$ have finite exponential moments near zero and (uniformly) bounded first and second moment. This together with \eqref{eq:RewriteH} and \cite[p.\@ 59, 14.]{petrov2012sums} implies that there is a $C_{\varepsilon}>0$ depending only on $c_{\varepsilon}$ such that $\mathbb{P}$-a.s.\@ 
\begin{equation}\label{eq:HkxExp}
E^{\zeta^{(x)}}\left[H_{k-x}^{c_\varepsilon}\right] \le C_\varepsilon (k-x)^{c_{\varepsilon}}.
\end{equation}
Now, we deal with the second summand on the right-hand side of \eqref{eq:MainThingyHkx}. By \cite[p.\@ 52--54, Th. 15, L. 5]{petrov2012sums} we have that there are $c, C>0$ (independent of $\zeta$) such that $\mathbb{P}$-a.s.
\begin{equation}
P^{\zeta^{(x)}}[H_{k-x}\ge \log(n)^2-y_0] \le e^{-\min\left\{\frac{(\log(n)^2-y_0)^2}{2C(k-x)},\frac{c(\log(n)^2-y_0)}{2}\right\}}. \label{eq:PzetetHkxBig}
\end{equation}
Finally, we recall that by Lemma~\ref{Lem:LBpn} there is a $\gamma>0$ such that $\mathbb{P}$-a.s.\@ we have $\limsup_{n\to\infty} \frac{|\log(p_n)|}{\log(n)}\le \gamma$. Plugging this and \eqref{eq:RewriteH}, \eqref{eq:PzetetHkxBig} into \eqref{eq:AimSecMomSecPart} yields that
\begin{equation} \label{eq:SecMomSecPartFastFinal}
\begin{aligned}
&1-\varepsilon \\
&\le \liminf_{n\to\infty}\mathbb{P}\left[ \sum_{l=0}^\infty p_n^{-1}\frownacc{p}_n^{y_0}(y_0+l)P^{\zeta^{(x)}}[H_{k-x}\in[l,l+1]] \le C_\varepsilon (k-x)^{c_\varepsilon}+n^{\gamma}e^{-\min\left\{\frac{(\log(n)^2-y_0)^2}{2C(k-x)},\frac{c(\log(n)^2-y_0)}{2}\right\}}\right],
\end{aligned}
\end{equation}
Now, we can combine \eqref{eq:SecMomSecPartFastFinal} with the fact that there are $C_{L^\zeta}, c_{L^\zeta}>0$ such that $\mathbb{P}$-a.s.\@ for all $k\in\N$ we have $L_k^{\zeta}(\overline{\eta}) \in[-C_{L^\zeta},-c_{L^\zeta}]$ to conclude \eqref{eq:AimSecMomSecPart}, thus finishing the proof. 
\end{proof}
\subsection[\texorpdfstring{Upper bound on the right tail of $\mathbf{H}_n-m_n$ -- cutting the tree}{Upper bound on the right tail -- cutting the tree}]{Upper bound on the right tail of $\mathbf{H}_n-m_n$ -- cutting the tree}\label{Sec:UseTree}

This section aims  to provide an upper bound on the right tail of $\mathbf{H}_n-m_n$ supplementing Theorem~\ref{Theo:UBLT} and thus overall providing all ingredients necessary for the proof of Theorem~\ref{Theo:Main}.
\begin{Sa}\label{Theo:UBUT}
For all $\varepsilon>0$ there is a $y_\varepsilon\ge 0$ such that
\[
\liminf_{n\to \infty} \mathbb{P}\left[P_0^\xi\left[\mathbf{H}_n-m_n\ge y_\varepsilon\right]\le \varepsilon\right]\ge 1-\varepsilon.
\]
\end{Sa}
The idea of the proof is to let the process run for some time $t_\varepsilon>0$ and use that at time $t_\varepsilon$ there is an exponential in $t_\varepsilon$ number of particles at location $0$.  Applying Theorem~\ref{Theo:LBLT} to these particles shows that we can lower bound the probability of the children of any of these particles to reach $n$ by time $m_n+y_\varepsilon$.  However, since under the quenched law the children of different particles at time $t_\varepsilon$ behave independently, it is quite likely that one of these children reaches $n$ by time $m_n+y_\varepsilon$.  First, we cite a lower bound from \cite{CernyDrewitz} on the size of $N(t,0) := \{Y\in N(t) : Y_t = 0\}$, i.e.\@ the number of particles at the origin at time $t$. 
\begin{Lem}[{\cite[Lemma 6.8]{CernyDrewitz}}] \label{Lem:Part0}
There exist $t_0<\infty$ and $C_0>1$ such that for all $t\ge t_0$ $\mathbb{P}$-a.s.
\begin{equation}
P_0^\xi\left[ \# N(t,0)\le C_0^t\right]\le C_0^{-t}. \label{eq:Part0}
\end{equation}
\end{Lem}
\begin{proof}[Proof of Theorem~\ref{Theo:UBUT}]
Fix $\varepsilon>0$,  $y\ge 0$. We take $t_0$ as in Lemma~\ref{Lem:Part0}. By conditioning on $\mathcal{F}_t$ we see that for $t\ge t_0$%
\begin{align*}
P_0^\xi\left[\mathbf{H}_n-m_n \ge y \right]&\le E^\chi_0\smash[t]{\left[ \prod_{Y\in N(t,0)} P_0^{\xi}\left[ H_n( (Y_s)_{s\ge t})-m_n\ge y-t\right] \right]}\\
&\le P_0^\xi[\#N(t,0)\le C_0^t]+P_0^\xi\left[\mathbf{H}_n-m_n\ge y-t\right]^{C_0^t}\\
&\stackrel{\eqref{eq:Part0}}{\le} C_0^{-t}+P_0^\xi\left[\mathbf{H}_n-m_n\ge y-t\right]^{C_0^t}.
\end{align*}
Applying Theorem~\ref{Theo:LBLT} to the last display yields that for $t\ge y$ there exists a $c_\varepsilon>0$, such that 
\begin{equation}
1-\varepsilon \le \liminf_{n\to\infty} \mathbb{P}\left[ P_0^\xi\left[\mathbf{H}_n-m_n \ge y \right] \le  C_0^{-t}+\left(1-e^{-\vartheta^\ast(y_0-y+t)} (y_0-y+t)^{-c_\varepsilon}\right)^{C_0^t} \right]. \label{eq:ZwCut}
\end{equation}
Now fix $y_\varepsilon = t_\varepsilon$ such that $C_0^{-t_\varepsilon}\le \varepsilon/2$ and $\left(1- e^{-\vartheta^\ast y_0}y_0^{-c_\varepsilon}\right)^{C_0^{t_\varepsilon}} \le \varepsilon/2$. For this choice \eqref{eq:ZwCut} implies that
\[
1-\varepsilon \le \liminf_{n\to \infty} \mathbb{P}\left[ P_0^\xi[\mathbf{H}_n-m_n\ge y_{\varepsilon}]\le \varepsilon/2+\varepsilon/2]\right]
\]
concluding the proof of Theorem~\ref{Theo:UBUT}.
\end{proof}

\subsection{Combining Theorems~\ref{Theo:UBLT},~\ref{Theo:UBUT} to prove Theorem~\ref{Theo:Main}}
Let $\varepsilon>0$ and $y\ge 0$. Set
\[
A_{n,\varepsilon}(y) := \{P_0^\xi[|\mathbf{H}_n-m_n|\ge y]\le \varepsilon/2\}.
\]
By applying Theorems~\ref{Theo:UBLT} and~\ref{Theo:UBUT} there is a $y_\varepsilon\ge y_0$ such that
\begin{equation}
\limsup_{n\to\infty}\mathbb{P}[A_{n,\varepsilon}(y_\varepsilon)^c]\le \varepsilon/2.\label{eq:Anyeps}
\end{equation}
Furthermore, we have that
\begin{align}
P[|\mathbf{H}_n-m_n|\ge y_\varepsilon] &= \mathbb{E}[P_0^\xi[\mathbf{H}_n-m_n|\ge y_{\varepsilon}]\mathbf{1}_{A_{n,\varepsilon}(y_\varepsilon)}]+\mathbb{E}[P_0^\xi[\mathbf{H}_n-m_n|\ge y_{\varepsilon}]\mathbf{1}_{\left(A_{n,\varepsilon}(y_\varepsilon)\right)^c}] \notag\\
&\le \varepsilon/2+\mathbb{P}[\left(A_{n,\varepsilon}(y_\varepsilon)\right)^c] \label{eq:Anyeps2}
\end{align}
Combining \eqref{eq:Anyeps} and \eqref{eq:Anyeps2} yields that
\[
\limsup_{n\to\infty} P[|\mathbf{H}_n-m_n|\ge y_\varepsilon]\le \varepsilon,
\]
which gives that $(\mathbf{H}_n-m_n)_{n\ge 0}$ is tight with regard to the annealed measure $P$. \hfill \qed
\subsection{Proof of Corollary~\ref{Kor:Max}}\label{Sec:CorProof}
In this section we show that there is a centering $\wt{m}_t$ such that $(M_t-\wt{m}_t)_{t\ge 0}$ is tight. Since we know that $(\mathbf{H}_n-m_n)_{n\in\N}$ is tight, and $M_{\mathbf{H}_n} = n$ by definition, it is reasonable to set $\wt{m}_t := \inf\{k\in\N : m_{k+1}\ge t\}$. If $M_t \gg \wt{m}_t$, that means that $k+C$ was hit earlier than $m_{k+C}$, which is unlikely since $(\mathbf{H}_n-m_n)_{n\in\N}$ is tight. If on the other hand $M_t\ll \wt{m}_t$ that means either that $k$ is hit way later than $m_k$ or that $k$ is hit before time $t$, but all descendants of the particle having hit $k$ are below $k-1$. The first of these options is unlikely because of the tightness of $(\mathbf{H}_n-m_n)_{n\in\N}$, the second is unlikely because \cite[Lemma 6.8]{CernyDrewitz}, which says that the number of descendants which are at the original position of the ancestor grows exponentially in time. Before we make the above outline precise we need one more growth estimate on $(m_n)_{n\in\N}$, ensuring that $m_{k+C}-m_k$ can't be too small.

\begin{Lem}\label{Lem:mkComp}
For all $\varepsilon>0$, $y>0$ there is a (deterministic) $C_{\varepsilon}>0$ such that
\begin{equation}\label{eq:GoalmkComp}
 \liminf_{k\to\infty} \mathbb{P}\left[ m_{k+C_\varepsilon}-m_k \ge y \right] \ge 1-\varepsilon.
\end{equation}
\end{Lem}
\begin{proof}
By definition \eqref{Def:pn}, \eqref{Def:mn}, we have that for $k,C\in\N$ arbitrary
\[
\vartheta^\ast(m_{k+C}-m_k) =  K_{k+C}^\zeta-K_k^\zeta+\log(p_k)-\log(p_{k+C}) = -\sum_{j=k+1}^{k+C} L_j^\zeta(\overline{\eta})+\log(p_k/p_{k+C}).
\]
Recalling that $\mathbb{P}$-a.s.\@ we have $L_k^\zeta(\overline{\eta})\in [-C_{L^\zeta},-c_{L^\zeta}]$ this yields that there is some fixed constant $C'$ such that 
\begin{equation} \label{eq:mkCmk}
m_{k+C}-m_k \ge C'\cdot C+(\vartheta^\ast)^{-1}\log(p_k/p_{k+C}). 
\end{equation}
Thus we need to establish a lower bound on $p_k/p_{k+C}$. The idea is to decompose $p_{k+C}$ according to the position of $B_k^\zeta-W_k^\zeta$, which yields
\begin{equation}\label{eq:pk/pkCLB}
\begin{aligned}
p_{k+C}&\le P^{\zeta}\left[\mathcal{B}_{[0,k],-W_\cdot}^{y_0,\R}(B_\cdot^\zeta), y_0+B_{k+C}^\zeta-W_{k+C} \in [y_0-1,y_0]\right] \\
&\le P^{\zeta}\left[\mathcal{B}_{[0,k],-W_{\cdot}}^{y_0, [0,y_0-1]}(B_\cdot^\zeta)\right]+\sum_{l=y_0}^\infty \mathbb{P}\left[\mathcal{B}_{[0,k],-W_{\cdot}}^{y_0, J_l}(B_\cdot^\zeta)\right]\cdot \max_{z\in J_l} P^{\zeta}\left[z+B_{k+C;k}^\zeta-W_{k+C;k} \in[y_0-1,y_0] \right].
\end{aligned}
\end{equation}
The first summand in the last display is like $Cp_k$ by a monotonicity argument, the second summand can be bound since we know that $\mathbb{P}\left[\mathcal{B}_{[0,k],-W_{\cdot}}^{y_0, J_l}(B_\cdot^\zeta)\right] \approx l^cp_k$, while $\max_{z\in J_l} P^{\zeta}\left[z+B_{k+C;k}^\zeta-W_{k+C;k} \in[y_0-1,y_0] \right]$ decays in $l$ like $e^{-cl^2}l^c$ so that overall the sum is like $Cp_k$. We next make this outline rigorous.

For the first summand in \eqref{eq:pk/pkCLB}, we use $g_{\mu,\sigma^2}$ to denote the Gaussian density with mean $\mu$ and variance $\sigma^2$, and write
\begin{equation} \label{eq:BrownBridgeKorMax}
P^{\zeta}\left[\mathcal{B}_{[0,k],-W_{\cdot}}^{y_0,[0,y_0-1]}(B_\cdot^{\zeta})\right] = \int_0^{y_0-1} P^{\zeta}\left[\mathcal{B}_{[0,k],-W_{\cdot}}^{y_0}(B_\cdot^{\zeta})|B_k^\zeta = z\right] g_{y_0-W_k,\sigma_k^2}(z)\;\mathrm{d}z.
\end{equation}
Now we use that the probability of barrier events is monotone in the endpoint of a Brownian bridge to get from \eqref{eq:BrownBridgeKorMax}
\begin{equation} \label{eq:pk0y0UB}
P^{\zeta}\left[\mathcal{B}_{[0,k],-W_{\cdot}}^{y_0,[0,y_0-1]}(B_{\cdot}^\zeta)\right] \le  P^{\zeta}\left[\mathcal{B}_{[0,k],-W_{\cdot}}^{y_0}(B_{\cdot}^\zeta)| B_k^\zeta = y_0-1\right] \int_0^{y_0-1} g_{y_0-W_k,\sigma_k^2}(z);\mathrm{d}z.
\end{equation}
By an analogous argument we conclude
\begin{equation} \label{eq:pkLB}
p_k \ge P^{\zeta}\left[\mathcal{B}_{[0,k],-W_{\cdot}}^{y_0}(B_{\cdot}^\zeta)| B_k^\zeta = y_0-1\right] \int_{y_0-1}^{y_0} g_{y_0-W_k,\sigma_k^2}(z)\;\mathrm{d}z.
\end{equation}
Combining \eqref{eq:pk/pkCLB}, \eqref{eq:pk0y0UB}, \eqref{eq:pkLB} yields that there is a $\wt{C}>0$ so that 
\begin{equation}\label{eq:pk/pkCLBHalfWay}
p_{k+C} \le \wt{C}p_k+\sum_{l=y_0}^\infty \mathbb{P}\left[\mathcal{B}_{[0,k],-W_{\cdot}}^{y_0, J_l}(B_\cdot^\zeta)\right]\cdot \max_{z\in J_l} P^{\zeta}\left[z+B_{k+C;k}^\zeta-W_{k+C;k} \in[y_0-1,y_0] \right].
\end{equation}
Next, we deal with the second summand in \eqref{eq:pk/pkCLBHalfWay}. Using the definition \eqref{Def:Wnzeta} of $W_n$  and \cite[Lemma A.1]{CernyDrewitz} we know that there is a $C''$ so that $|W_{k+C;k}| \in C\cdot C''$ so that Gaussian tail estimates and $\sigma_{k+C;k}^2 \in[c_\chi C, C_\chi C]$ yield that there is a $c>0$ such that 
\begin{align*}
\max_{z\in J_l} P^{\zeta}\left[z+B_{k+C;k}^\zeta-W_{k+C;k} \in[y_0-1,y_0] \right] \le e^{-c\frac{\min\{l-1-y_0-CC'',0\}^2}{C^2}}.
\end{align*}
Now, we fix $k$ with $\log(k)^2-1\ge 2(y_0+CC'')$ so that for $l\ge \log(k)^2$ we have $\min\{l-1-y_0-CC'',0\}^2 \ge z^2/4$ and get 
\begin{align*}
p_{k+C} &\le \sum_{l = y_0}^{\lfloor \log(k)^2\rfloor} P^{\zeta}\left[\mathcal{B}_{[0,k],-W_{\cdot}}^{y_0, J_l}(B_\cdot^\zeta)\right] e^{-c\frac{\min\{l-1-y_0-CC'',0\}^2}{C^2}}+\sum_{l=\lceil \log(k)^2\rceil}^{\infty} e^{-\frac{cl^2}{4C^2}}.
\end{align*}
By Lemma~\ref{Lem:FMomSimplBarrComp} we know that there is a $c_{\varepsilon}>0$ so that 
\begin{equation}
\liminf\limits_{k\to\infty} \mathbb{P}\left[\bigcap_{l\in [y_0,\log(k)^2]\cap \N} \frac{P^{\zeta}\left[\mathcal{B}_{[0,k],-W_{\cdot}}^{y_0, J_l}(B_\cdot^\zeta)\right]}{p_k} \le l^{c_\varepsilon} \right] \ge 1-\varepsilon,
\end{equation}
which yields
\begin{align*}
1-\varepsilon \le \liminf_{k\to\infty} \mathbb{P}\left[ p_{k+C} \le  \wt{C}p_k+p_k\sum_{l=y_0}^\infty e^{-c\frac{\min\{l-1-y_0-CC'',0\}^2}{C^2}}l^{c_\varepsilon}\right],
\end{align*}
and bounding $\sum_{l=y_0}^\infty e^{-c\frac{\min\{l-1-y_0-CC'',0\}^2}{C^2}}l^{c_\varepsilon} \le (CC''+y_0)^{c_{\varepsilon}+2}$ this means that there is a constant $c_{\varepsilon}'>0$ so that (for $C\ge 2$)
\begin{equation} 
1-\varepsilon \le \liminf_{k\to\infty} \mathbb{P}\left[ p_{k+C} \le  C^{c_{\varepsilon}'}p_k\right].
\end{equation}
Plugging this into \eqref{eq:mkCmk} yields that 
\[
1-\varepsilon\le \liminf_{k\to\infty} \mathbb{P}\left[ m_{k+C}-m_k \ge C'\cdot C-(\vartheta^\ast)^{-1}c_{\varepsilon}'\log(C) \right].
\]
The last display yields \eqref{eq:GoalmkComp}, by choosing $C$ such that $C'\cdot C-(\vartheta^\ast)^{-1}c_{\varepsilon}'\log(C) \ge y$.
\end{proof}

\begin{proof}[\hypertarget{Proof:CorollaryMax}{Proof of Corollary \ref{Kor:Max}}]
We recall $\wt{m}_t := \inf\{k\in\N : m_{k+1}\ge t\}$.  The aim is to show that $(M_t-\wt{m}_t)_{t\ge0}$ is tight. Thus, we fix $\varepsilon>0$. 

 We start with a bound on the right tail of $M_t-\wt{m}_t$. Fix $t$ with $\wt{m}_t = k$ and $C\in\N$. We have
 \begin{equation}\label{eq:Mt-mtRT}
\begin{aligned}
P\left[M_t-\wt{m}_t\ge C \right] &= P\left[M_t \ge k+C\right] \le P\left[\mathbf{H}_{k+C}\le t\right] \le P\left[\mathbf{H}_{k+C} \le m_{k+1}\right]\\
&\le P\left[ \mathbf{H}_{k+C}-m_{k+C} \le m_{k+1}-m_{k+C}\right].
\end{aligned}
\end{equation}
By Theorem \ref{Theo:Main} there is a $y_\varepsilon>0$ such that
\begin{equation}\label{eq:HktightCor}
\limsup_{k\to\infty}P\left[\mathbf{H}_{k}-m_{k}\le -y_\varepsilon\right]\le \varepsilon/2.
\end{equation}
Furthermore, by Lemma \ref{Lem:mkComp} there is a $C_{\varepsilon}>0$ such that 
\begin{equation}\label{eq:mkBoundCor}
\limsup_{k\to\infty} \mathbb{P}\left[m_{k+1}-m_{k+C_{\varepsilon}} \le -y_{\varepsilon}\right]\le \varepsilon/2.
\end{equation}
Combining \eqref{eq:Mt-mtRT}, \eqref{eq:HktightCor}, \eqref{eq:mkBoundCor} yields that there is a $C_{\varepsilon}>0$ such that
\begin{equation}\label{eq:MtmtRT}
\limsup_{t\to\infty} P\left[M_t-\wt{m}_t\ge C_{\varepsilon}\right] \le \varepsilon.
\end{equation}

Next, we give a corresponding bound on the left tail of $M_t-\wt{m}_t$. We have
\begin{align*}
P\left[M_t-\wt{m}_t \le -2C \right] &\le P\left[\mathbf{H}_{k-C} \ge t  \right]+P\left[\mathbf{H}_{k-C} \le t, M_t \le k-2C\right]\\
&\le P\left[ \mathbf{H}_{k-C}-m_{k-C} \ge m_k-m_{k-C}\right]+P_{k-C}[ M_{t-H_{k-C}} \le k-2C]\\
&\le P\left[ \mathbf{H}_{k-C}-m_{k-C} \ge m_k-m_{k-C}\right]+\smash{\sup_{s\le t}}P[M_{s} \le -C],
\end{align*}%
where in the last step we used that the distribution of $x+M_k$ under $P_0$ is equal to the distribution of $M_k$ under $P_x$, since our environment is i.i.d.

The first summand in the last display can be controlled analogously to the bound on \eqref{eq:Mt-mtRT} which yields that there is a $C_{\varepsilon}'$ such that
\begin{equation} \label{eq:TightMT-mtLTFirstStep}
\limsup_{t\to\infty} P[M_t-\wt{m}_t \le -2C_{\varepsilon}' ] \le \varepsilon/2+\limsup_{t\to\infty} \sup_{s\le t} P[M_s\le -C_{\varepsilon}'].
\end{equation}
It remains to show that we can choose $C_{\varepsilon}'$ so that the second summand in the last display is smaller than $\varepsilon/2$. For this, we use that we can couple different environments so that the maximum $M_t$ is monotonically increasing in the environment, which yields 
\begin{equation} \label{eq:MoveToHom}
\sup_{s\le t} P[M_s\le -C_{\varepsilon}'] \le \sup_{s\le t} P_0^{\mathrm{ei}}[M_s\le -C_{\varepsilon}'],
\end{equation}
where we recall that $P_0^{\mathrm{ei}}$ is used to denote the law of the BRW with $\xi(x) = \mathrm{ei}$ for all $x\in\Z$, i.e.\@ the law of homogeneous branching random walk with branching rate $\mathrm{ei}$. The last two displays together with tightness for the maximum of homogeneous random walk around its mean imply that there is a $C_{\varepsilon}'$ so that
\begin{equation} \label{eq:MTmtLT}
\limsup_{t\to\infty}P[M_t-\wt{m}_t \le -2C_{\varepsilon}'] \le \varepsilon.
\end{equation}
We combine \eqref{eq:MtmtRT}, \eqref{eq:MTmtLT} to conclude that $(M_t-\wt{m}_t)_{t\ge 0}$ is tight.
\end{proof}
\section{Barrier computations} \label{Sec:BarrComps}
This section aims to prove the barrier computations used in Section~\ref{Sec:Mainbody}, i.e.\@ Lemmata~\ref{Lem:FMomSimplBarrComp},~\ref{Lem:FMomHardBarrComp},~\ref{Lem:SecMomBarrCompLB}--\ref{Lem:SecMomBarrCompComplicated}. In \cite{TimeInhKriech} we have proved statements of this form in the case that $(H_k)_{k\in\N}$ is a Brownian motion. Thus, it suffices to approximate the sequence $(H_k-E^{\zeta}[H_k])_{k\in\N}$ by a Brownian motion under the measure $P^{\zeta}$.   The main tool for this will be a version of the KMT-Theorem for sums of independent but not identically distributed random variables from \cite{Sakhanenko}.  We need a few additional definitions before we proceed.
\begin{Def}[Nice Functions]\label{Def:NiceFunction}
We call a family of functions $(h_n)_{n\in\N}$ nice if $h_n(0) = 0$ for all $n\in\N$ and
\begin{equation}
C_h :=1\vee \sup_{n\in\N}\sup_{k\le n} \sup_{0\le l\le n-k} \frac{\max\{|h_n(n-k-l)-h_n(n-k)|, |h_n(k+l)-h_n(k)|\}}{\sqrt{l}}<\infty. \label{Def:Ch}
\end{equation}
\end{Def}

\begin{Def}[Interpolation adapted to $(B_s^\zeta)_{s\ge0}$]
The interpolation between $x_1, x_2\in \R$, on time-frame $[t_1,t_2]$, $t_1,t_2\in\N$, which we call adapted to $B_s^\zeta$ is given by
\begin{equation}\label{eq:AdaptedInter}
\begin{aligned}
L_{t_1,t_2}^{x_1,x_2} : [t_1,t_2] &\to \R,\\
 s&\mapsto x_1\cdot\left(\frac{\sum_{k=\lceil s\rceil+1}^{t_2}\xi_k^2+(\lceil s\rceil-s)\xi_{\lceil s\rceil}^2}{\sum_{k=t_1+1}^{t_2}\xi_k^2} \right)+x_2\left(\frac{(s-\lfloor s\rfloor)\cdot \xi_{\lceil s\rceil}^2 +\sum_{k= t_1+1}^{\lfloor s\rfloor} \xi_k^2}{\sum_{k=t_1+1}^{t_2} \xi_k^2}\right).
\end{aligned}
\end{equation}
\end{Def}
\begin{Def}[Barrier probabilities, dyadic scale]\label{Def:KMTBasicDef}
We set $\bar{\tau}_j := \tau_j-E^{\zeta}[\tau_j]$, $\bar{H}_k = \sum_{j=1}^k \bar{\tau}_j$ and note that $W_k = (\vartheta^\ast)^{-1}K_k^\zeta-\sum_{j=1}^k E^{\zeta}[\tau_j]$, compare \eqref{Def:Wnzeta} and \cite[Eq.~(A.2)]{CernyDrewitz}. We also recall the definitions of $\chi_k^2 := \mathrm{Var}_{P^{\zeta}}[\bar{\tau}_j]$ and $\sigma_k^2 := \sum_{j=1}^k \chi_j^2$. We set $c_0 = 0$, $c_k := 2^{k+1}$, $k\ge 1$, $k_0 := k_0(n) := \sup\{k\in\N : c_k\le \min\{n^{1/6},n/2\}\}$. Furthermore, we set
\begin{align}
C_{\log} := 1\vee \sup_{k\in\N_0} \sup_{l\ge \log(1\vee c_k)^3} \frac{|W_{c_k+l}-W_{c_k}|}{\sqrt{l\sqrt{\log(l)}}}<\infty,  \label{eq:Clog}
\end{align}
which is $\mathbb{P}$-a.s.\@ finite, which can be seen similarly to \cite[Eq.\@ (5.25)]{CernyDrewitz}. 

For a family of functions $(h_n)_{n\in\N}$ we set $W_k^h := W_k+h_n(k)$. We set $B_k^\zeta := \sum_{j=1}^k b_j$ for $(b_j)_{j\in\N}$ independent with $b_j\sim \mathcal{N}(0,\chi_j^2)$ and recall the notation $B_{t_2;t_1}^\zeta := B_{t_2}^\zeta-B_{t_1}^\zeta$. Furthermore, we introduce the following shorthand notations to make displays more readable. Let $t_1,t_2\in\N$, $y\in\R$, $I\subseteq \R$, $h: \R\to \R$ and define
\begin{equation}
\begin{aligned}
p_{t_1,t_2}^h(y,I) &:= P^{\zeta}\big[\mathcal{B}_{\{t_1+1,\dots, t_2\},-(W_{\cdot}^h-W_{t_1}^h)}^{y,I}(\bar{H}_{\cdot}-\bar{H}_{t_1})\big],\\
p_{t_1,t_2}^{B,h}(y,I) &:= P^{\zeta}\big[\mathcal{B}_{\{t_1+1,\dots, t_2\},-(W_{\cdot}^h-W_{t_1}^h)}^{y,I}(B^\zeta_{\cdot}-B_{t_1}^\zeta)\big], \\
p_{t_1,t_2}^{B,h}(y|x) &:= P^{\zeta}\big[ \mathcal{B}_{\{t_1+1,\dots, t_2\}, -W^h_{\cdot;t_1}}^{y}(B_{\cdot;t_1}^\zeta)| B_{t_2;t_1}^\zeta= x \big]
\end{aligned}\label{eq:pathprobsDef}
\end{equation}
We note that for $y$ fixed, all of these are (probability) measures on $\R$ and use $p_{t_1,t_2}(y,z)$, $z\in\R$ to denote the corresponding density function (with respect to the Lebesgue density). We use $\widehat{p}_{t_1,t_2}^{h,g,(n)}$ for the time-reversed version of these probabilities, e.g.\@
\begin{equation}
\widehat{p}_{t_1,t_2}^{h,g,(n)}(y,I) := P^{\zeta}\big[\mathcal{B}_{\{t_1+1,\dots, t_2\}, -(W_{n-\cdot}^h-W_{n-t_1}^h+g(\cdot))}^{y,I}(\bar{H}_{n-\cdot}-\bar{H}_{n-t_2})\big]. \label{eq:pt1t2TimeReversed}
\end{equation}
The reason that we allow two functions as superscripts in the last definition, is that this allows us to specify one function, which will be reversed in time, namely $h$, and one, which will be added to the barrier as is, without being reversed, namely $g$.
\end{Def}

We now state the theorem we use to approximate our random walk by a Brownian motion.\footnote{The result originates from \cite{Sakhanenko}, it can also recovered from \cite[Corollary 2]{SakhanenkoTilde} by making there the choice $y_0 = \log(n)$ (note that this is not equal to our $y_0$ throughout this paper).  The condition \cite[Eq.~(1.2)]{Zaitsev} that there is a $\lambda_0>0$ such that for all natural numbers $k$ we have $E^{\zeta}[|\bar{\tau}_k^3e^{\lambda_0^{-1}|\bar{\tau}_k|}] \le \lambda_0 E[|\bar{\tau}_k|^2]$ can be checked using the H\"{o}lder-inequality as well as the facts that the moments of $\bar{\tau}_j$ w.r.t.\@ $P^{\zeta}$ are uniformly bounded away from $0$ and $\infty$ and that $E^{\zeta}[e^{\lambda \tau_i}]<\infty$ for all $\lambda\le |\overline{\eta}|$, see \cite[Eq.\@ (4.9), (A.1)]{CernyDrewitz}. }

\begin{Sa}[Theorem 1.1,  \cite{Zaitsev}]\label{Def:KMT}
We can construct $(B_k^\zeta)_{k\in\N}$ on the same probability space as $(\bar{\tau}_k)_{k\in\N}$ such that there are $C,c\ge 0$ for which for all $n\in\N$
\begin{align*}
P^{\zeta}\bigg[\exists_{k\le n}\ \bigg|\sum_{j=1}^k \bar{\tau}_j-B_k^\zeta\bigg| \ge \log(n)^{1.5}\bigg]\le Ce^{-c\log(n)^{1.5}}.
\end{align*}
\end{Sa}

\subsection[\texorpdfstring{Uniform lower bound on barrier probabilities on different timescales}{Uniform lower bound on barrier probabilities on different intervals} ]{Lower bound on $p_{c_k,c_{k+1}}^h(y,I)$ and $p_{c_k,c_{k+1}}^{B,h-g}(y,I)$}
We need lower bounds on several barrier probabilities, which need to be uniform over different time-scales. We outline one of the applications for these lower bounds next, which motivates \eqref{eq:UnifPolyLB} below. We recall that $c_k := 2^{k+1}$, $k_0(n) := \sup\{k\in\N : c_k\le \min\{n/2, n^{1/6}\}$.

When we apply Theorem~\ref{Def:KMT} we get error terms of order $e^{-c\log(n)^{1.5}}$. Thus, we need to show that these are negligible compared to the size of the barrier probabilities we work with. For this purpose, it suffices to prove a polynomial lower bound on $p_{c_k,c_{k+1}}^h(y,I)$. One subtlety is, that we only get $p_{c_k,c_{k+1}}^h(y,I) \ge C((W_l)_{l\in \{c_k+1,\dots, c_{k+1}\}}) c_k^{-\gamma( (W_l)_{l\in\{c_k+1,\dots, c_{k+1}\}})}$, i.e.\@ the constants in the lower bound we get depend on the (random) walk $W$ at times between $c_k+1,c_{k+1}$. Since we need to apply the lower bound for all $k\le k_0(n)$ simultaneously, we thus need to show that there are tight sequences of random variables $\mathcal{C}_n$, $\gamma_n$ such that uniformly for $k\le k_0(n)$ we have $p_{c_k,c_{k+1}}^h(y,I)\ge \mathcal{C}_n^{-1} c_k^{-\gamma_n}$. In Lemma~\ref{Lem:UnifCrudeLB} we state this in terms of the notation $\wt{O}_{\text{unif},g}^n(f(k,n))$ introduced in Section~\ref{Sec:Notation}. We recall that we call a family of functions $(h_n)_{n\in\N}$ nice if $C_h<\infty$ with $C_h$ as in \eqref{Def:Ch}.

\begin{Lem}\label{Lem:UnifCrudeLB}
For $k\in\N$ let $y_k>0$ be a positive real number and $I_k = [y_k+a_k,y_k+b_k]$, with $a_k\ge -y_k$, $|a_k|,|b_k|\le C'\sqrt{c_k\log(c_k)}$, $b_k-a_k\ge \log(c_k)^3$, $C'\in (0,\infty)$, be an interval. Let $(h_n)_{n\in\N}$ be nice. There is a sequence $(\gamma_n)_{n\in\N}$ such that (uniformly in $y_k, I_k$) $\gamma_n\le \wt{O}_n(1)$ and 
\begin{equation} \label{eq:UnifPolyLB}
(p_{c_k,c_{k+1}}^h(y_k,I_k) \wedge p_{c_k,c_{k+1}}^{B,h}(y_k,I_k))^{-1} \le \wt{O}_{\text{unif},k\le k_0(n)-1}^{n}(c_k^{\gamma_n+(C'-1)^2+6(C'-1)(C_{\log}+C_h)^2c_{\chi}^{-1}}),
\end{equation}
where the implied constant and $\gamma_n$ are allowed to depend on $(h_n)_{n\in\N}$.

 Furthermore, there is a tight sequence $(\widehat{\gamma}_n)_{n\in\N}$ such that for $C''>0$ fixed, $y>0$, $I_k = [y+a_k,y+b_k]$, $-y_k\le a_k$, $|a_k| ,|b_k|\le \max\{\sqrt{c_k\log(c_k)},C''\}$, $b-a\ge \log(c_k)^3$, we have
\begin{equation} \label{eq:PolyLBonp0n}
(p_{0,c_k}^h(y,I_k))^{-1}\le \wt{O}_{\text{unif},k\le k_0(n)}^n(c_k^{\widehat{\gamma}_n}),
\end{equation}
where the implicit constant in $\wt{O}_{\text{unif},k\le k_0(n)}^n$ may depend on $C''$, but given $C''$ is uniform over the choice of $I_k$. Additionally, we have 
\begin{equation}\label{eq:p0nlbFin}
p_{0,n}^h(y,I_n)^{-1}\le \wt{O}_n(n^{\widehat{\gamma}_n}),
\end{equation}
for $I_n = [y+a_n,b_n]$, $-y\le a_n$, $|a_n| ,|b_n|\le \sqrt{n\log(n)}$, $b-a\ge \log(n)^3$. 

Finally, there is a tight sequence $(\wt{\gamma}_n)_{n\in\N}$ such that uniformly in $k\le k_0(n)$, $y_k \in[c_k^{1/3},\sqrt{(n-c_k)}]$,  $l\in \{y_k^2,\dots, \lceil (n-c_k)/2\rceil\}$ we have
\begin{equation}\label{eq:PolyLBonctn}
p_{c_k,c_k+l}^{B,h}(y_k, [2C_h\sqrt{l},4(C_h-1)\sqrt{l}]) \ge l^{-\wt{\gamma}_n}.
\end{equation}
Analogous statements to \eqref{eq:UnifPolyLB}--\eqref{eq:PolyLBonctn} hold for the time-reversed probabilities $\widehat{p}_{c_k,c_{k+1}}^{h,(n)}(y_k,I_k)$, $\widehat{p}_{0,c_k}^{h,(n)}(y_k,I_k)$, and $\widehat{p}_{c_k,l}^{B,h,(n)}(y_k, [2C_h\sqrt{l},4(C_h-1)\sqrt{l}]$.
\end{Lem}
\begin{Bem}
The reason, that we make the $C'$ dependency of the exponent in the right hand side of \eqref{eq:UnifPolyLB} explicit, is that later in this paper we will define a constant $\wt{C}$ depending on $\gamma_n$ and then apply \eqref{eq:UnifPolyLB} for $C' = \wt{C}$, see \eqref{eq:RWCond1C2} and \eqref{eq:DemonstrateReasonInLB}. We have formulated the $C'$ dependency to vanish for $C'=1$ for convenience.

We also point out that \cite[Eq. (6.8)]{CernyDrewitz} has a similar polynomial lower bound on a barrier probability, in fact, our $W_n$ corresponds to their $R_n'$. We cannot just use their result, since we need some uniformity over different time scales. We still use some of their ideas to deal with the fact that the random walk $W$ does not have independent increments. 
\end{Bem} 
\begin{proof}[Proof of Lemma~\ref{Lem:UnifCrudeLB}]
We start with proving the lower bound for $p_{c_k,c_{k+1}}^h(y_k,I_k)$. Fix $n\in\N$ and $k\le k_0(n)$. We want to apply Theorem~\ref{Def:KMT} to compare to a Gaussian probability, however the KMT error of $\log(c_k)^{1.5}$ may dominate the starting height $y_k$, which would produce a bad bound. This is circumvented by adding an event that the random walk distances itself from the wall early similarly to the argument below Lemma~6.4 in \cite{MMComp}. For the event involving only the beginning of the random walk we will be able to use a much cruder bound. Explicitly,  we set $K_{c_k} =\lfloor \log(c_k)^6\rfloor$ and $A_{c_k} := \log(c_k)^2$. By monotonicity and the Markov property
\begin{align} \label{eq:SetupCrudeLB}
p_{c_k,c_{k+1}}^h(y_k,I_k) \ge p_{c_k,c_k+K_{c_k}}(y_k, [y_k+A_{c_k}, y_k+c_k^{1/3}]) \inf_{y\in [y_k+A_{c_k},y_k+c_k^{1/3}]} p_{c_k+K_{c_k},c_{k+1}}^h(y,I_k).
\end{align}
We first deal with $p_{c_k,c_k+K_{c_k}}^h(y_k, [y_k+A_{c_k},y_k+c_k^{1/3}])$. It turns out that the time horizon $\{c_k,\dots, c_k+K_{c_k}\}$ is still slightly too long to bound this naively. Thus, we iterate this construction once and set $\widetilde{K}_{c_k} := \lfloor \log(K_{c_k})^6\rfloor$ and $\widetilde{A}_{c_k} := \widetilde{K}_{c_k}^{1/3}$. Another application of the Markov property yields
\begin{multline*}
p_{c_k,c_k+K_{c_k}}^h(y_k,[y_k+A_{c_k},y_k+c_k^{1/3}]) \ge p_{c_k,c_k+\widetilde{K}_{c_k}}(y_k,[y_k+\widetilde{A}_{c_k},y_k+K_{c_k}^{1/3}])\\
\times \inf_{y\in [y_k+\widetilde{A}_{c_k}, y_k+K_{c_k}^{1/3}]} p_{c_k+\widetilde{K}_{c_k},c_K+K_{c_k}}^h(y,[y_k+A_{c_k},y_k+c_k^{1/3}]).
\end{multline*}
To control the first factor on the right-hand side of the last display we use the fact that both $W$ and $h$ have $\mathbb{P}$-a.s.\@ bounded increments, which means that for some constants $C_{\chi,h}, c_{\chi,h}>0$ we have
\begin{multline*}
p_{c_k,c_k+\widetilde{K}_{c_k}}(y_k,[y_k+\widetilde{A}_{c_k},y_k+K_{c_k}^{1/3}])\\
\begin{aligned}&\ge P^{\zeta}[\forall_{l\in \{c_k,\dots, c_k+\widetilde{K}_{c_k}\}}\, \bar{H}_{l;c_k}\ge (l-c_k)\cdot C_{\chi,h}, \bar{H}_{c_k+\widetilde{K}_{c_k};c_k} \in [\widetilde{A}_{c_k}+C_{\chi,h}\widetilde{K}_{c_k},K_{c_k}^{1/3}-C_{\chi,h}\widetilde{K}_{c_k}  ] ]\\
&\ge e^{-c_{\chi,h}\widetilde{K}_{c_k}}.
\end{aligned}
\end{multline*}
Here the last inequality is achieved by forcing $\bar{\tau}_i \in [2C_{\chi,h},3 C_{\chi,h}]$, $i\in\{c_k+1,\dots, c_k+\widetilde{K}_{c_k}\}$, which has an exponential in $\widetilde{K}_{c_k}$ cost and is possilbe since the $\bar{\tau}_i$ have support on all of $\R^+$. We note that by definition $e^{-c_{\chi,h}\widetilde{K}_{c_k}}$ decays slower than polynomially in $c_k$. 

It remains to control
\begin{align}
&\inf_{y\in [y_k+A_{c_k},y_k+c_k^{1/3}]} p_{c_k+K_{c_k},c_{k+1}}^h(y,I_k),\quad\text{and} \label{eq:CrudeBarrGoodKMTTerm}\\
&\inf_{y\in [y_k+\widetilde{A}_{c_k}, y_k+K_{c_k}^{1/3}]} p_{c_k+\widetilde{K}_{c_k},c_K+K_{c_k}}^h(y,[y_k+A_{c_k},y_k+c_k^{1/3}]).
\end{align}

Here we have ensured that the starting height is high enough, i.e., for the first line we have ensured that $y-\log(c_k)^{1.5}\gg 0$ while for the second we have ensured that $y-\log(K_{c_k}-\widetilde{K}_{c_k})^{1.5}\gg 0$. Thus, we can apply Theorem~\ref{Def:KMT}. We do this in some detail for \eqref{eq:CrudeBarrGoodKMTTerm}, the other term can be handled analogously. For \eqref{eq:CrudeBarrGoodKMTTerm} Theorem~\ref{Def:KMT} together with $c_{k+1}-c_k = c_k$ imply
\begin{multline} \label{eq:CrudeLBKMT}
\inf_{y\in [y_k+A_{c_k},y_k+c_k^{1/3}]} p_{c_k+K_{c_k}}^h(y,I_k) \\
\ge \inf_{y\in [y_k+A_{c_k},y_k+c_k^{1/3}]\atop x\in [y+a+\log(c_k)^{1.5}, y+b-\log(c_k)^{1.5}]}P^{\zeta}\left[\mathcal{B}_{\{c_k+K_{c_k},\dots,c_{k+1}\}, -W_{\cdot;c_k+K_{c_k}}^h}^{y-\log(c_k)^{1.5},J_x}(B_{\cdot;c_k}^\zeta)\right]-Ce^{-c\log(c_k)^{1.5}}.
\end{multline}
We note that $y_k+A_{c_k}-\log(c_k)^{1.5}>0$. Next, we want to argue similarly to \cite[Theorem 10.2]{TimeInhKriech}. For ease of notation we will provide a lower bound on 
\begin{equation}\label{eq:AltLB}
\inf_{y>0, x\in [y+a+2\log(c_k)^{1.5}, y+b]} P^\zeta\left[\mathcal{B}_{\{c_k,\dots, c_{k+1}\}, -W_{\cdot;c_k}^h}^{y,J_x}(B_{\cdot;c_k}^{\zeta})\right]-Ce^{-c\log(c_k)^{1.5}},
\end{equation} i.e., we won't properly keep track of the fact that we work on the time horizon $\{c_k+K_{c_k},\dots, c_{k+1}\}$ instead of $\{c_k,\dots, c_{k+1}\}$. Many of the steps from \cite{TimeInhKriech} need some adaptation to our setup, since our underlying Brownian motion has a time-dependent variance and our barrier $W$ doesn't have independent increments.  We spell out those adaptations here once, and later in similar situations only refer to the results from \cite{TimeInhKriech}.

 Set $k_1(c_k) := \lfloor \log_2(c_k/3)\rfloor$, $t_1(c_k) := c_k+2^{k_1}$, $t_2(c_k) := c_{k+1}-k_1(c_k)$, we omit the dependence on $c_k$ from the notation. Additionally fix a pair $(y,x)$ admissable for the infimum in \eqref{eq:AltLB}.  For $j\le k_1-1$ we define
\begin{align*}
\Xi_{j,k}^{\text{start}} &:= \max_{s\in [c_k,c_k\PlusUps 2^{j+1}]} W_{s;c_k}^h =\max_{s\in [c_k,c_k\PlusUps 2^{j+1}]} (W_s-W_{c_k}+h(s)-h(c_k)),\\
\Xi_{j,k}^{\text{end}} &:= \max_{s\in [c_{k+1}\MinusUps 2^{j+1}, c_{k+1}]} W_{s;c_k}^h-W_{c_{k+1};c_k}^h = \max_{s\in [c_{k+1}\MinusUps 2^{j+1}, c_{k+1}]} (-W_{c_{k+1};s}^h)\\
\Xi_{j,k} &:= \max\{\Xi_{j,k}^{\text{start}}, \Xi_{j,k}^{\text{end}}\},\quad \Xi_{k_1,k} := \max\{\max_{s\in[c_k,c_{k+1}]} W_s, \max_{s\in[c_k,c_{k+1}]} (W_{s;c_k}^h-W_{c_{k+1};c_k}^h)\}.
\end{align*}
We recall the definition of $L_{t_1,t_2}^{x_1,x_2}$ from~\eqref{eq:AdaptedInter} and define $g_k: [c_k,c_{k+1}]\to \R$ as follows
\begin{equation}
\begin{alignedat}{2}
g_k(s) &:=L_{c_k,c_k\PlusUps1 }^{0,\Xi_{0,k}}(s),&&\quad\text{for}\ s\in [c_k,c_k\PlusUps1],\\
g_k(s) &:=L_{c_k\PlusUps 2^j, c_k\PlusUps 2^{j+1}}^{\Xi_{j,k},\Xi_{j+1,k}}(s),&&\quad \text{for}\ s\in [c_k\PlusUps2^j,c_k\PlusUps2^{j+1}],\ j < k_1(c_k),\\
g_k(s) &:=L_{t_1,t_2}^{\Xi_{k_1,t}, \Xi_{k_1,t}+W_{c_{k+1};c_k}^h}(s),&&\quad\text{for}\ s\in[c_k+t_1,c_k+t_2],\\
g_k(s) &:=L_{c_{k+1}\MinusUps2^{j+1},c_{k+1}\MinusUps 2^j}^{\Xi_{j+1,k}+W_{c_{k+1};c_k}^h, \Xi_{j,k}+W_{c_{k+1};c_k}^h}(s),&&\quad\text{for}\  s\in[c_{k+1}\MinusUps2^{j+1}, c_{k+1}\MinusUps 2^j],\ j < k_1(c_k),\\
g_k(s) &:=L_{c_{k+1}\MinusUps1, c_{k+1}}^{\Xi_{0,k}+W_{c_{k+1};c_k}^h, W_{c_{k+1};c_k}^h}(s),&&\quad\text{for}\  s\in[c_{k+1}\MinusUps1,c_{k+1}].
\end{alignedat}
\end{equation}
By definition, $g_k(s)\ge W_{s;c_k}^h$ for all $s\in[c_k,c_{k+1}]$ and $g_k(c_{k+1}) = W_{c_{k+1};c_k}^h$, such that by monotonicity
\[
P^{\zeta}\big[\mathcal{B}_{\{c_k+1,\dots, c_{k+1}\}, -W_{\cdot;c_k}^h}^{y,J_x}(B_{\cdot;c_k}^\zeta)\big] \ge P^{\zeta}\big[\mathcal{B}_{\{c_k+1,\dots, c_{k+1}\}, -g_k(\cdot)}^{y,J_x}(B_{\cdot;c_k}^\zeta)\big].
\]
Next, we tilt $B_{s+1}^\zeta-B_s^\zeta$ so that $B_{s+1}^\zeta-B_s^\zeta-g_k(s+1)+g_k(s)$ is $\mathcal{N}(0,\chi_s^2)$ distributed, we introduce additional notation for this purpose. For $j\le k_1$ we set
\begin{equation}
\begin{alignedat}{3}
\Delta_{0,k}^{\text{start}} &:= \frac{\Xi_{0,k}}{\chi_{c_k+1}^2},&\quad \Delta_{0,k}^{\text{end}} &:= \frac{\Xi_{0,k}}{\chi_{c_{k+1}}^2},&\quad\Delta_{j,k}^{\text{start}} &:= \frac{\Xi_{j}-\Xi_{j-1,k}}{\sigma_{c_k\PlusUps2^{j};c_k\PlusUps2^{j-1}}^2},\\
\Delta_{j,k}^{\text{end}} &:= \frac{\Xi_{j,k}-\Xi_{j-1,k}}{\sigma_{c_{k+1}\MinusUps2^{j-1};c_{k+1}\MinusUps2^j}^2},&\quad\Delta^{\text{mid}}_k &:= \frac{W_{c_{k+1};c_k}^h}{\sigma_{t_2;t_1}^2}.
\end{alignedat}
\end{equation}
We also use the notation $g_k(s+1;s) := g_k(s+1)-g_k(s)$. The Radon-Nikodym derivative for the tilting which makes $B_{l+1}^\zeta-B_l^\zeta-g_k(l+1;l)$ into a mean zero normal random variable for all $l\in\{c_k+1,\dots, c_{k+1}-1\}$ is given by
\begin{align*}
Z_k &:= \frac{dQ^{\zeta}}{dP^{\zeta}} = \exp\bigg( \sum_{l=c_k}^{c_{k+1}-1} \frac{ B_{l+1;l}^\zeta g_k(l+1;l)}{\chi_{k+1}^2}-\frac{g_k(l+1;l)^2}{2\chi_{k+1}^2} \bigg).
\end{align*}
Next, we note that on the intervals $\{c_k\PlusUps2^j,\dots, c_k\PlusUps2^{j+1}-1\}$ respectively $\{c_{k+1}\MinusUps2^{j+1}, c_{k+1}\MinusUps2^j-1\}$ the quantity $\frac{g_k(l+1;l)}{\chi_{l+1}^2}$ is constant and equal to  $\Delta_{j+1,k}^{\text{start}}$ respectively $\Delta_{j+1,k}^{\text{end}}$ while  for $l = c_k$ respectively $l = c_{k+1}\MinusUps1$ it is equal to $\Delta_{0,k}^{\text{start}}$ respectively $\Delta_{0,k}^{\text{end}}$. Furthermore, for $l\in \{t_1,\dots, t_2-1\}$ it is equal to $\Delta_{k}^{\text{mid}}$. Finally, we drop the $g_k(l+1;l)^2$ terms, since we are after an upper bound. Thus, overall we get
\begin{equation}\label{eq:SeeAdaptedGood}
\begin{aligned}
Z_k &:= \frac{dQ^{\zeta}}{dP^{\zeta}} = \exp\bigg( \sum_{j=1}^{k_1} \left(\Delta_{j,k}^{\text{start}}(B_{c_k\PlusUps2^{j}}^\zeta-B_{c_k\PlusUps2^{j-1}}^{\zeta})+\Delta_{j,k}^{\text{end}}(B_{c_{k+1}\MinusUps2^{j-1}}^\zeta-B_{c_{k+1}\MinusUps2^j}^\zeta\right) \bigg)\\
&\qquad\cdot \exp\bigg(\Delta^{\text{mid}}_kB_{t_2;t_1}^\zeta+\Delta_{0,k}^{\text{start}}B_{c_k\PlusUps1;c_k}^\zeta+\Delta_{0,k}^{\text{end}}B_{c_{k+1};c_{k+1}-1}^\zeta \bigg).
\end{aligned}
\end{equation}
From here we can proceed just as in \cite[Theorem 10.2]{TimeInhKriech}, i.e.\@ we additionally force the Gaussian random walk $B_{c_k\PlusUps2^j;c_k}^\zeta$ to be roughly at height $y+g_k(c_k\PlusUps2^j)+\sigma_{c_k\PlusUps2^j;c_k}$ and $B_{c_{k+1};c_{k+1}\MinusUps2^j}$ to be roughly at height $x+\sigma_{c_{k+1};c_{k+1}\MinusUps2^j}+g_k(c_{k+1}\MinusUps2^j)$ both for $j\le k_1$. On the event that the Gaussian random walk behaves this way, we can control the random part of $Z_k^{-1}$ using the last display. To shorten the upcoming display we set 
\begin{equation}\label{eq:Ajstartend}
\begin{aligned}
 A_{j,\text{start}}(B^\zeta-g_k) &:= \left\{B_{c_k+2^j;c_k}^\zeta-g_k(c_k+2^j) \in [y+\sigma_{c_k+2^j;c_k}, y+2\sigma_{c_k+2^j;c_k}]\right\} \\
  A_{j,\text{end}}(B^\zeta-g_k) &:= \left\{B_{c_{k+1};c_{k+1}-2^j}- g_k(c_{k+1}-2^j) \in x+[\sigma_{c_{k+1};c_{k+1}-2^j}, 2\sigma_{c_{k+1};c_{k+1}-2^j}]\right\}.
\end{aligned}
\end{equation} Going through the procedure outlined above, using the control \eqref{eq:Ajstartend} gives on increments of $B_{\cdot}^\zeta$ and also keeping in mind~\eqref{Def:Ch},~\eqref{eq:Clog}, $|y-x|\le C'\sqrt{c_k\log(c_k)}+\log(c_k)^{1.5}$, $x\in [-\log(c_k)^{1.5}, (C'+1)\sqrt{c_k\log(c_k)}]$, $\sigma_{u;v}^2 \in [c_{\chi}(u-v), C_{\chi}(u-v)]$ as well as $\Xi_{0,k}^{\text{start}}, \Xi_{0,k}^{\text{end}}\in [c,C]$ for some universal positive constants, since $W$, $h$ have uniformly bounded increments, yields
\begin{align}
P^{\zeta}\big[\mathcal{B}_{\{c_k+1,\dots,c_{k+1}\}, -W_{\cdot;c_k}^h}^{y,J_x}(B_{\cdot;c_k}^\zeta)\big]&\ge Q^{\zeta}\bigg[\hskip-0.65pt \mathcal{B}_{\{c_{k}+1,\dots, c_{k+1}\}}^{y,J_x}(B_{\cdot;c_k}^\zeta-g_k(\cdot)), \bigcap_{j\le k_1} (A_{j,\text{start}}(B^\zeta-g_k)\cap A_{j,\text{end}}(B^\zeta-g_k)) \hskip-0.65pt\bigg]\notag
\\&\quad\times \exp\bigg(-4\sum_{j=1}^{k_1} \left((c_{\chi} 2^j)^{-1}(\Xi_{j}-\Xi_{j-1,k})^2+c_{\chi}^{-1}2^{-j/2+1}(\Xi_{j,k}-\Xi_{j-1,k})\right)\bigg) \notag\\
&\qquad\times\exp\bigg(-6C'(C_{\log}+C_h)^2c_{\chi}^{-1}\log(c_k)^{3/4}\bigg). \label{eq:ForceBMToBesqrtGood}
\end{align}
Next, we note that under $Q^\zeta$ the sequence $(B_{r;c_k}-g_k(r))_{r\in \{c_k+1,\dots, c_{k+1}\}}$ is a Gaussian random walk with variance structure inherited from $B^\zeta$, such that
\begin{equation}\label{eq:BrownianProbAfterTilt}
\begin{aligned}
&Q^{\zeta}\left[\mathcal{B}_{\{c_{k}+1,\dots, c_{k+1}\}}^{y,J_x}(B_{\cdot;c_k}^\zeta-g_k(\cdot)), \bigcap_{j\le k_1} (A_{j,\text{start}}(B^\zeta-g_k)\cap A_{j,\text{end}}(B^\zeta-g_k)) \right]\\
&= P^{\zeta}\left[\mathcal{B}_{\{c_{k}+1,\dots, c_{k+1}\}}^{y,J_x}(B_{\cdot;c_k}^\zeta), \bigcap_{j\le k_1} (A_{j,\text{start}}(B^\zeta)\cap A_{j,\text{end}}(B^\zeta)) \right] \ge c_k^{-\gamma-(C'^2)},
\end{aligned}
\end{equation}
for some $\gamma>0$ (independent of $c_k$, $W$, $y$, $x$), where the last step uses a ballot theorem for Brownian motion.\footnote{Since we force the Brownian motion to be at roughly height $\sqrt{k}$ this is slightly more general than the classical ballot theorem for Brownian motion. The version we use can be proved by hand, by using the Markov property of Brownian motion at times $c_k+2^j$, $c_{k+1}-2^j$, $j\le k_1$ to write the barrier probability as the product of roughly $2\log(c_k)$ uniformly lower bounded factors.} The $-(C')^2$ in the exponent corresponds to the cost for the Brownian motion to travel from $y$ to $x$, which can be bound by $c_k^{-(C')^2}$, since we have $|y-x|\le C'\sqrt{c_k\log(c_k)}+\log(c_k)^{1.5}$.

We combine \eqref{eq:CrudeLBKMT},~\eqref{eq:ForceBMToBesqrtGood} and \eqref{eq:BrownianProbAfterTilt} and do some minor further simplification to get
\begin{equation}
\label{eq:LowerBoundGotten}
\begin{aligned}
&p_{c_k,c_{k+1}}^h(y_k,I_k)\ge c_k^{-\gamma-(C')^2-6C'(C_{\log}+C_h)^2c_{\chi}^{-1}}\\
&\qquad\times c_k^{-\log(c_k)^{-1}C\sum_{j=1}^{k_1} 2^{-j}(\Xi_{j,k}-\Xi_{j-1,k})^2-\log(c_k)^{-1}C\sum_{j=1}^{k_1}2^{-j/2}|\Xi_{j,k}-\Xi_{j-1,k}|}-Ce^{-c\log(c_k)^{1.5}}.
\end{aligned}
\end{equation}

Thus, to conclude \eqref{eq:UnifPolyLB} we only need to prove that
\begin{equation}\label{eq:Tightness}
\sup_{k\le k_0(n)} \log(c_k)^{-1}\sum_{j=1}^{k_1}2^{-j}(\Xi_{j,k}-\Xi_{j-1,k})^2\hskip10pt\text{and}\hskip8pt \sup_{k\le k_0(n)} \log(c_k)^{-1}\sum_{j=1}^{k_1}2^{-j/2}|\Xi_{j,k}-\Xi_{j-1,k}|
\end{equation}
 are tight. We emphasize, that we are able to absorb the $Ce^{-c\log(c_k)^{1.5}}$ error-term in \eqref{eq:LowerBoundGotten} into the polynomial lower bound we prove on $p_{c_k,c_{k+1}}^h(y_k,I_k)$ for all but finitely many $k$. Call the set of $k$ for which this is impossible $I_{\text{bad}}$. We have that $|I_{\text{bad}}|$ is a continuous function of the exponent on the right hand side of \eqref{eq:LowerBoundGotten}, which we will prove to be tight in $n$. Assuming this tightness, and using that all increments $\bar{\tau}$ have a density with support on $\R_+$, we get
\[
\left(\min_{k\in I_{\text{bad}}} p_{c_k,c_{k+1}}^h(y_k, I_k)\right)^{-1} \le \wt{O}_n(1),
\]
which can be absorbed into the upper bound we are after.

We proceed with proving \eqref{eq:Tightness}, as a first step towards this we note that
\begin{equation}\label{eq:Upper}
\begin{aligned}
|\Xi_{j,k}-\Xi_{j-1,k}| &\le \max\big\{\sup_{s\in c_k+\{2^j,\dots, 2^{j+1}\}} |W_{s;c_k+2^j}^h|, \sup_{s\in\{c_{k+1}-2^{j+1},\dots, c_{k+1}-2^j\}} |W_{s;c_{k+1}-2^{j+1}}^h|\big\} \\
&=:\max\{a_{1,j,c_k}, a_{2,j,c_k}\}.
\end{aligned}
\end{equation}
Because $(h_n)_{n\in\N}$ is nice, we have 
\[
\sup_{s\in \{c_k+2^j,\dots, c_k+2^{j+1}\}} |W_{s;c_k+2^j}^h| \le \sup_{s\in \{c_k+2^j,\dots, c_k+2^{j+1}\}} |W_{s;c_k+2^j}|+ C_h2^{-j/2},
\]
so that we can assume from now on w.l.o.g.\@ that $h_n \equiv 0$.

Secondly, we set 
\begin{equation}\label{eq:EssSuppedGuys}
\widehat{a}_{1,j,c_k} := \esssup_{(\xi(x))_{x\le c_k+2^{j-1}}} a_{1,j,c_k}\qquad \text{and}\qquad \widehat{a}_{2,j,c_k} := \esssup_{(\xi(x))_{x\le c_{k+1}-2^{j+2}}} a_{2,j,c_k},
\end{equation} i.e.\@ for $\widehat{a}_{1,j,c_k}$ we take the essential supremum with respect to $(\xi(x))_{x\le c_k+2^{j-1}}$ and leave the remaining environment random. We note that by definition $a_{r,j,c_k}\le \widehat{a}_{r,j,c_k}$, $r\in\{1,2\}$. To establish \eqref{eq:Tightness} it suffices to show that for $r\in\{1,2\}$
\begin{equation}\label{eq:ImprovedTight}
\sup_{k\le k_0(n)} \log(c_k)^{-1} \sum_{j=1}^{k_1} \big(2^{-j/2}\widehat{a}_{r,j,c_k}+2^{-j}\widehat{a}_{r,j,c_k}^2\big)  \le \wt{O}_n(1).
\end{equation} This has the advantage, that $\widehat{a}_{1,j,c_k} \in \sigma\left(\xi(x) : x\in \{c_k+2^{j-1}+1,\dots,c_k+2^{j+1} \}\right)$ which implies that the sequence $(\widehat{a}_{1,j,c_k})_{j\le k_1}$ is $1$-dependent and similarly for $\widehat{a}_{2,j,c_k}$. To prove \eqref{eq:ImprovedTight} we need tail bounds on $\sum_{j=1}^{k_1} 2^{-j/2}\widehat{a}_{r,j,c_k}$.  As a first step towards this, we note that from \cite[Lemma A.2]{CernyDrewitz} we know that there is a constant $C_{\overline{\eta}} \in(0,\infty)$ such that for $r\in\{1,2\}$ $\mathbb{P}$-a.s.
\begin{equation}\label{eq:HatNonHatDistSmall}
\widehat{a}_{r,j,c_k}-a_{r,j,c_k} \le \sum_{r=2^j}^{2^{j+1}} C_{\overline{\eta}}e^{-(r-2^{j-1})/C_{\overline{\eta}}} \le 2^jC_{\overline{\eta}}e^{-2^{j-1}/C_{\overline{\eta}}} \le \wt{C}_{\overline{\eta}}
\end{equation}
is smaller than some universal constant. Furthermore, \cite[Lemma A.5]{CernyDrewitz} and \cite[Theorem 1]{10.3150/10-BEJ304} yield that there are constants $C,c>0$ such that uniformly in $r\in\{1,2\}$, $k\le k_0(n)$, $j\le k_1$
\begin{equation}\label{eq:MaxABound}
\mathbb{P}[2^{-j/2}a_{r,j,c_k}> x] \le Ce^{-cx^2}\quad\text{and}\quad\mathbb{P}[2^{-j}a_{r,j,c_k}^2> x] \le Ce^{-cx}.
\end{equation}
We combine \eqref{eq:HatNonHatDistSmall} and \eqref{eq:MaxABound} to get that there are (different) $C,c>0$ such that uniformly in $r\in\{1,2\}$, $k\le k_0(n)$, $j\le k_1$,
\begin{equation}
\mathbb{P}[2^{-j/2}\widehat{a}_{r,j,c_k}> x] \le Ce^{-cx^2}\quad\text{and}\quad\mathbb{P}[2^{-j}\widehat{a}_{r,j,c_k}^2> x] \le Ce^{-cx},
\end{equation}
which by integration implies that for $C_{\lambda} := Ce^{c\lambda^2}$ uniformly in $r\in\{1,2\}$, $k\le k_0(n)$, $j\le k_1$ 
\begin{align*}
\mathbb{E}[e^{\lambda 2^{-j/2}\widehat{a}_{r,j,c_k}}] \le C_{\lambda}\quad\text{and}\quad\mathbb{E}[e^{\lambda 2^{-j}\widehat{a}_{r,j,c_k}^2}] \le C_{\lambda}.
\end{align*}
Now we use a union bound to split the sums into $j$ even and $j$ odd, thus getting independent sums, and we apply the exponential Chebyshev inequality to conclude that
\[
\mathbb{P}\bigg[\log(c_k)^{-1} \sum_{j=1}^{k_1} \big(2^{-j/2}\widehat{a}_{r,j,c_k}+2^{-j}\widehat{a}_{r,j,c_k}^2\big) > x\bigg] \le 4e^{-\lambda x\log(c_k)}e^{C_{\lambda}\log(c_k)},
\]
and a union bound over $k$ gives \eqref{eq:ImprovedTight}, since $\lim\limits_{x\to\infty} \sup_{n\in\N}\smash[t]{\sum_{k=1}^{k_0(n)}} 4e^{-\lambda x\log(c_k)}e^{C_{\lambda}\log(c_k)} = 0$.

The inequalities \eqref{eq:PolyLBonp0n}--\eqref{eq:PolyLBonctn} can be proven analogously, with minor changes in the sequences of random variables which need to be tight in the last step of the proof.  For this reason, we do not copy the whole calculation but move directly to the analogues of \eqref{eq:EssSuppedGuys}. For \eqref{eq:PolyLBonp0n} we set
\begin{align*}
\widehat{b}_{1,j} &:= \esssup_{(\xi(x))_{x\le 2^{j-1}}}\sup_{s\in \{2^j,\dots,2^{j+1}\}} |W_{s;2^j}|\quad\text{and}\\
\widehat{b}_{2,j,k} &:= \esssup_{(\xi(x))_{x\le c_k-2^{j+2}}}\sup_{s\in \{c_k-2^{j+1},\dots, c_k-2^j\}} |W_{s;c_k-2^{j+1}}|
\end{align*}
and we need to show that
\begin{equation} \label{eq:UnifTIghtnessTwo}
\sup_{k\le k_0(n)} \log(c_k)^{-1}\sum_{j=1}^{\lfloor \log_2(c_k/3)\rfloor} 2^{-j/2}(\widehat{b}_{1,j}+\widehat{b}_{2,j,k})+\log(c_k)^{-1}\sum_{j=1}^{\lfloor \log_2(c_k/3)\rfloor}2^{-j}(\widehat{b}_{1,j}^2+\widehat{b}_{2,j,k}^2)\le \wt{O}_n(1).
\end{equation}
Exactly the same calculation as for $\widehat{a}_{r,j,c_k}$ yields that there are $\lambda', C_{\lambda}'>0$ such that
\begin{multline*}
\mathbb{P}\left[\log(c_k)^{-1}\sum_{j=1}^{\lfloor \log_2(c_k/3)\rfloor} 2^{-j/2}(\widehat{b}_{1,j}+\widehat{b}_{2,j,k})+\log(k)^{-1}\sum_{j=1}^{\lfloor \log_2(c_k/3)\rfloor}2^{-j}(\widehat{b}_{1,j}^2+\widehat{b}_{2,j,k}^2) \right] \\
\le Ce^{-\lambda' x\log(k)}e^{C_\lambda'\log(c_kk)},
\end{multline*}
so that a union bound yields \eqref{eq:UnifTIghtnessTwo}, since $\lim_{x\to\infty} \sup_{n\in\N} \sum_{k\le k_0(n)} Ce^{-\lambda' x\log(c_k)}e^{C_{\lambda}'\log(c_k)} = 0$ and a similar calculation yields \eqref{eq:p0nlbFin} .

For \eqref{eq:PolyLBonctn} we don't need $\Xi_{j,k}^{\text{end}}$ since we don't require the Brownian motion to end up at distance $O(1)$ to the barrier. An analogous argument to the one leading to \eqref{eq:LowerBoundGotten} allows us to establish
\begin{align*}
&p_{c_k,c_k+l}^{B,h}(y_k, [2C_h\sqrt{l},4(C_h-1)\sqrt{l}])\\
& \ge l^{-\widetilde{\gamma}-\log(l)^{-1}C\sum_{j=1}^{\log_2(l/3)} 2^{-j}\left(\Xi_{j,k}^{\text{start}}-\Xi_{j-1,k}^{\text{start}}\right)^2-\log(l)^{-1}C\sum_{j=1}^{\log_2(l/3)}2^{-j/2}\left|\Xi_{j,k}^{\text{start}}-\Xi_{j-1,k}^{\text{start}}\right|}.
\end{align*}
 Thus we define 
\[
\widehat{c}_{j,k} := \esssup_{(\xi(x))_{x\le c_k+2^{j-1}}} \sup_{s\in \{c_k+2^j,\dots, c_k+2^{j+1}\}} |W_{c_k;c_k+s}|
\] 
and to conclude \eqref{eq:PolyLBonctn} it suffices to show 
\begin{equation}\label{eq:LastTight}
\sup_{k\le k_0(n)} \sup_{l\in \{c_k^{2/3},\dots, \lceil (n-c_k)/2\rceil\}} \log(l)^{-1}\sum_{j=1}^{\log_2(l/3)} \left(2^{-j/2}\widehat{c}_{j,k}+2^{-j}\widehat{c}_{j,k}^2\right) \le \wt{O}_n(1).
\end{equation}
The same calculation as above yields that there are $\lambda'', C_{\lambda}''>0$ such that
\begin{align*}
\mathbb{P}\left[\log(l)^{-1}\sum_{j=1}^{\log_2(l/3)} \left(2^{-j/2}\widehat{c}_{j,k}+2^{-j}\widehat{c}_{j,k}^2\right) \ge x\right] \le e^{-\lambda'' x\log(l)}e^{C_{\lambda}''\log(l)},
\end{align*}
but this time after unit bounding the quantity which needs to decay as $x$ goes to infinity is
\begin{equation}\label{eq:SomeSum}
\sup_{n\in\N} \sum_{k=1}^{\lceil\log(n)\rceil} \sum_{l=c_k^{2/3}}^{\lceil(n-c_k)/2\rceil} e^{-\lambda'' x\log(l)}e^{C_{\lambda}''\log(k)}
\end{equation}
For $x>C_{\lambda}''/\lambda''$ we can upper bound \eqref{eq:SomeSum} by 
\[
 C\sum_{k=1}^\infty c_k^{2/3(C_{\lambda}''-\lambda'' x)+1}  = C \left(\frac{1}{1-2^{2/3(C_{\lambda}''-\lambda'' x)+1}}-1\right),
 \]
 which goes to $0$ as $x\to\infty$, yielding \eqref{eq:LastTight}.
\end{proof}

\subsection[\texorpdfstring{Comparison to Brownian barrier probabilities -- endpoint in an interval of logarithmic length}{Comparison to Brownian barrier probabilities -- endpoint in an interval of logarithmic length}]{Comparison to Brownian barrier probabilities -- endpoint in an interval of logarithmic length} \label{Sec:RWtoBMHard}

In this section, we give upper and lower bounds on $p_{0,n}^h(y,J)$ in terms of Brownian motion probabilities $p_{0,n}^{B,h\mp g}(y,\wt{J})$, where both $\wt{J}$ and $g$ will be made explicit further down. These bounds are achieved, by approximating the random walk $(\bar{H}_k)_{k\in\N}$ with a Brownian motion using Theorem~\ref{Def:KMT}. We emphasize that here due to the use of Theorem~\ref{Def:KMT} we have $|\wt{J}|\ge |J|+\log(n)^{1.5}$, and that we assume $|J|\ge C\log(n)^3$ so that the bounds we achieve will not be sufficient to immediately control the barrier probabilities in Lemmata~\ref{Lem:FMomSimplBarrComp},~\ref{Lem:SecMomBarrCompLB},~\ref{Lem:SecMomBarrCompUB}. The additional step needed to get a bound which can be used in the proofs of these lemmata is done in Section~\ref{Sec:BrownComparCompact}.

Next, we define the accumulated KMT error, which is
\begin{equation}\label{eq:gknjPre}
g_{k,\text{pre}}^{(n)}(j) := \sum_{l=k}^{k_0-1} \left(\mathbf{1}_{j\ge c_l+1} \log(c_l\vee 1)^{1.5}\right)+\mathbf{1}_{j\ge c_{k_0}+1}\log(n-c_{k_0})^{1.5},\quad j\in\N,
\end{equation} where the superscript $(n)$ is omitted most of the times. To be able to apply Lemma~\ref{Lem:UnifCrudeLB} we need a nice version of $g_{k,\text{pre}}^{(n)}(j)$, which is given by a linear interpolation, see \eqref{eq:AdaptedInter},
\begin{equation}\label{eq:gknj}
g_{k}^{(n)}(j) := L_{c_r+1,c_{r+1}+1}^{g_{k,\text{pre}}^{(n)}(c_r+1),g_{k,\text{pre}}^{(n)}(c_{r+1}+1)}(l),\quad \text{for}\ l\in \{c_r+1,\dots, c_{r+1}+1\}, r\le k_0,
\end{equation}
and $g_k^{(n)}(j) = g_k^{(n)}(c_{k_0}+1)$ for $j\ge c_{k_0}+1$. We note that $g_k^{(n)}(j)\ge g_{k,\text{pre}}^{(n)}(j)$ by definition.

 We recall that $y_0\ge e+1$ is some fixed constant.

\begin{Sa} \label{Sa:IndArgUpperBound}
Fix $C>0$, $(h_n)_{n\in\N}$ nice. There is a family of random variables $C_{k,n}(\zeta)>0$ so that the sequence $\left(\sup_{k\le k_0} \prod_{r=k}^{k_0} C_{r,n}(\zeta)\right)_{n\in\N}$ is tight and for $n$ big enough (not depending on $W$), all $y\in[y_0,\log(n)^2]\cap \N$, $k\in \{0,\dots, k_0\}$, $J = [a,b]$ with $0\le a<C\sqrt{n\log(n)}$, $|b-a|>\log(n)^3$ and $\lambda>0$
\begin{equation}
\begin{aligned}
p_{0,n}^h(y,J) \le \prod_{r=k}^{k_0} C_{r,n}(\zeta)\int_0^\infty p_{0,c_k}^h(y,z) p_{c_k,n}^{B,h-g_k}(z,\wt{J})\;\mathrm{d}z+\wt{O}_n((k_0-k)\cdot e^{-\lambda\log(n)}),
\end{aligned}
 \label{eq:UBBrownInd}
\end{equation}
where $\wt{J} := [a,b+2\log(n-c_{k_0})^{1.5}]$.
\end{Sa}

\begin{Sa}\label{Sa:IndArgLowerBound}
Fix $C>0$, $(h_n)_{n\in\N}$ nice. There is a family of random variables $\wt{C}_{k,n}(\zeta) \in (0,1]$ so that the sequence $\big(\big(\inf_{k\le k_0}\prod_{r=k}^{k_0} \wt{C}_{r,n}(\zeta)\big)^{-1}\big)_{n\in\N}$ is tight and for $n$ big enough (not depending on $W$), all $y\in[y_0,\log(n)^2]\cap \N$, $k\in \{0,\dots, k_0\}$, $K = [a,b]$ with $0\le a<C\sqrt{n\log(n)}$, $|b-a|>\log(n)^3$, and all $\lambda>0$
\begin{equation}
p_{0,n}^h(y,K)\ge \prod_{r=k}^{k_0} \wt{C}_{r,n}(\zeta)\int_0^\infty p_{0,c_k}^h(y,z)p_{c_k,n}^{B,h+g_k}(z,\wt{K})\;\mathrm{d}z-\wt{O}_n((k_0-k)e^{-\lambda \log(n)}), \label{eq:LBBrownInd}
\end{equation}
where $\wt{K} := [a, b-2\log(n-c_{k_0})^{1.5}]$. 
\end{Sa}

\begin{Bem}
\begin{itemize}
\item The idea of both proofs is to first replace the random walk $(\bar{H}_k)_{k\in\N}$ by a Brownian motion on $[c_{k_0},n]$, where the KMT error produces the blow-up of the interval to $\wt{J}$. Then we inductively apply KMT on $[c_{k-1},c_k]$. This will produce an error of order $\log(c_{k-1})^{1.5}$ on the starting position of the barrier event on $[c_k,n]$, however there we have already moved to a Brownian motion such that we can use the results on moving the starting point from \cite{TimeInhKriech}. The main technical problem is that we can only apply the KMT if $|\bar{H}_{c_{k-1}}-\bar{H}_{c_k}| \le C\sqrt{c_{k-1}\log(c_{k-1})}$, since otherwise the error probability in the KMT is bigger than the probability that $\bar{H}_\cdot$ connects these points. Thus, the proof has two steps, showing that $\bar{H}_{c_k}$ not being at roughly $y+\sqrt{c_k}$ is unlikely, and applying KMT to move to the Brownian probability in the case $\bar{H}_{c_k}\approx y+\sqrt{c_k}$.%
\item We mention that for $k=0$ the upper and the lower bound provided by Theorems~\ref{Sa:IndArgUpperBound},~\ref{Sa:IndArgLowerBound} match up to a shift, which is a logarithmic power, in the barrier and a change of order $\log(n)^{1.5}$ in the length of the ending interval $\wt{J}$. However, for $|\wt{J}|\gg\log(n)^{1.5}$ we know that these changes are insignificant for the barrier event in the case of Brownian motions, such that this gives upper and lower bounds, which are within a constant factor of each other.%
\item The proofs go through as well for $\widehat{p}_{0,n}^{h,(n)}(y,I) := P^{\zeta}\big[\mathcal{B}_{\{1,\dots, n\}, -(W_{n-\cdot}^h-W_n^h) }^{y,I}(\bar{H}_{n-\cdot}-\bar{H}_{n})\big]$, i.e.\@ the time-reversed analogue of the probability we care for, after redefining all constants depending on the barrier to depend on the time-reversed version of the barrier. We note that the Gaussian term on the right hand side will be $\widehat{p}_{c_k,n}^{B,h,\mp \widehat{g}_k,(n)}(z,\widetilde{J})$, with $\widehat{g}_k$ given by 
\begin{equation}\label{eq:gknjReversed}
\widehat{g}_{k}^{(n)}(j) := L_{n-c_{r+1}-1, n-c_r-1}^{g_{k,\text{pre}}^{(n)}(c_{r+1}+1),g_{k,\text{pre}}^{(n)}(c_{r+1}+1)}(l),\quad \text{for}\ l\in \{c_r+1,\dots, c_{r+1}+1\}, r\le k_0,
\end{equation}
and $\widehat{g}_k^{(n)}(j) = g_k^{(n)}(c_{k_0}+1)$ for $j\ge c_{k_0}+1$.
\end{itemize}
\end{Bem}
\subsubsection{Technical prerequisites}
Next, we prepare the proofs of Theorems~\ref{Sa:IndArgUpperBound},~\ref{Sa:IndArgLowerBound}.   For this, we need to understand how $p_{c_k,n}^{B,h\pm g_k}(z,\wt{J})$ varies in $z$. One subtlety in this is that like in \cite{TimeInhKriech} the bounds will depend on the barrier $W$ near $c_k$. We need this dependence to be uniform in $k$.  As a first step, we need to extend the result from \cite[Th.s~8.2,~8.4]{TimeInhKriech} that $p_{c_k,n}^{B,h\pm g_k}(z,\wt{J})$ grows polynomially in $z$ to the case that $\{z\}\cup\wt{J} \not \subseteq [0,\log(n)]$. In a second step we show that for $z,w \ge c_k^{1/3}$ with $|z-w|\le \log(c_k)^3$ we have $p_{c_k,n}^{B,h\pm g_k}(z,\wt{J}) \approx p_{c_k,n}^{B,h\pm g_k}(w,\wt{J})$.

\begin{Lem}\label{Lem:StartMove}
Fix $n\in\N$ big enough (independent of $W$), $k\le k_0(n)$, $0\le a<b$, $a<C\sqrt{n\log(n)}$, $J =  [a,b]$. Let $n^{1/6}>z\ge c_k^{1/3}$. Then for all $\lambda>0$ there is a tight sequence of random variables $(C_2(W,h,n))_{n\in\N}$ and a constant $C(\lambda)>0$ (which importantly both do not depend on $k$, $z$ nor $J$) such that
\begin{align*}
\frac{p_{c_k,n}^{B,h\pm g_k}(z,J)}{\inf_{x\in [1,z]} p_{c_k,n}^{B,h\pm g_k}(x,J)} \le  C_2(W,h,n) z^{C_2(W,h,n)}+C(\lambda)e^{-\lambda \log(n)}\big(\inf_{x\in [1,z]}p_{c_k,n}^{B,h\pm g_k}(x)\big)^{-1}.
\end{align*}
\end{Lem}
\begin{Bem}
Neither of the restrictions $z\ge c_k^{1/3}$, $z<n^{1/6}$ is optimal. The restriction $z\ge c_k^{1/3}$ is in place to get some control on e.g.\@ $\sup_{k\le k_0(n)} W_{c_k+z^2;c_k}$. In particular some lower bound on $z$ is needed to get a statement, which is uniform in $k$. On the other hand the statement should be true for $z$ all the way up to $\infty$. For $z\in [n^{\varepsilon}, \sqrt{n\log(n)}]$ it is true by upper bounding the numerator by $1$ and using Lemma~\ref{Lem:UnifCrudeLB} to get a polynomial in $n$ (and thus in $z$) lower bound on the denominator. For $z>\sqrt{n\log(n)}$ just the cost of going from $z$ to $J$ should force the infimum in the denominator to be taken near $z$, again giving the statement. However, for our applications $c_k^{1/3}\le z<n^{1/6}$ suffices and thus we don't pursue full generality.

It is worth pointing out that we will also apply Lemma~\ref{Lem:StartMove} in the time-reversed situation, where the Brownian motion starts at time $n-c_k$ and runs (backwards in time) until time $0$. In that situation the analogue of $C_{\log}$ will depend on $n$ and won't be uniformly bounded in $n$, which is why in Lemma~\ref{Lem:StartMove} we only claim that the sequence of coefficients $C_2(W,h,n)$ is tight and not almost surely bounded. 
\end{Bem}

\begin{proof}[Proof of Lemma~\ref{Lem:StartMove}]
For $z\le 2\log(n)^2$ we want to reduce this to \cite[Th.s 8.2, 8.4]{TimeInhKriech}. These theorems are only stated for ending in intervals of length 1, however going through the argument one sees that this restriction is arbitrary.\footnote{The only place in which $J$ matters is in the proof of \cite[Theorem 8.2]{TimeInhKriech}. In their equation (8.12) one would need to replace the integral over an interval of length 1 by an integral over $J$, but since in \cite[Lemma 8.6]{TimeInhKriech} lower bound does not depend on $z$ at all, this produces no trouble. (Note that their $z$ is not the same as our $z$, theirs corresponds to the precise location the Brownian motion ends up in, while ours is the starting point).}  The uniform tightness of the involved sequences of constants, follows after noting the following three things:
\begin{enumerate}[(a)]
\item The analogue to $\gamma$ from \cite[Eq.~(8.4)]{TimeInhKriech} is given by
\begin{equation}\label{eq:GammaEquiv}
\sup_{k\le k_0(n)} \inf\left\{\gamma\ge 0 : \forall_{z\le \log(n)^2} \forall_{s\in \{c_k+z^2,\dots, \lceil (n-c_k)/2\rceil\}} P^{\zeta}\left[\mathcal{B}_{\{c_k+1,\dots, s\},-W_{\cdot;c_k}^{h}\pm g_k(\cdot;c_k)}^{z,[2(C_h +C_{g_k})\sqrt{s}, 4(C_h+C_{g_k}-1)\sqrt{s}]}(B_{\cdot}^\zeta) \right] \ge s^{-\gamma}\right\}.
\end{equation}
The supremum over $k$ corresponds to the fact, that Lemma \ref{Lem:StartMove} makes a statement uniformly over different time-scales. Tightness of the expression in \eqref{eq:GammaEquiv} is ensured by \eqref{eq:PolyLBonctn} from Lemma~\ref{Lem:UnifCrudeLB}.
\item To apply the arguments of \cite[Th. 8.2--8.4]{TimeInhKriech} to $p_{c_k,n}^{B,h\pm g_k}(z,J)$ with $z\ge c_k^{1/3}$ one can replace \cite[Eq.~(8.2)]{TimeInhKriech}~by 
\[
\wt{\mathcal{C}}_{\log}(W,[c_k,n]) = \max\left\{1, \sup_{s\ge c_k^{2/3}} \frac{|W_{c_k+s;c_k}|}{\sqrt{s\log(s)}}\right\}. 
\]
The reason for this is, that in the proofs of \cite[Th. 8.2--8.4]{TimeInhKriech} the barrier $(W_{s+c_k;c_k})_{s\ge0}$ is always dropped for times of size at least $z^2\ge c_k^{2/3}$ (where $z$ is the starting height of the random walk, which in \cite{TimeInhKriech} is called $y$) so that we only need to control its behaviour for $s\ge c_k^{2/3}$. However the sequence  $
\left(\sup_{k\le k_0(n)} \wt{\mathcal{C}}_{\log(W,[c_k,n])}\right)_{n\in\N}
$
is a tight sequence of random variables, which can be seen similar to the finiteness of $C_{\log}$ from \eqref{eq:Clog}.
\item The constant depending on the non-random part of the barrier from \cite[Eq. (8.3)]{TimeInhKriech}, which is used in the bounds on $p_{c_k,n}^{B,h\pm g_k}(z,J)$ is given by
\begin{align}
&\max\left\{1, \sup_{n\in\N} \max_{l\in \{c_k+1,\dots, n\}} \frac{\max\{|h_n(l)-h_n(c_k)|, |h_n(n)-h_n(n-l)|\}}{\sqrt{l}}\right\}\label{eq:hPartBarr}\\ 
&\qquad+\max\left\{1, \sup_{n\in\N} \max_{l\in \{c_k+1,\dots, n\}} \frac{\max\{|g_k^{(n)}(l)-g_k^{(n)}(c_k)|, |g_k^{(n)}(n)-g_k^{(n)}(n-l)|\}}{\sqrt{l}}\right\}. \label{eq:gPartBarr}
\end{align}
The contribution of \eqref{eq:hPartBarr} is controlled by our assumption on $(h_n)_{n\in\N}$ that $C_h$ from \eqref{Def:Ch} is bounded. The contribution of \eqref{eq:gPartBarr} can be seen to be finite directly from the definition of $g_k^{(n)}(j)$ in \eqref{eq:gknj}.
\end{enumerate}

Thus we have dealt with the case that $z\le 2\log(n)^2$. 

Next, we handle $z\in [2\log(n)^2, n^{1/6}]$. In this case, we drop the barrier for time $z^3$ and use the Markov-property at time $c_k+z^3$ to get
\begin{align}
p_{c_k,n}^{B,h\pm g_k}(v,J) &\le \int_0^\infty g_{z-W_{c_k+z^3;c_k}^{h\pm g_k},\sigma_{c_k+z^3;z}^2}(r)\cdot p_{c_k+z^3,n}^{B,h\pm g_k}(r,J)\;\mathrm{d}r. \label{eq:UppBoundStartPointBigz}
\end{align}
We note that dropping the barrier for time $z^3$ is what forces us to consider $z\le n^{1/6}$ to be able to use \cite[Theorem 8.2]{TimeInhKriech}. Combining the argument leading to \cite[Theorem 8.2]{TimeInhKriech} and Lemma~\ref{Lem:UnifCrudeLB} yields that there is a tight sequence of random variables $(\gamma_n')_{n\in\N}$ (independent of $k$, $z$, $J$), such that 
\begin{equation}
\inf_{x\in [1,z]} p_{c_k,n}^{B,h-g_k}(x,J) \ge z^{-3\gamma_n'}\inf_{x\in[1,z]}\int_0^\infty g_{x-W_{c_k+z^3;c_k}^{h\pm g_k},\sigma^2_{c_k+z^3;c_k}}(r)\cdot p_{c_k+z^3,n}^{B,h\pm g_k}(r,J)\;\mathrm{d}r. \label{eq:LowBoundStartPointBigz}
\end{equation}
A direct calculation using the explicit form of $g_k$ and the bounds \eqref{Def:Ch}, \eqref{eq:Clog} yields that for $r\le z^2$, $x\le z$ we have
\begin{equation}
g_{z-W_{c_k+z^3;c_k}^{h\pm g_k},\sigma_{c_k+z^3;z}^2}(r)^{-1}g_{x-W_{c_k+z^3;c_k}^{h\pm g_k},\sigma^2_{c_k+z^3;c_k}}(r) \le C z^c, \label{eq:DensComp}
\end{equation}
where $C$ and $c$ only depend on $C_{\log}$, $C_h$ from \eqref{Def:Ch}, \eqref{eq:Clog}. For $r\ge z^2$ we can bound $p_{c_k+z^3,n}^{B,h\pm g_k}(r,\wt{J})$ by $1$ to get that
\begin{align}
\int_{z^2}^\infty g_{z-W_{c_k+z^3;c_k}^{h\pm g_k},\sigma_{c_k+z^3;z}^2}(r)p_{c_k+z^3,n}^{B,h\pm g_k}(r,\wt{J})\;\mathrm{d}r \le Ce^{-cz}, \label{eq:HighRegionIsSmall}
\end{align}
where again $C$ and $c$ are functions of $C_{\log}$, $C_h$. Since we have $z\ge 2\log(n)^2$ in particular we have that there is a $C(\lambda)>0$ so that $Ce^{-cz} \le C(\lambda) e^{-\lambda \log(n)}$ for all $n\in\N$. Combining \eqref{eq:UppBoundStartPointBigz}, \eqref{eq:LowBoundStartPointBigz}, \eqref{eq:DensComp} and \eqref{eq:HighRegionIsSmall} yields the claim of the Lemma for $z\in[ 2\log(n)^2, n^{1/6}]$.
\end{proof}

\begin{Lem} \label{Lem:StartPert} 
Fix $n\in\N$ big enough (independent of $W$), $k\le k_0(n)$, $0\le a<b$, $a<C\sqrt{n\log(n)}$, $\wt{J} = [a ,b\mp 2\log(n-c_{k_0})^{1.5}]$. Fix $n^{1/11}> z,w>c_k^{1/3}$ with $|z-w|\le 2\log(c_k)^3$. Then we have that for all $\lambda>0$
\begin{equation}
\begin{aligned}
p_{c_k,n}^{B,h\pm g_k}(w,\wt{J})\big(1-\wt{O}_{\text{unif}, k\le k_0(n)}(\log(c_k)^{-2})\big)-\wt{O}_n\big(e^{-\lambda\log(n)}\big)&\le p_{c_k,n}^{B,h\pm g_k}(z,\wt{J}),\\
p_{c_k,n}^{B,h\pm g_k}(w,\wt{J})\big(1+\wt{O}_{\text{unif},k\le k_0(n)}^n(\log(c_k)^{-2})\big)+\wt{O}_n\big(e^{-\lambda\log(n)}\big)&\ge p_{c_k,n}^{B,h\pm g_k}(z,\wt{J}).
\end{aligned}
 \label{eq:smallStartChange}
\end{equation}
\end{Lem}
The heuristic for Lemma~\ref{Lem:StartPert} being true is that the barrier will only start to matter from time $c_k+\min\{z,w\}^2$ onwards, such that it can be dropped until time $c_k+z^{1.9}$. However, a Brownian motion starting from $z$ and one starting from $w$ at time $c_k$ will have similar densities (near $z$) by time $c_k+z^{1.9}$.
\begin{Bem}
The upper bound $z,w<n^{1/11}$ is not optimal, we expect the statement to hold up to roughly height $\sqrt{n}$. We don't think that our method of proof is applicable up to $\sqrt{n}$ though, and $z,w<n^{1/11}$ is good enough for our applications.
\end{Bem}
\begin{proof}[Proof of Lemma~\ref{Lem:StartPert}]
Since $z,w$ are interchangeable it is enough to prove the first inequality in \eqref{eq:smallStartChange}. We recall that $\sigma_l^2 = \sum_{j=1}^l \Var_{P^{\zeta}}\left[\bar{\tau}_j\right]$,  $\sigma_{r;c_k}^2 := \sigma_l^2-\sigma_{c_k}^2$ and $\sigma_{u;v}^2\in [c_\chi(u-v),C_\chi(u-v)]$. We set $\wt{z} :=z^{1.9}/C_{\log}^2$. First, we give a lower bound on $p_{c_k,n}^{B,h\pm g_k}(w,\wt{J})$. Similarly to the proof of \cite[Th.~8.2.]{TimeInhKriech} we see that there is some $\eta>0$ small enough such that
\begin{multline*}
p_{c_k,n}^{B,h\pm g_k}(w,\wt{J}) \\
 \ge (1-\widetilde{O}_n(n^{-\eta}))P^{\zeta}\left[\mathcal{B}_{\{c_k,\dots,c_k+\widetilde{z}\}, -W_{c_k+\cdot;c_k}^{h\pm g_k}}^{w,[w-\sqrt{\widetilde{z}}\log(\widetilde{z}),w+\sqrt{\widetilde{z}}\log(\widetilde{z})]}(B_{\cdot;c_k}^\zeta) \right]\int_0^\infty g_{w-W_{c_k+\wt{z};c_k}^{h\pm g_k},\sigma_{c_k+\wt{z};c_k}^2}(r)p_{c_k+\wt{z},n}^{B,h\pm g_k}(r,\wt{J})\;\mathrm{d}r. 
\end{multline*}
Next, we control the middle probability in the last term. Using tail-bounds on the minimum of Brownian motion to control the probability that $B_{\cdot;c_k}^\zeta$ breaches the lower barrier\footnote{In this step we use that $\widetilde{z} = z^{1.9}/C_{\log}^2 \ll z-2\log(c_k)^3\le w$, i.e.\@ that the Brownian motion travels by less than $w$ by time $\widetilde{z}$ and $W_{c_k+\widetilde{z},c_k} \le \sqrt{z^{1.9}}\log(z^{1.9}/C_{\log}^2)$, i.e.\@ that the barrier won't move higher than our starting point by time $c_k+\widetilde{z}$.} while also  keeping in mind \eqref{Def:Ch}, \eqref{eq:Clog} and \eqref{eq:gknj} yields that
\[
P^{\zeta}\left[\mathcal{B}_{\{c_k,\dots,c_k+\widetilde{z}\}, -W_{c_k+\cdot;c_k}^{h\pm g_k}}^{w,[w-\sqrt{\widetilde{z}}\log(\widetilde{z}),w+\sqrt{\widetilde{z}}\log(\widetilde{z})]}(B_{\cdot;c_k}^\zeta) \right] \ge 1-\widetilde{O}_{\text{unif},k\le k_0(n)}(\log(c_k)^{-2}),
\]
which together with $1-\widetilde{O}_n(n^{-\eta})\ge 1-\widetilde{O}_{\text{unif},k\le k_0(n)}(\log(c_k)^{-2})$ yields that overall
\begin{equation}\label{eq:PertFKG}
p_{c_k,n}^{B,h\pm g_k}(w,\widetilde{J}) \ge (1-\widetilde{O}_{\text{unif}, k\le k_0(n)}(\log(c_k)^{-2}))\int_0^\infty g_{w-W_{c_k+\wt{z};c_k}^{h\pm g_k},\sigma_{c_k+\wt{z};c_k}^2}(r)p_{c_k+\wt{z},n}^{B,h\pm g_k}(r,\wt{J})\;\mathrm{d}r. 
\end{equation}
We proceed with an upper bound on $p_{c_k,n}^{B,h\pm g_k}(z,\widetilde{J})$. For this, we drop the barrier on $\{c_k,\dots, c_k+\wt{z}\}$ and use the Markov-property at time $c_k+\wt{z}$ to get
\begin{equation}
\begin{aligned}
p_{c_k,n}^{B,h\pm g_k}(z,\wt{J}) &\le \int_0^\infty g_{z-W_{c_k+\wt{z};c_k}^{h\pm g_k},\sigma_{c_k+\wt{z};c_k}^2}(r)p_{c_k+\wt{z},n}^{B,h\pm g_k}(r,\wt{J})\;\mathrm{d}r\\
&= \int_0^\infty e^{\frac{-2r(w-z)-2(w-z)W_{c_k+\wt{z};c_k}^{h\pm g_k}+w^2-z^2}{2\sigma_{c_k+\wt{z};c_k}^2}} g_{w-W_{c_k+\wt{z};c_k}^{h\pm g_k},\sigma_{c_k+\wt{z};c_k}^2}(r)p_{c_k+\wt{z},n}^{B,h\pm g_k}(r,\wt{J})\;\mathrm{d}r\\
&\stackrel{\eqref{eq:Clog}}{\le} \int_0^\infty e^{\frac{2r\log(c_k)^3+2(C_{\log}+C_h)\log(c_k)^4z^{0.95}+(w-z)(w+z)}{\sigma_{c_k+\wt{z};c_k}^2}}g_{w-W_{c_k+\wt{z};c_k}^{h\pm g_k},\sigma_{c_k+\wt{z};c_k}^2}(r)p_{c_k+\wt{z},n}^{B,h\pm g_k}(r,\wt{J})\;\mathrm{d}r. \label{eq:PertUBMark}
\end{aligned}
\end{equation}
Now we split into $r \le z^{1.8}$, $r\in[z^{1.8},n^{1/6}]$, $r>n^{1/6}$. In the region $r\le z^{1.8}$ we have
\begin{equation}
e^{\frac{2r\log(c_k)^3+2(C_{\log}+C_h)\log(c_k)^4c_k^{3/10}+(w-z)(w+z)}{\sigma_{c_k+\wt{z};c_k}^2}} \le e^{\frac{C(W,h,n)}{c_k^{1/30}}}\le 1+\wt{O}_{\text{unif}, k\le k_0(n)}^n(\log(c_k)^{-2}). \label{eq:PertUBLow}
\end{equation}

Next, we show that the region $r\in [z^{1.8}, n^{1/6}]$ is negligible. For this we want to apply Lemma~\ref{Lem:StartMove} to replace $p_{c_k+\wt{z},n}^{B,h\pm g_k}(r,\widetilde{J})$ by $p_{c_k+\wt{z},n}^{B,h\pm g_k}(w,\widetilde{J})$, paying a cost that it polynomial in $r$. However, Lemma~\ref{Lem:StartMove} is only applicable for probabilities of the form $p_{c_k;n}^{\dots}(\dots)$ such that first we need to compare $p_{c_k+\wt{z},n}^{B,g\pm g_k}(r,\widetilde{J})$ with $p_{c_k,n}^{B,h\pm g_k}(r,\widetilde{J})$. Similarly to \eqref{eq:BarrSplit1}, we get
\begin{equation}\label{eq:OneSplit}
\inf_{x \in [r,r+1]} p_{c_k+\wt{z},n}^{B,h\pm g_k}(x,\widetilde{J}) \le p_{c_k,n}^{B,h\pm g_k}(r,\widetilde{J})\cdot p_{c_k,c_k+\wt{z}}^{B,h\pm g_k}(r-1,[r-1,r])^{-1}.
\end{equation}
A direct computation, which we postpone to the end of this proof yields that 
\begin{equation}\label{eq:InfNice}
\inf_{x \in [r,r+1]} p_{c_k+\wt{z},n}^{B,h\pm g_k}(x,\widetilde{J})\le \widetilde{O}_n(1) p_{c_k+\wt{z},n}^{B,h\pm g_k}(x,\widetilde{J}).
\end{equation} 
We once more use $\wt{z} \ll w$, i.e.\@ that breaching the barrier before time $c_k+\wt{z}$ is unlikely, to get 
\begin{equation}\label{eq:SplitCons}
p_{c_k,c_k+\wt{z}}^{B,h\pm g_k}(r-1,[r-1,r])^{-1} \le \widetilde{O}_{\text{unif},k\le k_0(n)}(\wt{z}^{1/2}).
\end{equation}
Plugging \eqref{eq:OneSplit}, \eqref{eq:InfNice}, \eqref{eq:SplitCons} and Lemma~\ref{Lem:StartMove} into \eqref{eq:PertUBMark} yields
\begin{equation}\label{eq:PerturbUBHigh}
\begin{aligned}
&\int_{z^{1.8}}^{n^{1/6}}  g_{z-W_{c_k+\wt{z};c_k}^{h\pm g_k},\sigma_{c_k+\wt{z};c_k}^2}(r)p_{c_k+\wt{z},n}^{B,h\pm g_k}(r,\wt{J})\;\mathrm{d}r  \\
&\le \widetilde{O}_n(1)p_{c_k,n}^{B,h\pm g_k}(w,\wt{J})\int_{z^{1.8}}^\infty \frac{\wt{z}^{1/2}}{\sqrt{2\pi}\sigma_{c_k+\wt{z};c_k}}e^{-\frac{(r-z+W_{c_k+\wt{z};c_k}^{h\pm g_k})^2}{\sigma_{c_k+\wt{z};c_k}^2}} r^{C_2(W,h,n)} \;\mathrm{d}r+C(\lambda)e^{-\lambda \log(n)}\\
&\le \wt{O}_{\text{unif}, k\le k_0(n)}^n(\log(c_k)^{-2})p_{c_k,n}^{B,h\pm g_k}(w,\wt{J})+C(\lambda)e^{-\lambda \log(n)}.
\end{aligned}
\end{equation}

Finally, we deal with $r>n^{1/6}$ in this region we upper bound $p_{c_k+\wt{z},n}^{B,h\pm g_k}(r,\wt{J})\le 1$ and use Gaussian tail bounds to see
\begin{equation}\label{eq:UnreasonableRegion}
\begin{aligned}
\int_{n^{1/6}}^{\infty} g_{z-W_{c_k+\wt{z};c_k}^{h\pm g_k},\sigma_{c_k+\wt{z};c_k}^2}(r)p_{c_k+\wt{z},n}^{B,h\pm g_k}(r,\wt{J})\;\mathrm{d}r &\le \int_{n^{1/6}}^{\infty} g_{z-W_{c_k+\wt{z};c_k}^{h\pm g_k},\sigma_{c_k+\wt{z};c_k}^2}(r)\;\mathrm{d}r\\
&\le \wt{O}_n(e^{-cn^{1/60}}),
\end{aligned}
\end{equation}
where we also have used \eqref{Def:Ch}, \eqref{eq:Clog} and $\sigma_{c_k+\wt{z};c_k}^2\in [c_\chi z^{1.9}, C_\chi^{1.9}z^{1.9}]$.

Combining \eqref{eq:PertFKG}, \eqref{eq:PertUBMark}, \eqref{eq:PertUBLow}, \eqref{eq:PerturbUBHigh} and \eqref{eq:UnreasonableRegion} yields that for all $\lambda>0$
\begin{align*}
\frac{p_{c_k,n}^{B,h\pm g_k}(z,\wt{J})}{p_{c_k,n}^{B,h\pm g_k}(w,\wt{J})} \le 1+\wt{O}_{\text{unif}, k\le k_0(n)}^n(\log(c_k)^{-2})+p_{c_k,n}^{B,h\pm g_k}(w,\wt{J})^{-1}\wt{O}_n(e^{-\lambda\log(n)}),
\end{align*}
which is exactly the first inequality in \eqref{eq:smallStartChange}.

We finish the proof by showing \eqref{eq:InfNice}. We recall \eqref{eq:pathprobsDef} and use it to write for $x\in [r,r+1]$
\begin{align*}
p_{c_k+\wt{z},n}^{B,h\pm g_k}(x,\wt{J}) &= \int_{\wt{J}} g_{x-W_{n;c_k+\wt{z}}^{h\pm g_k}, \sigma_{n;c_k+\wt{z}}^2}(u) p_{c_k+\wt{z},n}^{B,h\pm g_k}(x|u) \;\mathrm{d}u\\
&\ge \int_{\wt{J}} g_{x-W_{n;c_k+\wt{z}}^{h\pm g_k}, \sigma_{n;c_k+\wt{z}}^2}(u) p_{c_k+\wt{z},n}^{B,h\pm g_k}(r|u) \;\mathrm{d}u,
\end{align*}
so it suffices to show
\[
g_{x-W_{n;c_k+\wt{z}}^{h\pm g_k},\sigma_{n;c_k+\wt{z}}^2}(u) g_{r-W_{n;c_k+\wt{z}}^{h\pm g_k},\sigma_{n;c_k+\wt{z}}^2}(u)^{-1} \le \widetilde{O}_n(1),
\]
which can be seen by a computation keeping in mind \eqref{Def:Ch}, \eqref{eq:Clog}, \eqref{eq:gknj}, $\sigma_{n;c_k+\wt{z}}^2 \in [c_\chi (n-c_k-\wt{z}), C_{\chi}(n-c_k-\wt{z})]$ and our conditions that $c_k\le n^{1/6}$, $z\le n^{1/11}$, $r\le n^{1/6}$. 
\end{proof}

\subsubsection{Proof of Theorem~\ref{Sa:IndArgUpperBound}}
The proof runs mostly in parallel to the proof of Theorem \ref{Sa:IndArgLowerBound}, because of this we add to the majority of the lemmata used in the proof of Theorem~\ref{Sa:IndArgUpperBound} the statement needed for the proof of Theorem~\ref{Sa:IndArgLowerBound}. 
\paragraph{Induction Basis}
We work backwards inductively in  $k$ and start with the induction basis $k = k_0$.  We apply the Markov property at time $c_{k_0}$ to get
\begin{align}
p_{0,n}^h(y,J) &= \int_0^\infty p_{0,c_{k_0}}^h(y,z)p_{c_{k_0},n}^h(z,J)\;\mathrm{d}z
\end{align}
Applying Theorem~\ref{Def:KMT} to $p_{c_{k_0},n}(z,J)$ yields that 
\begin{align}
p_{0,n}^h(y,J)\le \int_0^{\infty} p_{0,c_{k_0}}(y,z)\left(p_{c_{k_0},n}^{B,h-g_{k_0}}(y,\wt{J})+Ce^{-c\log(n)^{1.5})}\right)\;\mathrm{d}z, \label{eq:IBKMT}
\end{align}
where we used that $n/2\le n-c_{k_0}\le n$ and that $\log(n-c_{k_0})^{1.5} = g_{k_0}(j)$ for all $j\in \{c_{k_0}+1,\dots, n\}$.  By \eqref{eq:p0nlbFin} from Lemma~\ref{Lem:UnifCrudeLB}  we know that $p_{0,n}^h(y,J)^{-1}e^{-c\log(n)^{1.5}}\in \wt{O}_n(n^{-2})$. We combine this with \eqref{eq:IBKMT} to get
\begin{align*}
p_{0,n}^h(y,J) \le (1-\wt{O}_n(n^{-2}))^{-1}\int_0^\infty p_{0,c_{k_0}}(y,z)p_{c_{k_0},n}^{B,h-g_{k_0}}(y,\wt{J}),
\end{align*}
which yields the induction basis with $C_{k_0,n}(\zeta) = (1-\wt{O}_n(n^{-2}))^{-1}$.

\paragraph{Induction Step} Next, we prepare the induction step $k\to k-1$.  For this, we use the Markov property at time $c_{k-1}$ on the right-hand side of \eqref{eq:UBBrownInd} to get
\begin{equation}
\begin{aligned}
&p_{0,n}^h(y,J)\le \left(\prod_{r=k}^{k_0}C_{r,n}(\zeta)\right)\cdot \int_0^\infty \int_0^\infty p_{0,c_{k-1}}^h(y,u)p_{c_{k-1},c_k}^h(u,v)p_{c_k,n}^{B,h-g_k}(v,\wt{J})\;\mathrm{d}u\mathrm{d}v+\wt{O}_n((k_0-k)\cdot e^{-\lambda\log(n)}). \label{eq:KMTISbase}
\end{aligned}
\end{equation}
Our strategy is to apply Theorem~\ref{Def:KMT} to the term $p_{c_{k-1},c_k}^h(u,v)$. In contrast to the induction basis, the error we accrue doing this will be of order $Ce^{-c\log(c_{k-1})^{1.5}}$, and thus not negligible compared to $p_{0,n}^h(y,J)$. For $|u-v| \le C\sqrt{c_k\log(c_k)}$ the error will be negligible compared to $p_{c_{k-1},c_k}^h(u,v)$. Another difference to the induction basis is, that we can't ignore that the position at time $c_k$ can be moved by (up to) $\log(c_k)^{1.5}$ in the process of replacing $\bar{H}$ by $B^{\zeta}$ using Theorem \ref{Def:KMT}. To deal with this we will need to use that $p_{c_k,n}^{B,h-g_k}(v,\wt{J})$ is not sensitive to changes of order of size $\log(c_k)$ in $v$, which can be ensured by in addition to $|u-v|\le C\sqrt{c_k\log(c_k)}$ restricting to $v\ge C_{\log}c_k^{1/3}$, thus we need to show that the region $v \le C_{\log}c_k^{1/3}$ is negligible. The region $|u-v|\ge \sqrt{c_k}\log(n)^{1.1}$, can be handled with tail bounds for $\bar{H}_k$.

Recall the tight sequences of random variables $(\gamma_n)_{n\in\N}$, $(\widehat{\gamma}_n)_{n\in\N}$ from Lemma~\ref{Lem:UnifCrudeLB} and $C_2 := C_2(W,h,n)$ from Lemma \ref{Lem:StartMove}. We fix $2\le C_3 := C_3(n) \le \wt{O}_n(1)$ such that the following conditions hold
\begin{align}
\sup_{\mathmakebox[0.5em]{k\le k_0}} \log\big(P^{\zeta}\big[\bar{H}_{c_k;c_{k-1}}\hskip-0.1cm-W_{c_k;c_{k-1}}^h\le - C_3\sqrt{c_k\log(c_k)}  \big]\big) &\le -1-\gamma_n-\widehat{\gamma}_n,\label{eq:RWCond1C2}\\
\sup_{\mathmakebox[0.5em]{k\le k_0}} \log\big(P^{\zeta}\big[B_{c_k;c_{k-1}}^{\zeta}\hskip-0.1cm-W_{c_k;c_{k-1}}^{h+g_{k-1}}\le - C_3\sqrt{c_k\log(c_k)}  \big]\big) &\le -1-\gamma_n-\widehat{\gamma}_n,\label{eq:BMCond1C2} \\
\sum_{\mathmakebox[0.5em]{l=\lfloor C_3\sqrt{c_{k}\log(c_{k})}\rfloor}}^{\infty} P^{\zeta}\big[\bar{H}_{c_k;c_{k-1}}\hskip-0.1cm-W_{c_k;c_{k-1}}^h \ge l\big] (2(l+1))^{C_2} &\le \wt{O}_{\text{unif},k\le k_0}^n\big(c_k^{-\gamma_n-1} \big) \label{eq:Consumable2RW},\\
\sum_{\mathmakebox[0.5em]{l=\lfloor C_3\sqrt{c_{k}\log(c_k)}\rfloor}}^{\infty} P^{\zeta}\big[B_{c_k;c_{k-1}}^\zeta\hskip-0.1cm-W_{c_k;c_{k-1}}^{h+g_{k-1}} \ge l\big] (2(l+1))^{C_2} &\le \wt{O}_{\text{unif},k\le k_0}^n\big(c_k^{-\gamma_n-1} \big) \label{eq:Consumable2BM},\\
\sum_{\mathmakebox[0.5em]{l = \lfloor 2C_3\sqrt{c_k\log(c_k)}\rfloor}}^\infty P^{\zeta}\big[\bar{H}_{c_k}\hskip-0.1cm-W_{c_k}^h \ge l\big]\hskip-0.035cm (l+1)^{C_2} &\le \hskip-0.00cm\wt{O}_{\text{unif},k\le k_0}\big(c_k^{-\gamma_n-\widehat{\gamma}_n-1}\big). \label{eq:CondC2LastRW}
\end{align}
This choice is possible by using \eqref{Def:Ch}, \eqref{eq:Clog}, \cite[Theorem 15, p.\@ 52]{petrov2012sums}, Gaussian tail bounds as well as the facts that $E^{\zeta}\left[e^{\lambda \bar{\tau}_j}\right]<\infty$ for all $j\in\N$, $\lambda< |\eta|$, $\Var_{P^{\zeta}}[\tau_j] \in [c_{\chi},C_{\chi}]$ for all $j\in\N$, see\@ \cite[Eq.~(4.9), Lemma~A.1]{CernyDrewitz}. Conditions \eqref{eq:RWCond1C2}, \eqref{eq:BMCond1C2} basically ensure that the probability of $(\bar{H}_l)_{l\in\N}$ or $(B_l^\zeta)_{l\in\N}$ moving by more than $C_3\sqrt{c_k\log(c_k)}$ between $c_{k-1}$, $c_k$ is polynomially small in $c_k$. Conditions~\eqref{eq:Consumable2RW},~\eqref{eq:Consumable2BM} ensure that the probability that $\bar{H}_{c_k;c_{k-1}}$ or $B_{c_k;c_{k-1}}^\zeta$  is bigger than $l$ decays fast enough. Both of these are useful for dealing with $\bar{H}_{c_{k-1}}-W_{c_{k-1}}^h \approx y+\sqrt{c_{k-1}}$, but $\bar{H}_{c_k}-W_{c_k}^h \not \approx y+\sqrt{c_k}$ and vice versa. To shorten the next display we set $C_3\sqrt{c_k\log(c_k)} := \Delta_{c_k}$  The above outline motivates the following definition
\begin{equation}
\begin{alignedat}{2}
\wt{A}_0 &:= [0,C_{\log}c_k^{1/3} \vee (y-\sqrt{c_k\log(c_k)})],&\quad \wt{B}_0(u) &:= [0,C_{\log}c_k^{1/3}\vee(u-\Delta_{c_k})],\\
\wt{A}_1 &:= [C_{\log}c_k^{1/3}\vee (y-\sqrt{c_k\log(c_k)}) ,y+\Delta_{c_k}] ,&\quad \wt{B}_1(u) &:= [C_{\log}c_k^{1/3}\vee (u-\Delta_{c_k}),u+\Delta_{c_k}]\\
\wt{A}_2 &:= [y+\Delta_{c_k},y+\sqrt{c_k}\log(n)^{1.1}],&\quad \wt{B}_2(u) &:= [u+\Delta_{c_k},2\sqrt{c_k}\log(n)^{1.1}]\ \\
\wt{A}_3 &:= [y+\sqrt{c_k}\log(n)^{1.1},\infty),&\quad \wt{B}_3 &:= [2\sqrt{c_k}\log(n)^{1.1},\infty),\\
\end{alignedat} \label{eq:Decomposition}
\end{equation}
and 
\begin{equation} \label{eq:CutPartition}
\begin{alignedat}{2}
A_i &:= \wt{A}_i \cap [0,n^{1/11}],\quad\quad &B_i(u) &:=\wt{B}_i(u)\cap [0,n^{1/11}],\\
A_3&:=\wt{A}_3\cup (n^{1/11},\infty),\quad\quad &B_3 &:= \wt{B}_3\cup (n^{1/11},\infty),
\end{alignedat}
\end{equation}
which allows us to rewrite \eqref{eq:KMTISbase} as 
\begin{multline} \label{eq:KMTISbaseSplit}
p_{0,n}^h(y,J) \le \wt{O}_n((k_0-k)\cdot e^{-\lambda\log(n)})+\bigg(\prod_{r=k}^{k_0} C_{r,n}(\zeta)\bigg)\\
\times\smash[t]{\sum_{m=0}^3 \int_{A_m} \sum_{l=0}^3 \int_{B_l(u)}} p_{0,c_{k-1}}^h(y,u)p_{c_{k-1},c_k}^h(u,v)p_{c_k, n}^{B,h-g_k}(v,\wt{J})\;\mathrm{d}v\mathrm{d}u.
\end{multline}
Next, we state a number of lemmata handling the different summands, showing that the summands for $l\neq 1$ are negligible while for $l=1$ we can replace the random walk by a Brownian motion. For now, we only give one-sentence summaries of the proofs, the details will take up the majority of this Section. All lemmata have the same assumptions/definitions as Theorem~\ref{Sa:IndArgUpperBound} respectively~\ref{Sa:IndArgLowerBound}.We mention that the succinct summary of Lemmata \ref{Lem:IndUBk3l3}--\ref{Lem:IndUBl2k2} is that a Brownian motion respectively random walk which stays above the rugged barrier $W^h-y$ is at height roughly $\sqrt{c_k}$ above $W_{c_k}^h-y$ at time $c_k$.

\subparagraph{Comparable summand}

In this lemma, we handle $l=1$. The main ingredient of the proof is Theorem~\ref{Def:KMT}. To enable this we will need to discretize the integral over $v$. Additionally, applying Theorem~\ref{Def:KMT} to $p_{c_{k-1},c_k}^h(u,I)$, $I\subseteq B_1(u)$, will add an error of size roughly $\log(c_k)^2$ to the length of $I$, which we can remove by using Lemma~\ref{Lem:StartPert}.
\begin{Lem}[Summands with $l=1$]\label{Lem:IndUBl1}
For $m\in\{0,1,2\}$ we have
\begin{equation}
\begin{aligned}
&\int_{A_m}\int_{B_1(u)} p_{0,c_{k-1}}^h(y,u)p_{c_{k-1},c_k}^h(u,v)p_{c_k,n}^{B,h-g_k}(v,\wt{J})\;\mathrm{d}v\mathrm{d}u \\
&\le (1+\wt{O}_{\text{unif}, k\le k_0}^n(\log(c_k)^{-1.5}))\int_{A_m}\int_{B_1(u)}p_{0,c_{k-1}}^h(y,u)p_{c_{k-1},c_k}^{B,h-g_{k-1}}(u,v)p_{c_k,n}^{B,h-g_k}(v,\wt{J})\;\mathrm{d}v\mathrm{d}u\\
&\qquad+\wt{O}_{\text{unif}, k\le k_0}^n(\log(c_k)^{-1.5}) \int_{A_m}p_{0,c_{k-1}}^h(y,u)p_{c_{k-1},n}^{B,h-g_{k-1}}(u,\wt{J})\;\mathrm{d}u+\wt{O}_n(e^{-\lambda\log(n)}).
\end{aligned}
\end{equation}
\end{Lem}
\subparagraph{Negligible summands}

The next lemma handles $\max\{m,l\} = 3$, we omit its proof since it follows from Lemma 5 in \cite[pp 52--54]{petrov2012sums}. We note that by our choice of $k_0$ we have $\sqrt{c_{k_0}} \le n^{1/12} \ll n^{1/11}$, which justifies the addition of the region $(n^{1/11},\infty)$ in \eqref{eq:CutPartition}.
\begin{Lem}[Summands with $m=3$ or $l = 3$] \label{Lem:IndUBk3l3}
If $m=3$ or $l=3$ we have
\begin{align}
\int_{A_m}\int_{B_l(u)} p_{0,c_{k-1}}^h(y,u)p_{c_{k-1},c_k}^h(u,v)p_{c_k,n}^{B,h-g_k}(v,\wt{J})\;\mathrm{d}v\mathrm{d}u \le O_n(e^{-\log(n)^{1.1}}), \label{eq:k3l3UB}\\
\int_{A_m}\int_{B_l(u)} p_{0,c_{k-1}}^h(y,u)p_{c_{k-1},c_k}^{B,h+g_{k-1}}(u,v)p_{c_k,n}^{B,h+g_k}(v,\wt{K})\;\mathrm{d}v\mathrm{d}u \le O_n(e^{-\log(n)^{1.1}}). \label{eq:k3l3LB}
\end{align}
\end{Lem}
To handle the summands with $l=0$ we have two strategies. If $B_0(u) = [0,C_{\log}c_k^{1/3}]$, we can use that both $\bar{H}_{c_k;c_{k-1}}$ and $B_{c_k;c_{k-1}}^\zeta$ have variance proportional to $\sqrt{c_k}$, such that the probability that they end in an interval of length $C_{\log}c_k^{1/3}$ is small. This then allows to show that this region is negligible. If on the other hand $B_0(u) = [0,u-\Delta_k]$, the random walk respective Brownian motion has to make a big jump between $c_{k-1}$, $c_k$, which is unlikely enough to give negligibility. Motivated by this, we split $A_m = A_{m,1}\cup A_{m,2}$ with
\begin{equation}\label{eq:AmSplit}
\begin{aligned}
A_{m,1} := \{u\in A_m : B_0(u) = [0,C_{\log}c_k^{1/3}]\}\hskip8.5pt\text{and}\hskip8.5pt
A_{m,2} := \{u\in A_m : B_0(u) = [0,u-\Delta_k]\}.
\end{aligned}
\end{equation}
We note that both $A_{m,1}$, $A_{m,2}$ are intervals and that their intersection is at most a single point.

\begin{Lem}[Summands with $l = 0$, $m\in\{0,1,2\}$, $A_{m,1}$ part] \label{Lem:IndUBl0k01}
There is a $c>0$ such that
\begin{align}
\begin{aligned}
&\int_{A_{m,1}}\int_{B_0(u)} p_{0,c_{k-1}}^h(y,u)p_{c_{k-1},c_k}^h(u,v)p_{c_k,n}^{B,h-g_k}(v,\wt{J})\;\mathrm{d}u\mathrm{d}v\\
&\qquad\qquad\le \wt{O}_{\text{unif}, k\le k_0(n)}^n(c_k^{-c}) \int_{A_m} p_{0,c_{k-1}}^h(y,u)p_{c_{k-1},n}^{B,h-g_{k-1}}(u,\wt{J})\;\mathrm{d}u.
\end{aligned}\label{eq:IndUBl0k01State}\\
\begin{aligned}
&\int_{A_{m,1}}\int_{B_0(u)} p_{0,c_{k-1}}^h(y,u)p_{c_{k-1},c_k}^{B,h+g_{k-1}}(u,v)p_{c_k,n}^{B,h+g_k}(v,\wt{K})\;\mathrm{d}u\mathrm{d}v\\
&\qquad\qquad\le \wt{O}_{\text{unif}, k\le k_0(n)}^n(c_k^{-c}) \int_{A_m} p_{0,c_{k-1}}^h(y,u)p_{c_{k-1},n}^{B,h+g_{k-1}}(u,\wt{K})\;\mathrm{d}u.
\end{aligned}\label{eq:IndLBl0k01State}
\end{align}
\end{Lem}
\begin{Lem}[Summand with $l =0$, $m\in\{0,1,2\}$, $A_{m,2}$ part] \label{Lem:IndUBl0k2}
For all $\lambda>0$
\begin{align}
&\begin{aligned}
&\int_{A_{m,2}}\int_{B_0(u)} p_{0,c_{k-1}}^h(y,u)p_{c_{k-1},c_k}^h(u,v)p_{c_k,n}^{B,h-g_k}(v,\wt{J})\;\mathrm{d}v\mathrm{d}u\\
&\qquad\qquad\le \wt{O}_{\text{unif}, k\le k_0(n)}^n(c_k^{-1})\int_{A_{1}} p_{0,c_{k-1}}^h(y,u)p_{c_{k-1},n}^{B,h-g_{k-1}}(u,\wt{J})\;\mathrm{d}u+\wt{O}_n(n^{-\lambda}),
\end{aligned}\label{eq:l0k2LBTilde}\\
&\begin{aligned}
&\int_{A_{m,2}}\int_{B_0(u)} p_{0,c_{k-1}}^h(y,u)p_{c_{k-1},c_k}^{B,h+g_{k-1}}(u,v)p_{c_k,n}^{B,h+g_k}(v,\wt{K})\;\mathrm{d}v\mathrm{d}u\\
&\qquad\qquad\le \wt{O}_{\text{unif}, k\le k_0(n)}^n(c_k^{-1})\int_{A_{1}} p_{0,c_{k-1}}^h(y,u)p_{c_{k-1},n}^{B,h+g_{k-1}}(u,\wt{K})\;\mathrm{d}u+\wt{O}_n(n^{-\lambda}).
\end{aligned} \label{eq:l0k2LB}
\end{align}
\end{Lem}
In the next lemma, we consider $l=2$, $m\in\{0,1\}$. Here, we can use tail bounds on $p_{c_{k-1},c_k}^h(u,v)$ respectively $p_{c_{k-1},c_k}^{B,h+g_{k-1}}(u,v)$, which then can be integrated against the polynomial in $v$ factor we pay when we use Lemma~\ref{Lem:StartMove} to replace $p_{c_k,n}^{B,h\mp g_k}(v,\wt{J})$ by $\min_{v\in B_1(u)} p_{c_k,n}^{B,h\mp g_k}(v,\wt{J})$.

\begin{Lem}[Summands with $l=2$, $m\in\{0,1\}$] \label{Lem:IndUBl2k01} For $m \in\{0,1\}$ and $\lambda>0$ arbitrary we have
\begin{align}
&\begin{aligned}
&\int_{A_m} \int_{B_2(u)} p_{0,c_{k-1}}^h(y,u)p_{c_{k-1},c_k}^h(u,v)p_{c_k,n}^{B,h-g_k}(v,\wt{J})\;\mathrm{d}v\mathrm{d}u\\
&\qquad\qquad\le \wt{O}_{\text{unif}, k\le k_0(n)}^n(c_k^{-1}) \int_{A_m}p_{0,c_{k-1}}^h(y,u)p_{c_{k-1},n}^{B,h-g_{k-1}}(u,\wt{J})\;\mathrm{d}u+\wt{O}_n(e^{-\lambda \log(n)}),
\end{aligned}\label{eq:l2m011}\\
&\begin{aligned}
&\int_{A_m} \int_{B_2(u)} p_{0,c_{k-1}}^h(y,u)p_{c_{k-1},c_k}^{B,h+g_{k-1}}(u,v)p_{c_k,n}^{B,h+g_k}(v,\wt{K})\;\mathrm{d}v\mathrm{d}u\\
&\qquad\qquad\le \wt{O}_{\text{unif}, k\le k_0(n)}^n(c_k^{-1}) \int_{A_m}p_{0,c_{k-1}}^h(y,u)p_{c_{k-1},n}^{B,h+g_{k-1}}(u,\wt{K})\;\mathrm{d}u+\wt{O}_n(e^{-\lambda \log(n)}).
\end{aligned}\label{eq:l2m012}
\end{align} 
\end{Lem}
Finally, in the next lemma, we handle $l =m =2$. This is done using Lemma~\ref{Lem:StartMove} and the fact that being at a height in $A_2$ at time $c_{k-1}$ respectively in $B_2(u)$ at time $c_k$ has a high cost.
\begin{Lem}[Summand with $l = m = 2$] \label{Lem:IndUBl2k2}
For $\lambda>0$ arbitrary
\begin{align}\label{eq:LemIndUBl2k2C0}
&\begin{aligned}
&\int_{A_2}\int_{B_2(u)} p_{0,c_{k-1}}^h(y,u)p_{c_{k-1},c_k}^h(u,v)p_{c_k,n}^{B,h-g_k}(v,\wt{J})\;\mathrm{d}v\mathrm{d}u\\
&\qquad\qquad\le \wt{O}_{\text{unif}, k\le k_0(n)}^n(c_k^{-1}) \int_{A_1} p_{0,c_{k-1}}^h(y,u)p_{c_{k-1},n}^{B,h-g_{k-1}}(u,\wt{J})\;\mathrm{d}u +\wt{O}_n(e^{-\lambda \log(n)}),
\end{aligned}\\ \label{eq:LemIndUBl2k2C0Tilde}
&\begin{aligned}
&\int_{A_2}\int_{B_0(u)} p_{0,c_{k-1}}^h(y,u)p_{c_{k-1},c_k}^{B,h+g_{k-1}}(u,v)p_{c_k,n}^{B,h+g_k}(v,\wt{K})\;\mathrm{d}u\mathrm{d}v\\
&\qquad\qquad\le \wt{O}_{\text{unif}, k\le k_0(n)}^n(c_k^{-1}) \int_{A_1} p_{0,c_{k-1}}^h(y,u)p_{c_{k-1},n}^{B,h+g_{k-1}}(u,\wt{K})\;\mathrm{d}u.
\end{aligned}
\end{align}
\end{Lem}
Before we dive into the proofs of the lemmata, we finish the induction step assuming Lemmata~\ref{Lem:IndUBl1}--\ref{Lem:IndUBl2k2}. Plugging the lemmata into \eqref{eq:KMTISbaseSplit} and simplifying yields that there is a $c>0$ such that
\begin{align*}
p_{0,n}^h(y,J)&\le \left(\prod_{r=k}^{k_0}C_{r,n}(\zeta)\right)(1+\wt{O}_{\text{unif}, k\le k_0(n)}^n(\log(c_k)^{-1.5})\\
&\quad\times \int_{\R} p_{0,c_{k-1}}^h(y,u)p_{c_{k-1},n}^{B,h-g_{k-1}}(u,\wt{J})\;\mathrm{d}u+\wt{O}_n\left((k_0-(k-1))e^{-\lambda \log(n)}\right),
\end{align*}
where we made use of $\left(\prod_{r=k}^{k_0}C_{r,n}(\zeta)\right)\le \wt{O}_n(1)$, which will be justified in the next paragraph, to move the $e^{-\lambda \log(n)}$ and $e^{-\log(n)^{1.1}}$ terms from the lemmata into the $\wt{O}_n((k_0-(k-1))e^{-\lambda \log(n)})$. However, by setting $C_{k-1,n}(\zeta) = (1+\wt{O}_{\text{unif}, k\le k_0(n)}^n(\log(c_k)^{-1.5}))$, the last display directly implies \eqref{eq:UBBrownInd}. We observe that since $\log(c_{r+1})^{1.5} = (r+1)^{1.5}$ we have 
\begin{align}
\notag\sup_{k\le k_0} \prod_{r=k}^{k_0} C_{r,n}(\zeta) &= (1-\wt{O}_{n}^n(n^{-2}))^{-1}\sup_{k\le k_0} \prod_{r=k}^{k_0-1} (1+\wt{O}_{\text{unif},r\le k_0(n)}^n(\log(c_{r+1})^{-1.5}) \\
&\le C(1-\wt{O}_n(n^{-2}))^{-1}\wt{O}_n(1). \label{eq:TightRef}
\end{align}
The right hand side of the last display is tight as a family in $n$, which allows us to conclude Theorem~\ref{Sa:IndArgUpperBound}.\hfill\qed

\paragraph{Proof of Lemma~\ref{Lem:IndUBl1}}\label{SubSec:IndUBl1}
In this paragraph, we take $m\in\{0,1,2\}$ and deal with
\[
I_{1,m} := \int_{A_m} \int^{\min\{n^{1/11}, u+C_3\sqrt{c_k\log(c_k)}\}}_{\max(u-C_3\sqrt{c_k\log(c_k)},C_{\log}c_k^{1/3})} p_{0,c_{k-1}}^h(y,u)p_{c_{k-1},c_k}^h(u,v)p_{c_k,n}^{B,h-g_k}(v,\wt{J})\;\mathrm{d}v\mathrm{d}u.
\]

We want to apply Theorem~\ref{Def:KMT}, which only works for events and not for densities. Thus, we set $u_0 := \max(u-C_3\sqrt{c_k\log(c_k)},C_{\log}c_k^{1/3})$, $u^+ := \min\{u+C_3\sqrt{c_k\log(c_k)},n^{1/11}\}$, $r_k := \sup\{r\in\N_0 : u_0+r\log(c_k)^3\le u^+\}$, $u_r := u_0+r\log(c_k)^3$ for $r<r_k$ and $u_{r_k} := u^+$. We discretize the position at time $c_k$ into intervals $[u_r,u_{r+1}]$, maximize the barrier probability on $\{c_k,\dots, n\}$ over its starting position in this interval and apply Theorem~\ref{Def:KMT} on the barrier probability  to see that  
\begin{align}
I_{1,m} &\le \int_{A_m} \sum_{r=0}^{r_k-1} p_{0,c_{k-1}}^h(y,u)\left(p_{c_{k-1},c_k}^{B,h-g_{k-1}}(u,[u_r,u_{r+1}+2g_{k-1}(c_k)])+Ce^{-c\log(c_{k-1})^{1.5}}\right)\notag\\
&\qquad\cdot \max_{v\in [u_r,u_{r+1}]} p_{c_k,n}^{B,h-g_k}(v,\wt{J})\;\mathrm{d}u. \label{eq:PostKMTUBl1}
\end{align}
By \eqref{eq:UnifPolyLB} from Lemma~\ref{Lem:UnifCrudeLB} there are tight sequences of random variables $(\wt{C}_n)_{n\in\N}$, $(\gamma_n)_{n\in\N}$, which in particular do not depend on $k$, $u$, $r$ such that
\begin{equation} \label{eq:DemonstrateReasonInLB}
p_{c_{k-1},c_k}^{B,h-g_{k-1}(c_k)}(u,[u_r,u_{r+1}+2g_{k-1}(c_k)]) \ge \wt{C}_n^{-1}c_k^{-\gamma_n-(C_3-1)^2-6(C_3-1)(C_{\log}+C_h)^2c_{\chi}^{-1}},
\end{equation}which implies 
\begin{equation}
\frac{p_{c_{k-1},c_k}^{B,h-g_{k-1}}(u,[u_r,u_{r+1}+2g_{k-1}(c_k)])+Ce^{-c\log(c_{k-1})^{1.5}}}{p_{c_{k-1},c_k}^{B,h-g_{k-1}(c_k)}(u,[u_r,u_{r+1}+2g_{k-1}(c_k)])} \le (1+\wt{O}_{\text{unif}, k\le k_0(n)}^n(\log(c_k)^{-2})). \label{eq:MiddlePoly}
\end{equation}
Furthermore, since $u^+\le n^{1/11}$, we can apply Lemma~\ref{Lem:StartPert} to get that for all $\lambda>0$
\begin{equation}
\max_{v\in [u_r,u_{r+1}]} p_{c_k,n}^{B,h-g_k}(v,\wt{J}) \le (1+\wt{O}_{\text{unif}, k\le k_0(n)}^n(\log(c_k)^{-2}))\min_{v\in[u_r,u_{r+1}]} p_{c_k,n}^{B,h-g_k}(v,\wt{J})+\wt{O}_n(e^{-\lambda \log(n)}). \label{eq:SmallPerturb}
\end{equation}
Plugging \eqref{eq:MiddlePoly}, \eqref{eq:SmallPerturb} into \eqref{eq:PostKMTUBl1} yields that
\begin{align*}
I_{1,m} &\le \wt{O}_n(e^{-\lambda \log(n)})+\big(1+\wt{O}_{c_k}(\log(c_k)^{-2})\big)\\
&\qquad\times\int_{A_m}\sum_{r=0}^{r_k-1}p_{0,c_{k-1}}^h(y,u)p_{c_{k-1},c_k}^{B,h-g_{k-1}}(u,[u_r,u_{r+1}+2g_{k-1}(c_k)])\min_{v\in[u_r,u_{r+1}]}p_{c_k,n}^{B,h-g_k}(v,\wt{J})\;\mathrm{d}u.
\end{align*}
If we could replace the interval $[u_r,u_{r+1}+2g_{k-1}(c_k)]$ in the middle factor in the integral by $[u_r,u_{r+1}]$ we'd be done, thus we write
\begin{equation}
\begin{aligned}
I_{1,m} &\le\wt{O}_n(e^{-\lambda \log(n)})+\big(1+\wt{O}_{\text{unif}, k\le k_0(n)}^n(\log(c_k)^{-2})\big)\\
&\quad\times\Bigg(\int_{A_m}\sum_{r=0}^{r_k-1}p_{0,c_{k-1}}^h(y,u)p_{c_{k-1},c_k}^{B,h-g_{k-1}}(u,[u_r,u_{r+1}])\min_{v\in[u_r,u_{r+1}]}p_{c_k,n}^{B,h-g_k}(v,\wt{J})\;\mathrm{d}u \\
&\quad+\int_{A_m}\sum_{r=0}^{r_k-1}p_{0,c_{k-1}}^h(y,u)p_{c_{k-1},c_k}^{B,h-g_{k-1}}(u,u_{r+1}+[0,2g_{k-1}(c_k)])\min_{v\in[u_r,u_{r+1}]}p_{c_k,n}^{B,h-g_k}(v,\wt{J})\;\mathrm{d}u\Bigg)
\end{aligned}\label{eq:UBl1I1k}
\end{equation}
and we only need to show that the last summand is negligible. First, we consider $r<r_k-1$.  We recall the notation $p_{c_{k-1},c_k}^{B,h-g_{k-1}}(u|z)$ from \eqref{eq:pathprobsDef} and use that it is monotone in $u$ and $z$ to get
\begin{align*}
&p_{c_{k-1},c_k}^{B,h-g_{k-1}}(u,[u_{r+1},u_{r+1}+2g_{k-1}(c_k)]) \\
&\qquad\le \int_{u_{r+1}}^{u_{r+1}+2g_{k-1}(c_k)}g_{u-W_{c_k;c_{k-1}}^{h-g_{k-1}},\sigma_{c_k;c_{k-1}}^2}(z)\cdot p_{c_{k-1},c_k}^{B,h-g_{k-1}}(u|u_{r+1}+2g_{k-1}(c_k))\;\mathrm{d}z,\\
&p_{c_{k-1},c_k}^{B,h-g_{k-1}}(u,[u_{r+1},u_{r+2}])\\
&\qquad\ge \int_{u_{r+1}+2g_{k-1}(c_k)}^{u_{r+2}} g_{u-W_{c_k;c_{k-1}}^{h-g_{k-1}},\sigma_{c_k;c_{k-1}}^2}(z)\cdot p_{c_{k-1},c_k}^{B,h-g_{k-1}}(u|u_{r+1}+2g_{k-1}(c_k))\;\mathrm{d}z,
\end{align*}
where we used that the probability of a Brownian bridge staying above a curve is monotone in the endpoint. From here a Gaussian calculation yields that
\begin{align}
\frac{p_{c_{k-1},c_k}^{B,h-g_{k-1}}(u,[u_{r+1},u_{r+1}+2g_{k-1}(c_k)])}{p_{c_{k-1},c_k}^{B,h-g_{k-1}}(u,[u_{r+1},u_{r+2}])} &\le \wt{O}_{\text{unif}, k\le k_0(n)}^n(\log(c_k)^{-1.5}). \label{eq:UBl1GoodSummands}
\end{align}
By definition we have $u_{r+2}\le u^+\le n^{1/11}$. Thus, we can apply Lemma~\ref{Lem:StartPert} to get
\begin{align}
\min_{v\in[u_{r},u_{r+1}]}p_{c_k,n}^{B,h-g_k}(v,\wt{J}) &\le (1+\wt{O}_{\text{unif}, k\le k_0(n)}^n(\log(c_k)^{-2})\min_{v\in[u_{r+1},u_{r+2}]}p_{c_k,n}^{B,h-g_k}(v,\wt{J})\notag\\
&\qquad+\wt{O}_{n}(e^{-\lambda \log(n)}) \label{eq:UBl1MinsComp}.
\end{align}
Plugging \eqref{eq:UBl1GoodSummands}, \eqref{eq:UBl1MinsComp} into \eqref{eq:UBl1I1k} yields that
\begin{equation}\label{eq:IndUBl1I1kzwischen}
\begin{aligned}
&I_{1,m}\le \wt{O}_n(e^{-\lambda \log(n)})+\left(1+\wt{O}_{\text{unif}, k\le k_0(n)}^n(\log(c_k)^{-1.5})\right)\\
&\times \Bigg(\int_{A_m} \sum_{r=0}^{r_k-1} p_{0,c_{k-1}}^h(y,u)p_{c_{k-1},c_k}^{B,h-g_{k-1}}(u,[u_r,u_{r+1}]) \min_{v\in[u_r,u_{r+1}]}p_{c_k,n}^{B,h-g_k}(v,\wt{J})\;\mathrm{d}u\Bigg)\\
& +\int_{A_m}p_{0,c_{k-1}}^h(y,u)p_{c_{k-1},c_k}^{B,h-g_{k-1}}(u,[u_{r_k},u_{r_k}+2g_{k-1}(c_k)])\min_{v\in[u_{r_k-1},u_{r_k}]}p_{c_k,n}^{B,h-g_k}(v,\wt{J})\;\mathrm{d}u
\end{aligned}
\end{equation}
It remains to deal with $r  = r_k-1$.  Minimizing over the starting position of the barrier probability on $\{c_{k+1},\dots, n\}$ yields
\begin{equation}\label{eq:UBl1mmt2}
\begin{aligned}
\int_{A_m}\hskip-2ptp_{0,c_{k-1}}^h(y,u)p_{c_{k-1},n}^{B,h-g_{k-1}}(u,\wt{J})\;\mathrm{d}u &\ge\int_{A_m}\hskip-2ptp_{0,c_{k-1}}^h(y,u) \bigg(p_{c_{k-1},c_k}^{B,h-g_{k-1}}(u,u_{r_k}+[0,2g_{k-1}(c_k)])\\
&\hskip1.5cm\times \min_{v\in u_{r_k}+[2g_{k-1}(c_k),\log(c_k)^3]} p_{c_k,n}^{B,h-g_k}(v,\wt{J})\bigg)\;\mathrm{d}u.
\end{aligned}
\end{equation}
Now, copying the argument for $r<r_k-1$, we can infer from \eqref{eq:UBl1mmt2} that 
\begin{equation}
\begin{aligned}
&\frac{\int_{A_m}p_{0,c_{k-1}}^h(y,u)p_{c_{k-1},c_k}^{B,h-g_{k-1}}(u,[u_{r_k},u_{r_k}+2g_{k-1}(c_k)])\min_{v\in [u_{r_k-1},u_{r_k}]}p_{c_k,n}^{B,h-g_k}(v,\wt{J})\;\mathrm{d}u}{\int_{A_m}p_{0,c_{k-1}}^h(y,u)p_{c_{k-1},n}^{B,h-g_{k-1}}(u,\wt{J})\;\mathrm{d}u}\\
&\le \wt{O}_{\text{unif},k\le k_0(n)}^n(\log(c_k)^{-1.5}). 
\end{aligned}\label{eq:BadSummand}
\end{equation}
Furthermore, using that the marginals of Brownian motion have a continuous density and additivity of measures, yields that
\begin{align}
&\int_{A_m} \int_{B_1(u)} p_{0,c_{k-1}}^h(y,u)p_{c_{k-1},c_k}^{B,h-g_{k-1}}(u,v)p_{c_k,n}^{B,h-g_k}(v,\wt{J})\;\mathrm{d}u\mathrm{d}v\notag\\
&\ge\int_{A_m}\sum_{r=0}^{r_k-1}p_{0,c_{k-1}}^h(y,u)p_{c_{k-1},c_k}^{B,h-g_{k-1}}(u,[u_r,u_{r+1}])\min_{v\in[u_r,u_{r+1}]} p_{c_k,n}^{B,h-g_k}(v,\wt{J})\;\mathrm{d}u.\label{eq:UBl1TrivAdd}
\end{align}
Now, plugging \eqref{eq:BadSummand}, \eqref{eq:UBl1TrivAdd} into \eqref{eq:IndUBl1I1kzwischen}
yields that
\begin{align*}
\pushQED{\qed}
&I_{1,m} \le \wt{O}_n(e^{-\lambda \log(n)})+\wt{O}_{\text{unif}, k\le k_0(n)}^n(\log(c_k)^{-1.5})\int_{A_m} \hskip-3ptp_{0,c_{k-1}}^h(y,u)p_{c_{k-1},n}^{B,h-g_{k-1}}(u,\wt{J})\;\mathrm{d}u\\
&+(1+\wt{O}_{c_k}(\log(c_k)^{-1.5})) \hskip-2pt\int_{A_m}\int_{B_1(u)} \hskip-3ptp_{0,c_{k-1}}^h(y,u)p_{c_{k-1},c_k}^{B,h-g_{k-1}}(u,v)p_{c_k,n}^{B,h-g_k}(v,\wt{J})\;\mathrm{d}u\mathrm{d}v. \hfill \qedhere\popQED
\end{align*}
\paragraph{Proof of Lemma~\ref{Lem:IndUBl0k01}} \label{Sec:PrfUBIndl0k01}
In this paragraph, we deal with
\[
I_{2,m}:=\int_{A_{m,1}}\int_{B_0(u)} p_{0,c_{k-1}}^h(y,u)p_{c_{k-1},c_k}^h(u,v)p_{c_k,n}^{B,h-g_k}(v,\wt{J})\;\mathrm{d}v\mathrm{d}u,\quad m\in \{0,1,2\}.
\]
By definition for $u\in A_{m,1}$ we have $B_0(u) = [0,C_{\log}c_k^{1/3}] =: B_0$, which is independent of $u$.  Thus, we can swap the order of integration, maximize the barrier probability on time-frame $\{c_k+1,\dots, n\}$ over its starting position and apply Theorem~\ref{Def:KMT} to see
\begin{multline}
I_{2,m} \le (1+\wt{O}_{\text{unif}, k\le k_0(n)}^n(\log(c_k)^{-2})\int_{A_{m,1}}\bigg( p_{0,c_{k-1}}^h(y,u)\\
\times p_{c_{k-1},c_k}^{B,h-g_{k-1}}(u,[0,C_{\log}c_k^{1/3}+2g_{k-1}(c_k)])\bigg)\;\mathrm{d}u\cdot  \max_{v\in[0,C_{\log}c_k^{1/3}]} p_{c_k,n}^{B,h-g_k}(v,\wt{J}), \label{eq:Oct1707}
\end{multline}
where the KMT-error is dealt with analogously to \eqref{eq:MiddlePoly} by applying Lemma~\ref{Lem:UnifCrudeLB}. One important thing to keep in mind for this, is that $u\in A_{m,1}$ implies $u-C_3\sqrt{c_k\log(c_k)}\le C_{\log}c_k^{1/3}$, i.e.\@ $u\le C_{\log}c_k^{1/3}+C_3\sqrt{c_k\log(c_k)}$.

Set $\wt{B}_1 := [C_{\log}c_k^{1/3},C_{\log}c_k^{2/5}]$. Recalling our aim \eqref{eq:IndUBl0k01State} and using the Markov-property at time $c_{k-1}$ on the RHS of it, it is enough to compare \eqref{eq:Oct1707} to 
\begin{multline}\label{eq:UBIndl0k01LBGauss}
\int_{A_{m,1}} \int_{\wt{B}_1} p_{0,c_{k-1}}^h(y,u) p_{c_{k-1},c_k}^{B,h-g_{k-1}}(u,v) p_{c_k,n}^{B,h-g_k}(v,\wt{J})\;\mathrm{d}v\mathrm{d}u\\
\ge \int_{A_{m,1}} p_{0,c_{k-1}}^h(y,u)p_{c_{k-1},c_k}^{B,h-g_{k-1}}(u,\wt{B}_1)\min_{v\in \wt{B}_1} p_{c_k,n}^{B,h-g_k}(v,\wt{J})\;\mathrm{d}u .
\end{multline}
In the proof of \eqref{eq:k3l3LB} we arrive at precisely the same situation, only replacing $-g_{k-1}$ with $+g_{k-1}$ everywhere except for $[0,C_{\log}c_k^{1/3}+2g_{k-1}(c_k)]$ in \eqref{eq:Oct1707}, which is replaced by $[0,C_{\log}c_k^{1/3}]$. The comparison of \eqref{eq:Oct1707} with \eqref{eq:UBIndl0k01LBGauss} uses two ingredients, which we state next.\begin{Lem}\label{Lem:IndUBl0k01Barr}
Fix $n\in\N$, $k<k_0$, $m\le 2$, $u\in A_{m,1}$. We have
\begin{equation}
\frac{p_{c_{k-1},c_k}^{B,h\pm g_{k-1}}(u, [0,C_{\log}c_k^{1/3}+2g_{k-1}(c_k)])}{p_{c_{k-1},c_k}^{B,h\pm g_{k-1}}(u,\wt{B}_1)} \le \wt{O}_{\text{unif}, k\le k_0(n)}^n(c_k^{-1/15}e^{\log(c_k)^{2/3}}).\label{eq:IndUBl0k01Barr}
\end{equation}
\end{Lem}
\begin{Lem}\label{Lem:IndUBl0k01MaxMin}
Fix $n\in\N$, $k<k_0$, and $J = [a,b]$ with $0\le a<(C+1)\sqrt{n\log(n)}$. We have 
\begin{equation}
\frac{\max_{v\in B_0(0)} p_{c_k,n}^{B,h\pm g_k}(v,J)}{\min_{v\in \wt{B}_1} p_{c_k,n}^{B,h\pm g_k}(v,J) } \le \wt{O}_{\text{unif}, k\le k_0(n)}^n\left(1\right). \label{eq:IndUBl0k01MaxMin}
\end{equation}
\end{Lem}

Before we prove Lemmata~\ref{Lem:IndUBl0k01Barr},~\ref{Lem:IndUBl0k01MaxMin}, we use them to finish the proof of Lemma~\ref{Lem:IndUBl0k01}. Plugging \eqref{eq:IndUBl0k01Barr}, \eqref{eq:IndUBl0k01MaxMin} into \eqref{eq:Oct1707} yields that
\begin{equation}
\begin{aligned} \label{eq:IndUBl0k01Fin}
I_{2,m} &\le (1+\wt{O}_{\text{unif}, k\le k_0(n)}^n(\log(c_k)^{-2})) \wt{O}_{\text{unif}, k\le k_0(n)}^n(c_k^{-1/15}e^{\log(c_k)^{2/3}})\\
&\times\wt{O}_{\text{unif}, k\le k_0(n)}^n(e^{\log(c_k)^{2/3}})\int_{A_{m,1}}p_{0,c_{k-1}}^h(y,u)p_{c_{k-1},c_k}^{B,h-g_{k-1}}(u,\wt{B}_1)\min_{v\in \wt{B}_1}p_{c_k,n}^{B,h-g_k}(v,\wt{J})\;\mathrm{d}u.
\end{aligned}
\end{equation}
Comparing \eqref{eq:IndUBl0k01Fin} with \eqref{eq:UBIndl0k01LBGauss} yields the claim of Lemma~\ref{Lem:IndUBl0k01}.\hfill \qed

For both of the following proofs we recall the notation $p_{c_{k-1},c_k}^{B,h\pm g_{k-1}}(u|z)$ from \eqref{eq:pathprobsDef}.
\begin{proof}[Proof of Lemma~\ref{Lem:IndUBl0k01Barr}]
Using that $p_{c_{k-1},c_k}^{B,h\pm g_{k-1}}(u|z)$ is monotone in $z$, we get
\begin{equation}\label{eq:IndUBl0k01BarrStep1}
\begin{aligned}
&p_{c_{k-1},c_k}^{B,h\pm g_{k-1}}(u,[0,C_{\log}c_k^{1/3}+2g_{k-1}(c_k)])\\
&\qquad\le \int_0^{C_{\log}c_k^{1/3}+2g_{k-1}(c_k)} g_{u\mp g_{k-1}(c_k)-W_{c_k;c_{k-1}}^h,\sigma_{c_k;c_{k-1}}^2}(z)\;\mathrm{d}z\cdot p_{c_{k-1},c_k}^{B,h\pm g_{k-1}}(u|C_{\log}c_k^{1/3}+2g_{k-1}(c_k)),\\[0.25cm]
&p_{c_{k-1},c_k}^{B,h\pm g_{k-1}}(u, \wt{B}_1)\\
&\qquad\ge \int_{C_{\log}c_k^{1/3}+2g_{k-1}(c_k)}^{C_{\log}c_k^{2/5}} g_{u\mp g_{k-1}(c_k)-W_{c_k;c_{k-1}}^h,\sigma_{c_k;c_{k-1}}}(z)\;\mathrm{d}z\cdot p_{c_{k-1},c_k}^{B,h\pm g_{k-1}}(u|C_{\log}c_k^{1/3}+2g_{k-1}(c_k)).
\end{aligned}
\end{equation}
Thus, we only need to compare the Gaussian integrals.  A direct calculation, using \eqref{Def:Ch}, \eqref{eq:Clog} yields that
\begin{equation}
\frac{\int_0^{C_{\log}c_k^{1/3}+2g_{k-1}(c_k)} g_{ u\mp g_{k-1}(c_k)-W^h_{c_k;c_{k-1}} , \sigma_{c_k;c_{k-1}}^2}(z)\;\mathrm{d}z}{\int_{C_{\log}c_k^{1/3}+2g_{k-1}(c_k)}^{C_{\log}c_k^{2/5}} g_{u\mp g_{k-1}(c_k)-W^h_{c_k;c_{k-1}},\sigma_{c_k;c_{k-1}}^2}(z)\;\mathrm{d}z}  \le \wt{O}_{\text{unif}, k\le k_0(n)}^n\left(c_k^{-1/15}e^{\log(c_k)^{2/3}} \right) \label{eq:Oct1706}
\end{equation}
uniformly for $u\in A_{m,1}$. Combining \eqref{eq:IndUBl0k01BarrStep1} and \eqref{eq:Oct1706} implies \eqref{eq:IndUBl0k01Barr}.
\end{proof}

\begin{proof}[Proof of Lemma~\ref{Lem:IndUBl0k01MaxMin}]
We note that $B_0(0) = [0,C_{\log}c_k^{1/3}]$, but continue using the notation $B_0(0)$ for the sake of brevity. 

We write
\[
p_{c_k,n}^{B,h\pm g_k}(v,J) = \int_{J} g_{v-W_{n;c_k}^h\mp g_k(n),\sigma_{n;c_k}^2}(z)p_{c_k,n}^{B,h\pm g_k}(v|z)\;\mathrm{d}z.
\]
Using this, and the monotonicity of $p_{c_k,n}^{B,h\pm g_k}(v|z)$ in $v$, we get that
\begin{equation}\label{eq:MinMaxToDens}
\begin{aligned}
\frac{\max_{v\in B_0(0)} p_{c_k,n}^{B,h\pm g_k}(v,J)}{\min_{v\in \wt{B}_1} p_{c_k,n}^{B,h\pm g_k}(v,J) } &\le \frac{ \max_{v\in B_0(0)} \int_J g_{v-W_{n;c_k}^h\mp g_k(n),\sigma_{n;c_k}^2}(z)p_{c_k,n}^{B,h\pm g_k}( \max(B_0(0))|z)\;\mathrm{d}z}{\min_{v\in \wt{B}_1}\int_{J} g_{v-W_{n;c_k}^h\mp g_k(n),\sigma_{n;c_k}^2}(z)p_{c_k,n}^{B,h\pm g_k}(\min(\wt{B}_1)|z)\;\mathrm{d}z}\\
&\le \frac{ \max_{v\in B_0(0)} \int_J g_{v-W_{n;c_k}^h\mp g_k(n),\sigma_{n;c_k}^2}(z)\;\mathrm{d}z}{\min_{v\in \wt{B}_1}\int_{J} g_{v-W_{n;c_k}^h\mp g_k(n),\sigma_{n;c_k}^2}(z)\;\mathrm{d}z}\\
&\le \frac{\int_J \max_{v\in B_0(0)} g_{v-W_{n;c_k}^h\mp g_k(n),\sigma_{n;c_k}^2}(z)\;\mathrm{d}z}{\int_J \min_{v\in \wt{B}_1} g_{v-W_{n;c_k}^h\mp g_k(n),\sigma_{n;c_k}^2}(z)\;\mathrm{d}z},
\end{aligned}
\end{equation}
where the second step used that $\max(B_0(0)) = \min(\wt{B}_1)$. Because of this we only need to compare Gaussian densities.  This will be split into multiple cases, to shorten displays we locally introduce one additional notation
\[
\psi_v(z) := g_{v-W_{n;c_k}^h\mp g_k(n),\sigma_{n;c_k}^2}(z).
\]
We also note that we only need to consider $z\ge -W_{n;c_k}^h\mp g_k(n)$, since $J \subseteq \R^+$.
\paragraph{Case 1: $z\ge \frac{C_{\log}c_k^{2/5}+C_{\log}c_k^{1/3}}{2}-W_{n;c_k}^h\mp g_k(n)$.} In this region, we have 
\[ 
\max_{v\in B_0(0)} \psi_v(z) = \psi_v(\max(B_0(0)) = \psi_v(\min(\wt{B}_1)) = \min_{v\in B_1(u)} \psi_v(z).
\]
\paragraph{Case 2: $-W_{n;c_k}^h\mp g_k(n)\le z\le C_{\log}c_k^{1/3}-W_{n;c_k}^h\mp g_k(n)$.}
In this region, we have that 
\begin{align*}
\max_{v\in B_0(0)} \psi_v(z) &= \psi_{z+W_{n;c_k}^h\pm g_k(n)}(z)  =\frac{1}{\sqrt{2\pi \sigma_{n;c_k}^2}},\\
\min_{v\in \wt{B}_1} \psi_v(z) &= \psi_{C_{\log}c_k^{2/5}}(z) = \frac{1}{\sqrt{2\pi \sigma_{n;c_k}^2}}\exp\left(-\frac{(C_{\log}c_k^{2/5}-W_{n;c_k}^h\mp g_k(n)-z)^2}{2\sigma_{n;c_k}^2}\right).
\end{align*}
Thus, using that $C_{\log}c_k^{1/3}\ge W_{n;c_k}^h\pm g_k(n)+z\ge 0$, we have
\begin{align}
\frac{\max_{v\in B_0(0)} \psi_v(z)}{\min_{v\in \wt{B}_1} \psi_v(z)} &\le \exp\left(\frac{C_{\log}^2c_k^{4/5}-(2C_{\log}c_k^{2/5}-z-W_{n;c_k}^h\mp g_k(n))(W_{n;c_k}^h\pm g_k(n)+z)}{2\sigma_{n;c_k}^2} \right)\notag\\
&\le \exp\left(\frac{C_{\log}^2c_k^{4/5}}{2\sigma_{n;c_k}^2}\right)\le \wt{O}_{\text{unif}, k\le k_0(n)}^n(1), \label{eq:MaxMinC2}
\end{align}
where the last step uses that $\frac{c_k^{4/5}}{\sigma_{n;c_k}^2}$ is bounded, since $\sigma_{n;c_k}^2\in [c_\xi (n-c_k), C_{\xi}(n-c_k)]$ and $n-c_k\ge c_k$ by definition of $k_0$ in Definition~\ref{Def:KMTBasicDef}.

\paragraph{Case 3: $C_{\log}c_k^{1/3}-W_{n;c_k}^h\mp g_k(n)\le z \le \frac{C_{\log}c_k^{1/3}+C_{\log}c_k^{2/5}}{2}-W_{n;c_k}^h\mp g_k(n)$.}
In this case, we have that 
\begin{align*}
\max_{v\in B_0(0)}\psi_v(z) &= \psi_{C_{\log}c_k^{1/3}}(z) = \frac{1}{\sqrt{2\pi \sigma_{n;c_k}^2}}\exp\left(-\frac{\left(C_{\log}c_k^{1/3}-W_{n;c_k}^h\mp g_k(n)-z\right)^2}{2\sigma_{n;c_k}^2} \right),\\
\min_{v\in\wt{B}_1} \psi_v(z) &= \psi_{C_{\log}c_k^{2/5}}(z) = \frac{1}{\sqrt{2\pi \sigma_{n;c_k}^2}}\exp\left(-\frac{(C_{\log}c_k^{2/5}-W_{n;c_k}^h\mp g_k(n)-z)^2}{2\sigma_{n;c_k}^2}\right).
\end{align*}
As in the last case we get
\begin{align*}
\frac{\max_{v\in B_0(0)}\psi_v(z)}{\min_{v\in \wt{B}_1}\psi_v(z)} &\le \exp\left(\frac{
2C_{\log}c_k^{1/3}(W_{n;c_k}^h\pm g_k(n)+z)+C_{\log}^2c_k^{4/5}-2C_{\log}c_k^{2/5}(W_{n;c_k}^h\pm g_k(n)+z)}{2\sigma_{n;c_k}^2} \right)\\
&\le \exp\left(\frac{C_{\log}^2c_k^{4/5}}{2\sigma_{n;c_k}^2}\right)\exp\left(\frac{2(C_{\log}c_k^{1/3}-C_{\log}c_k^{2/5})(W_{n;c_k}^h\pm g_k(n)+z)}{2\sigma_{n;c_k}^2}\right)\le \wt{O}_{\text{unif}, k\le k_0(n)}^n(1),
\end{align*}
where the last step uses that $C_{\log}c_k^{1/3}-C_{\log}c_k^{2/5}\le 0$ and that we only consider $z$ for which $W_{n;c_k}^h\pm g_k(n)+z\ge C_{\log}c_k^{1/3}\ge 0$. 

Overall, we have proved that uniformly in $z\ge -W_{n;c_k}^h\mp g_k(n)$
\[
\frac{\max_{v\in B_0(0)} g_{v-W_{n;c_k}^h\mp g_k(n),\sigma_{n;c_k}^2}(z)}{\min_{v\in \wt{B}_1} g_{v-W_{n;c_k}^h\mp g_k(n),\sigma_{n;c_k}^2}(z)} \le \wt{O}_{\text{unif}, k\le k_0(n)}^n(1).
\]
Plugging this into \eqref{eq:MinMaxToDens} yields the claim of Lemma~\ref{Lem:IndUBl0k01MaxMin}.
\end{proof}
\paragraph{Proof of Lemma~\ref{Lem:IndUBl0k2}}
In this paragraph, we deal with
\[
I_{2,2,m} := \int_{A_{m,2}}\int_{B_0(u)} p_{0,c_{k-1}}^h(y,u)p_{c_{k-1},c_k}^h(u,c)p_{c_k,n}^{B,h-g_k}(v,\wt{J})\;\mathrm{d}v\mathrm{d}u,\quad m\in\{0,1,2\}.
\]
We recall that for $u\in A_{m,2}$ we have $B_0(u) = [0,u-C_3\sqrt{c_k\log(c_k)}]$, which yields
\begin{align*}
&I_{2,2,m}\le \int_{A_2}p_{0,c_{k-1}}^h(y,u)p_{c_{k-1},c_k}^h(u,B_0(u))\max_{v\in B_0(u)} p_{c_k,n}^{B,h-g_k}(v,\wt{J})\;\mathrm{d}u,\\
&\le \int_{A_2}p_{0,c_{k-1}}^h(y,u)P^{\zeta}\big[\bar{H}_{c_k;c_{k-1}}-W_{c_k;c_{k-1}}^h \le -C_3\sqrt{c_k\log(c_k)}\big]\max_{v\in B_0(u)} p_{c_k,n}^{B,h-g_k}(v,\wt{J})\;\mathrm{d}u.
\end{align*}
Now using the condition \eqref{eq:RWCond1C2} on $C_3$ yields that for $\gamma_n$, $\widehat{\gamma}_n$ as in Lemma~\ref{Lem:UnifCrudeLB}
\begin{equation}
I_{2,2,m} \le C\int_{A_{m,2}} p_{0,c_{k-1}}^h(y,u)c_k^{-1-\gamma_n-\widehat{\gamma}_n}\max_{v\in B_0(u)}p_{c_k,n}^{B,h-g_k}(v,\wt{J})\;\mathrm{d}u.
\end{equation}
For the proof of \eqref{eq:l0k2LB} we arrive at the same upper bound only exchanging $g_k$ by $-g_k$ everywhere and using condition \eqref{eq:BMCond1C2} instead of \eqref{eq:RWCond1C2}.

The continuation of the proof differs for $m\in\{0,1\}$, $m = 2$. We first deal with $m\in\{0,1\}$. Here we have $\max(B_0(u))\le \max(A_1)-C_3\sqrt{c_k\log(c_k)} \le y$. Thus, for $m\in\{0,1\}$, we get
\begin{equation}\label{eq:IndUBl0k2I22UBSmall}
I_{2,2,m}\le Cp_{0,c_{k-1}}^h(y,A_{m,1}) c_k^{-1-\gamma_n-\widehat{\gamma}_n} \hskip-3.045pt\max_{v\in [0,y]} p_{c_k,n}^{B,h-g_k}(v,\wt{J}) \le Cc_k^{-1-\gamma_n-\widehat{\gamma}_n}\hskip-3.045pt\max_{v\in [0,y]}p_{c_k,n}^{B-g_k}(v,\wt{J}).
\end{equation}

For $m=2$ we use that $\max_{v\in B_0(u)} v\le u\le n^{1/11} $ and Lemma~\ref{Lem:StartMove} to get that there is a tight (in $n$) sequence of random variables $C_2(W,h,n)$ as well as $C(\lambda)$ such that for all $\lambda>0$
\begin{equation}\label{eq:IndUBl0k2I22UB}
\begin{aligned}
I_{2,2,2}&\le C\int_{A_{2,2}}p_{0,c_{k-1}}^h(y,u) c_k^{-1-\gamma_n-\widehat{\gamma}_n} C_2(W,h,n)u^{C_2(W,h,n)}\;\mathrm{d}u\cdot \max_{v\in [0,y]} p_{c_k,n}^{B,h-g_k}(v,\wt{J})+C(\lambda)n^{-\lambda}\\
&\le \wt{O}_n(1)c_k^{-1-\gamma_n-\widehat{\gamma}_n} \max_{v\in[0,y]} p_{c_k,n}^{B,h-g_k}(v,\wt{J})+C(\lambda)n^{-\lambda},
\end{aligned}
\end{equation}
where we also used $\min(A_{2,2})-y \ge C_3\sqrt{c_k\log(c_k)}$ and tail bounds  from \cite[pp.~52--54,~Th.~15,~L.~5]{petrov2012sums}.\footnote{Formally, we'd need to discretize the integral over $u$ as in the proof of Lemma~\ref{Lem:IndUBl1}, but that adds no further difficulty.}.

Next, we give a corresponding lower bound on the right hand side of \eqref{eq:l0k2LBTilde}. For this, we  set $\wt{B}_1(y) = [y,y+\sqrt{c_k\log(c_k)}]$ and $A_1(y) := A_1\cap [y,y+\sqrt{c_k\log(c_k)}]$. We have 
\begin{equation} \label{eq:IndUBl0k2UBGauss} 
\begin{aligned}
\wt{I}_{2,2} &:= \int_{A_1} p_{0,c_{k-1}}^h(y,u)p_{c_{k-1},n}^{B,h-g_{k-1}}(u,\wt{J})\;\mathrm{d}u \ge \int_{A_1(y)}p_{0,c_{k-1}}^h(y,u)p_{c_{k-1},n}^{B,h-g_{k-1}}(u,\wt{J})\;\mathrm{d}u\\
&\ge \int_{A_1(y)} p_{0,c_{k-1}}^h(y,u)p_{c_{k-1},c_k}^{B,h-g_{k-1}}(u,\wt{B}_1(y))\min_{v\in \wt{B}_1(y)} p_{c_k,n}^{B,h-g_k}(v,\wt{J})\;\mathrm{d}u.
\end{aligned}
\end{equation}
By \eqref{eq:UnifPolyLB} in Lemma~\ref{Lem:UnifCrudeLB} there exists a tight sequence $(\gamma_n)_{n\in\N}$ such that for $u\in A_1(y)$
\begin{equation}
p_{c_{k-1},c_k}^{B,h-g_{k-1}}(u,\wt{B}_1(y))^{-1} \le \wt{O}_{\text{unif}, k\le k_0(n)}^n(c_k^{\gamma_n}). \label{eq:IndUBl0k2trivUB}
\end{equation}
Thus, \eqref{eq:IndUBl0k2UBGauss} implies
\begin{equation}\label{eq:IndUBl0k2UBGauss2}
\begin{aligned}
\wt{I}_{2,2}^{-1} &\le \wt{O}_{\text{unif}, k\le k_0(n)}^n(c_k^{\gamma_n}) p_{0,c_{k-1}}^h(y,A_1(y))^{-1}\big(\min_{v\in \wt{B}_1(y)} p_{c_k,n}^{B,h-g_k}(v,\wt{J})\;\mathrm{d}u\big)^{-1}\\
&\stackrel{\eqref{eq:PolyLBonp0n}}{\le} \wt{O}_{\text{unif}, k\le k_0(n)}^n(c_k^{\gamma_n+\widehat{\gamma}_n})\big(\min_{v\in\wt{B}_1(y)} p_{c_k,n}^{B,h-g_k}(v,\wt{J})\big)^{-1}. 
\end{aligned}
\end{equation}
We need one additional Lemma similar to Lemma~\ref{Lem:IndUBl0k01MaxMin}, whose proof we'll omit since it runs parallel to the proof of Lemma~\ref{Lem:IndUBl0k01MaxMin}
\begin{Lem}
Fix $n\in\N$, $k<k_0$, $J = [a,b]$ with $0\le a<(C+1)\sqrt{n\log(n)}$. We have
\begin{equation}\label{eq:IndUBl0k2MaxMin}
\frac{\max_{v\in[0,y]} p_{c_k,n}^{B,h\pm g_k}(v,J)}{\min_{v\in \wt{B}_1(y)} p_{c_k,n }^{B,h\pm g_k}(v,J)}\le \wt{O}_{\text{unif},k\le k_0(n)}^n(1).
\end{equation}
\end{Lem}

Combining \eqref{eq:IndUBl0k2I22UBSmall}~,\eqref{eq:IndUBl0k2I22UB},~\eqref{eq:IndUBl0k2UBGauss2} and \eqref{eq:IndUBl0k2MaxMin} yields that there is a $C(\lambda)>0$ such that for $m\in\{0,1,2\}$ and all $\lambda>0$
\begin{align*}
I_{2,2,m} \le \wt{O}_{\text{unif}, k\le k_0(n)}^n(c_k^{-1}) +C(\lambda)n^{-\lambda} \le \wt{O}_{\text{unif}, k\le k_0(n)}^n(c_k^{-1})\wt{I}_{2,2}+C(\lambda)n^{-\lambda}.\hfill \qedhere
\end{align*}

\paragraph{Proof of Lemma~\ref{Lem:IndUBl2k01}}
In this paragraph, we handle
\[
I_{3,m} := \int_{A_m}\int_{B_2(u)} p_{0,c_{k-1}}^h(y,u)p_{c_{k-1},c_k}^h(u,v)p_{c_k,n}^{B,h-g_k}(v,\wt{J})\;\mathrm{d}v\mathrm{d}u,\quad m\in \{0,1\}.
\]
Because of \eqref{Lem:StartMove} we have $\max(B_2(u))\le n^{1/11}\le n^{1/6}$. Thus, we can apply Lemma~\ref{Lem:StartMove} to see that there is a tight sequence of random variables $C_2 := C_2(W,h,n)$ and a $C(\lambda)$ such that for $\lambda>0$
\begin{equation}\label{eq:I3k}
\begin{aligned}
&I_{3,k}-C(\lambda)e^{-\lambda \log(n)} \\
&\ \ \ \le C_2\int_{A_m}\int_{B_2(u)}p_{0,c_{k-1}}^h(y,u)p_{c_{k-1},c_k}^h(u,v) v^{C_2}\;\mathrm{d}v\mathrm{d}u\cdot p_{c_k,n}^{B,h-g_k}(C_{\log}c_k^{1/3},\wt{J})\\
&\ \ \ \stackrel{w = v-u}{\le} C_2\int_{A_m} \int_{B_2(0)} p_{0,c_{k-1}}^h(y,u) p_{c_{k-1},c_k}^h(u,w+u)(2w)^{C_2}\;\mathrm{d}w\mathrm{d}u\cdot p_{c_k,n}^{B,h-g_k}(C_{\log}c_k^{1/3},\wt{J}),
\end{aligned}
\end{equation}
where in the last inequality line we used that for $w\in B_2(0)$, $u\in A_0\cup A_1$ we get $w+u\le 2w$. 

Now, we discretize the integral over $B_2(0)$. For this let $j_1<\dots<j_N$ be the integer points inside of $B_2(0)$, $j_0 := \min B_2(u)$, $j_{N+1}:=\max B_2(u)$ such that $B_2(u) =\bigcup_{l=0}^{N} [j_l,j_{l+1}]$. Continuing on from \eqref{eq:I3k} we get
\begin{equation} \label{eq:Ik3Discret}
I_{3,k}-C(\lambda)e^{-\lambda \log(n)}
\le \wt{O}_n(1)\int_{A_m} \sum_{l=0}^N p_{0,c_{k-1}}^h(y,u)p_{c_{k-1},c_k}^h(u,u+[j_l,j_{l+1}])(2j_{l+1})^{C_2}\mathrm{d}u\cdot p_{c_k,n}^{B,h-g_k}(C_{\log}c_k^{1/3},\wt{J}).
\end{equation}
We need to upper bound the sum over $l$ (uniformly in $u$). We have for all $u\in A_m$
\begin{align*}
\sum_{l=0}^N p_{c_{k-1},c_k}^h(u,u+[j_l,j_{l+1}])(2j_{l+1})^{C_2} &\le \sum_{l=0}^N P^{\zeta}[\bar{H}_{c_k;c_{k-1}}-W_{c_k;c_{k-1}} \ge j_l](2j_{l+1})^{C_2}\\
&\le  \phantom{i_{l=\lceil C_r}}\sum_{\mathmakebox[0.5em]{l=\lfloor C_3\sqrt{c_k\log(c_k)} \rfloor }}^N P^{\zeta}[\bar{H}_{c_k;c_{k-1}}-W_{c_k;c_{k-1}} \ge j_l](2j_{l+1})^{C_2}\\
&\stackrel{\eqref{eq:Consumable2RW}}{\le} \wt{O}_{\text{unif},k\le k_0(n)}^n\big(c_k^{-\gamma_n-1}\big),
\end{align*} 
and plugging this into \eqref{eq:Ik3Discret} yields
\begin{equation}\label{eq:I3klastUB}
I_{3,k}-C(\lambda)e^{-\lambda \log(n)} \le \wt{O}_{\text{unif},k\le k_0(n)}^n\big(c_k^{-\gamma_n-1}\big)\int_{A_m} p_{0,c_{k-1}}^h(y,u)\;\mathrm{d}u\cdot p_{c_k,n}^{B,h-g_k}\big(C_{\log}c_k^{1/3},\wt{J}\big).
\end{equation}

On the other hand, using \eqref{eq:UBIndl0k01LBGauss} and applying Lemmata~\ref{Lem:UnifCrudeLB},~\ref{Lem:IndUBl0k01MaxMin} yields that there is a tight sequence of random variables $(\gamma_n)_{n\in\N}$ such that 
\begin{multline} \label{eq:I3kUBGauss}
\bigg(\int_{A_m}p_{0,c_{k-1}}^h(y,u)p_{c_{k-1},n}^{B,h-g_{k-1}}(u,\wt{J})\;\mathrm{d}u\bigg)^{-1}\\
\le \wt{O}_{\text{unif}, k\le k_0(n)}^n(c_k^{\gamma_n})\bigg(\int_{A_m}p_{0,c_{k-1}}^h(y,u)\;\mathrm{d}u\cdot p_{c_k,n}^{B,h-g_k}(C_{\log}c_k^{1/3},\wt{J})\bigg)^{-1}.
\end{multline}
We combine \eqref{eq:I3klastUB} and \eqref{eq:I3kUBGauss} to see that for all $\lambda>0$ there is a $C(\lambda)>0$ such that
\[
I_{3,m} \le \wt{O}_{\text{unif}, k\le k_0(n)}^n(c_k^{-1}) \int_{A_m}p_{0,c_{k-1}}^h(y,u)p_{c_{k-1},n}^{B,h-g_{k-1}}(u,\wt{J})\;\mathrm{d}u+C(\lambda)e^{-\lambda \log(n)},
\]
which implies \eqref{eq:l2m011}. The proof of \eqref{eq:l2m012} is analogous, only having $+g_k$ instead of $-g_k$ and using \eqref{eq:Consumable2BM} instead of \eqref{eq:Consumable2RW}.\hfill \qed
\paragraph{Proof of Lemma~\ref{Lem:IndUBl2k2}} \label{SubSec:IndUBl2k2}
In this paragraph, we handle

\[
I_4 := \int_{A_2} \int_{B_2(u)} p_{0,c_{k-1}}^h(y,u)p_{c_{k-1},c_k}^h(u,v)p_{c_k,n}^{B,h-g_k}(v,\wt{J})\;\mathrm{d}v\mathrm{d}u.
\]
As in the last Section, we have $\max(B_2(u))\le n^{1/6}$ and thus can apply Lemma~\ref{Lem:StartMove} to see that there exist a tight sequence of random variables $C_2:= C_2(W,h,n)$ and a constant $C(\lambda)>0$ such that for all $\lambda>0$ we have
\begin{align*}
I_4-C(\lambda)e^{-\lambda \log(n)}&\le \int_{A_2} \hskip-1pt\int_{B_2(u)} \hskip-5ptp_{0,c_{k-1}}^h\hskip-0.825pt(y,u)p_{c_{k-1},c_k}^h\hskip-0.825pt(u,v) C_2 v^{C_2}\, \mathrm{d}v\mathrm{d}u\cdot p_{c_k,n}^{B,h-g_k}(C_{\log}c_k^{1/3}\hskip-2pt,\wt{J})\\
&\le \int_{v\ge y+2C_3\sqrt{c_k\log(c_k)}} p_{0,c_k}^h(y,v)C_2v^{C_2}\, \mathrm{d}v\cdot p_{c_k,n}^{B,h-g_k}(C_{\log}c_k^{1/3},\wt{J}).
\end{align*}
As in the proof of Lemma \ref{Lem:IndUBl2k01} we discretize the integral over $v$ and then use condition \eqref{eq:CondC2LastRW} to see that
\begin{align*}
I_4-C(\lambda)e^{-\lambda \log(n)} \le \wt{O}_{\text{unif}, k\le k_0(n)}^n\big(c_k^{-\gamma_n-\widehat{\gamma}_n-1}\big) p_{c_k,n}^{B,h-g_k}(C_{\log}c_k^{1/3},\wt{J}),
\end{align*}
comparing this to \eqref{eq:IndUBl0k2UBGauss2} while keeping \eqref{eq:IndUBl0k2MaxMin} in mind, yields that for all $\lambda>0$ there is a $C(\lambda)>0$ such that
\[
I_4\le \wt{O}_{\text{unif}, k\le k_0(n)}^n(c_k^{-1})\int_{A_1} p_{0,c_{k-1}}^h(y,u)p_{c_{k-1},n}^{B,h-g_{k-1}}(u,\wt{J})+C(\lambda)e^{-\lambda \log(n)},
\]
which implies \eqref{eq:LemIndUBl2k2C0}. The proof of \eqref{eq:LemIndUBl2k2C0Tilde} is analogous, only using $+g_k$ instead of $-g_k$, which in particular slightly changes the precise statement of \eqref{eq:IndUBl0k2UBGauss2}.\hfill \qed

\subsubsection{Proof of Theorem~\ref{Sa:IndArgLowerBound}}
As in the proof of Theorem~\ref{Sa:IndArgUpperBound} we proceed by induction. First, we do the induction basis, $k = k_0$. We have that
\begin{equation}
p_{0,n}^h(y,K)= \int_0^\infty p_{0,c_{k_0}}(y,z)p_{c_{k_0},n}(z,K)\;\mathrm{d}z\stackrel{T.\,\ref{Def:KMT}}{\ge} \int_0^\infty p_{0,c_{k_0}}(y,z)\left(p_{c_{k_0},n}^{B,h+g_{k_0}}(z,\wt{K})-Ce^{-c\log(n)^{1.5}} \right)\;\mathrm{d}z, \label{eq:Nov201}
\end{equation}
where we used that $n/2\le n-c_{k_0}\le n$ and that $\log(n-c_{k_0})^{1.5} = g_{k_0}(j)$ for all $j\in \{c_{k_0}+1,\dots, n\}$. By Lemma~\ref{Lem:UnifCrudeLB} we know that $p_{0,n}^h(y,K)^{-1} e^{-c\log(n)^{1.5}} \in \wt{O}_n(n^{-2})$ such that \eqref{eq:Nov201} implies
\[
p_{0,n}^h(y,K) \ge (1+\wt{O}_n(n^{-2}))^{-1}\int_0^\infty p_{0,c_{k_0}}(y,z) p_{c_{k_0},n}^{B,h+g_{k_0}}(z,\wt{K}),
\]
which yields the induction basis with $C_{k_0,n}(\zeta) = (1+\wt{O}_n(n^{-2}))^{-1}$.

Next, we do the induction step $k\to k-1$. The setup is the same as in the proof of Theorem~\ref{Sa:IndArgUpperBound}, i.e.\@ we use the decomposition from \eqref{eq:Decomposition} and handle all of the regions individually. We recall that after applying the Markov-property at time $c_{k-1}$ in the induction hypothesis we get that
\begin{equation}
\begin{aligned}
&p_{0,n}^h(y,K)\\
&\ge \left(\prod_{r=k}^{k_0}\wt{C}_{r,n}(\zeta)\right)\sum_{m=0}^3\int_{A_m}\sum_{l=0}^3\int_{B_l(u)} p_{0,c_{k-1}}^h(y,u)p_{c_{k-1},c_k}^h(u,v)p_{c_k,n}^{B,h+g_k}(v,\wt{K})\;\mathrm{d}v\mathrm{d}u-\wt{O}_n((k_0-k)\cdot e^{-\lambda\log(n)}).
\end{aligned}\label{eq:IndLBPostMarkov}
\end{equation} We note, that since we are after a lower bound we can initially drop all of the summands in \eqref{eq:Decomposition} for which we can't apply KMT for free and then just need to show that the expression we  we get is comparable to the right-hand side of \eqref{eq:LBBrownInd} (with $k-1$ instead of $k$). Thus we only keep the summands with $l=1$, $m<3$,  in \eqref{eq:IndLBPostMarkov} to see that
\begin{equation}\label{eq:IndLBPostDrop}
\begin{aligned}
&p_{0,n}^h(y,K)\\
&\ge \left(\prod_{r=k}^{k_0}\wt{C}_{r,n}(\zeta)\right)\sum_{m=0}^2 \int_{A_m}\int_{B_1(u)} p_{0,c_{k-1}}^h(y,u)p_{c_{k-1},c_k}^h(u,v)p_{c_k,n}^{B,h+g_k}(v,\wt{K})\;\mathrm{d}v\mathrm{d}u-\wt{O}_n((k_0-k)e^{-\lambda \log(n)}),
\end{aligned}
\end{equation}
The next lemma will handle the remaining summands in \eqref{eq:IndLBPostDrop}.
\begin{Lem}\label{Lem:IndLBl1}
For $m\in\{0,1,2\}$ we have 
\begin{equation}
\begin{aligned}
&\int_{A_m}\int_{B_1(u)} p_{0,c_{k-1}}^h(y,u)p_{c_{k-1},c_k}^h(u,v)p_{c_k,n}^{B,h+g_{k-1}}(v,\wt{K})\;\mathrm{d}v\mathrm{d}u\\
&\quad\ge (1-\wt{O}_{\text{unif}, k\le k_0(n)}^n(\log(c_k)^{-1.5}))\int_{A_m}\int_{B_1(u)} p_{0,c_{k-1}}^h(y,u)p_{c_{k-1},c_k}^{B,h+g_{k-1}}(u,v)p_{c_k,n}^{B,h+g_k}(v,\wt{K})\;\mathrm{d}v\mathrm{d}u\\
&\qquad-\wt{O}_{\text{unif}, k\le k_0(n)}^n(\log(c_k)^{-1.5})\int_{A_m}p_{0,c_{k-1}}^h(y,u)p_{c_{k-1},n}^{B,h+g_{k-1}}(u,\wt{K})\;\mathrm{d}u-\wt{O}_n(e^{-\lambda \log(n)}).
\end{aligned}
\end{equation}
\end{Lem}
The proof of Lemma~\ref{Lem:IndLBl1} is extremely similar to the proof of Lemma~\ref{Lem:IndUBl1} for which reason we omit it.\footnote{We mention that the analogue to \eqref{eq:MiddlePoly} fails to hold for finitely many $k$, for which the numerator is negative. The number of $k$ for which this happens will be tight in $n$, such that controlling them is no issue.  A similar argument has been done to deal with the $-Ce^{-c\log(c_k)^{1.5}}$ term in \eqref{eq:LowerBoundGotten} for $k$ small.} We set $J_m := \{0,2,3\}$, $m\le 2$, $J_3 = \{0,1,2,3\}$ and apply Lemma~\ref{Lem:IndLBl1} to \eqref{eq:IndLBPostDrop} to see that 
\begin{equation}\label{eq:IndLBPostKMT}
\begin{aligned}
&p_{0,n}^h(y,K)+\wt{O}_n((k_0-k)\cdot e^{-\lambda \log(n)})\\
&\ge \left(\prod_{r=k}^{k_0}\wt{C}_{r,n}(\zeta)\right)\Bigg(\left(1-\wt{O}_{\text{unif}, k\le k_0(n)}^n(\log(c_k)^{-1.5})\right) \int_{\R} p_{0,c_{k-1}}^h(y,u)p_{c_{k-1},n}^{B,h+g_{k-1}}(v,\wt{K})\;\mathrm{d}u\\
&-\left(1+\wt{O}_{\text{unif}, k\le k_0}^n(\log(c_k)^{-1.5})\right) \sum_{m=0}^3\int_{A_m}\sum_{l\in J_m} \int_{B_l(u)} p_{0,c_{k-1}}^h(y,u)p_{c_{k-1},c_k}^{B,h+g_{k-1}}(u,v)p_{c_k,n}^{B,h+g_k}(v,\wt{K})\;\mathrm{d}v\mathrm{d}u\Bigg),
\end{aligned}
\end{equation}
and it remains to prove that the third line is negligible compared to the second line, this follows from Lemmata~\ref{Lem:IndUBk3l3}--\ref{Lem:IndUBl2k2}, 
since plugging them into \eqref{eq:IndLBPostKMT} yields that
\begin{equation}
\begin{aligned}
p_{0,n}^h(y,K)\ge \left(\prod_{r=k}^{k_0} \wt{C}_{r,n}(\zeta)\right) 
\left(1-\wt{O}_{\text{unif}, k\le k_0}^n(\log(c_k)^{-1.5})\right)\int_{\R}p_{0,c_{k-1}}^h(y,u)p_{c_{k-1},n}^{B,h+g_{k-1}}(v,\wt{K})\;\mathrm{d}u\\
&\qquad-\wt{O}_n((k_0-(k-1))e^{-\lambda \log(n)}),
\end{aligned} \label{eq:IndLBfinal}
\end{equation}
where we used that for all $c>0$ we have 
\[
\wt{O}_{\text{unif}, k\le k_0(n)}^n(c_k^{-c})\left(1+\wt{O}_{\text{unif}, k\le k_0(n)}^n(\log(c_k)^{-1.5})\right)\le \wt{O}_{\text{unif}, k\le k_0(n)}^n(\log(c_k)^{-1.5})
\] as well as the fact that $\left(\prod_{r=k}^{k_0} \wt{C}_{r,n}(\zeta)\right) \le 1$, such that $\left(\prod_{r=k}^{k_0} \wt{C}_{r,n}(\zeta)\right) C(\lambda)e^{-\lambda \log(n)}\le \wt{O}_n(e^{-\lambda \log(n)})$.  From \eqref{eq:IndLBfinal} we get \eqref{eq:LBBrownInd} by setting $\wt{C}_{k-1,n}(\zeta) = (1-\wt{O}_{\text{unif}, k\le k_0(n)}^n(\log(c_k)^{-1.5}))$. Now the observation that
\begin{align*}
\left(\inf_{k\le k_0}\prod_{r=k}^{k_0}\wt{C}_{r,n}(\zeta)\right)^{-1} &\le \sup_{k\le k_0}(1+\wt{O}_n(n^{-2})) \prod_{r=k}^{k_0-1} (1-\wt{O}_{\text{unif}, r\le k_0(n)}^n(\log(c_{r+1})^{-1.5}))^{-1}\\
&\le C(1+\wt{O}_n(n^{-2}))\sup_{r\le k_0}\left(\prod_{r=k}^{k_0}(1+r^{-1.5})\right)^{\wt{O}_n(1)}\\
&\le C(1+\wt{O}_n(n^{-2}))\wt{O}_n(1),
\end{align*} 
is tight as family in $n$, finishes the proof of Theorem~\ref{Sa:IndArgLowerBound}. \hfill \qed
\subsection{Comparison to Brownian barrier probabilities -- the endpoint in a bounded interval} \label{Sec:BrownComparCompact}
By applying Theorems~\ref{Sa:IndArgUpperBound},~\ref{Sa:IndArgLowerBound} for $k = 0$ we get that there are tight sequences of random variables $\mathcal{C}_{1,n}(\zeta)$, $\mathcal{C}_{2,n}(\zeta)^{-1}$ such that for all $n\in\N$ big enough, $J = [a,b]$, $0<a<C\sqrt{n\log(n)}$, $|a-b|>\log(n)^3$, $y\in [y_0, \log(n)^2]\cap\N$
\begin{align*}
p_{0,n}^h(y,J) \le \mathcal{C}_{1,n}(\zeta) p_{0,n}^{B,h-g_0}(y,\wt{J}^+)\quad\text{and}\quad p_{0,n}^h(y,J) \ge \mathcal{C}_{2,n}(\zeta) p_{0,n}^{B,h+g_0}(y,\wt{J}^-),
\end{align*}
where $\wt{J}^+ = [a,b+2 \log(n-c_{k_0})^{1.5}]$, $\wt{J}^- = [a,b-2\log(n-c_{k_0})^{1.5}]$.

However, we need to control probabilities of this form for $J = [a,b]$ with $|b-a| \approx O(1)$, to do this we use the usual strategy of dividing $[n]$ into three intervals of length proportional to $n$ and apply the above two inequalities on the first and last interval to ``distance'' the random walk from the barrier. In the middle piece, we use a local limit theorem and ignore the barrier, which works since we are high enough that ignoring the barrier does not introduce a relevant error. We recall Definitions~\ref{Def:HelpstuffFirstMoment},~\ref{Def:HelpstuffSecMoment} and in particular the definitions of $\smileacc{m}_{n,\smileacc{h}}$, $\frownacc{m}_{n,\frownacc{h}}$ therein. We also mention that from the definition and the fact that $\chi_k^2\in [c_{\chi},C_{\chi}]$ is bounded uniformly in $k$, it is immediate that both $\smileacc{m}_{n,\smileacc{h}}$ and $\frownacc{m}_{n,\frownacc{h}}$ are nice in the sense of Definition~\ref{Def:NiceFunction}.

\begin{Kor} \label{Kor:RWtoBM}
Let $h = (h_n)_{n\in\N} \in \{ (\smileacc{m}_{n,\smileacc{h}})_{n\in\N}, (\frownacc{m}_{n,\frownacc{h}})_{n\in\N} \}$. We have 
\begin{equation}
\left(\wt{O}_n(1)\right)^{-1}\le \frac{p_{0,n}^h(y,J)}{p_{0,n}^{B,h}(y,J)} \le \wt{O}_n(1), \label{eq:RWtoBM}
\end{equation}
uniformly in $y\in[y_0, \log(n)^2]\cap \N$, $J = [a,b]$ with $a,b\in [0,\log(n)^2]\cap\N$, $b-a\le O_n(1)$.
\end{Kor}
We note that denominator does not have the $g_0$ term, which is introduced by Theorems~\ref{Sa:IndArgUpperBound},~\ref{Sa:IndArgLowerBound}. The reason for this is that the denominator is of the form analyzed in \cite{TimeInhKriech}, in particular we know that adding a $\pm g_0(k)$ in the first and a $\pm \widehat{g}_0(n-k)$ in the last third of the barrier does not change the value of the barrier probability by more than a constant, which can be absorbed into $\wt{O}_n(1)$.

We also note that throughout this section $\wt{O}_n(1)$ is to be read uniformly  in $y\in[y_0, \log(n)^2]\cap \N$, $J = [a,b]$ with $a,b\in [0,\log(n)^2]\cap\N$, $b-a\le O_n(1)$.

 In the proof we will need to control the height of the barrier $W_k$ between $n/3$, $2n/3$, for this we set $q := \lfloor n/3\rfloor$ and
\begin{equation}
C_{\text{mid},n} := C_{\text{mid},n}(W^h) := \max\left\{1,\sup_{k\in \{q+1,\dots, n-q\}} \frac{|W_{k;q}^h|}{\sqrt{n}}\right\}, \label{eq:Cmid}
\end{equation}
and remark that since $(h_n)_{n\in\N}$ is nice, i.e.\@ the $C_h$ from \eqref{Def:Ch} is finite,  the sequence $(C_{\text{mid},n}(W^h))_{n\in\N}$ is tight, although not $\mathbb{P}$-a.s.\@ bounded. The tightness follows from \cite[p.~50,~Theorem~12]{petrov2012sums}.

Before we proceed to the proof of Corollary~\ref{Kor:RWtoBM} we need one additional tool to upper- and lowerbound $P^{\zeta}[\bar{H}_n \in I]$ by the corresponding Gaussian probability $P^{\zeta}[B^{\zeta}_n\in I]$, even if $\min(I)\gg \sqrt{n}$. The result we cite is stated for a normalized version of $\bar{H}_n$, thus set $Z_n := \bar{H}_n/\sigma_n^2$ and denote its density by $\varphi_{Z_n}$. We recall that $g_{0,1}$ denotes the density of a standard normal random variable.

\begin{Sa}[{\cite[Corollary~1]{doi:10.1137/1102014}}] \label{Sa:CramerAsympt}
Let $x>1$, $x = O_n(n^{1/6})$ as $n\to\infty$. Then
\[
\frac{\varphi_{Z_n}(x)}{g_{0,1}(x)} = e^{c_nx^3}\big(1+O\big(x\cdot n^{-1/2}\big)\big),
\]
where $c_n = (6\sigma_n^3)^{-1}\sum_{j=1}^{n} \mathbb{E}[\bar{H}_j^3]$.
\end{Sa}
In \cite{doi:10.1137/1102014} there are three conditions (A)--(C) for \ref{Sa:CramerAsympt} to hold. They are
\begin{enumerate}[(A)]
\item $\left|\mathbb{E}[e^{z \bar{\tau}_k}]\right|\in [c,C]$ for $z\in \mathbb{C}$, $|z|\le A$ for some $c,C,A$ independent of $k$, 
\item  $\sigma_n^2/n\ge \delta>0$ for some $\delta>0$ and all $n\in\N$ and
\item $|\mathbb{E}[e^{(v+it)\bar{\tau}_k}]| \le L|t|^{-\beta}$ for $|t|\ge N$, $|v|< A$ for ``enough'' $k\in\N$ and some $L,N,\beta>0$ independent of $k$.
\end{enumerate}
In our case, (A) is a consequence of \cite[Lemma~A.1.]{CernyDrewitz}, (B) holds since $\sigma_n^2\in [c_{\chi} n, C_{\chi} n]$ and (C) can be shown using that $\tau_k$ is a hitting time of a birth-death-process with uniformly bounded jump rates and that by \cite[Eq.~(6.19)]{CernyDrewitz} under $P^{\zeta}$ it has exponential tails, uniformly in $k$. As a consequence of Theorem~\ref{Sa:CramerAsympt} we get good normal approximation of $P^{\zeta}[\bar{H}_n\in J]$.

\begin{Kor}\label{Kor:GaussApprox}
Let $J = [a,b]$ be an interval with $b-a\ge 1$, $b,a\le \sqrt{n}\log(n)$. We have that $\mathbb{P}$-a.s.
\begin{align*}
P^{\zeta}[\bar{H}_{n-q;q}\in J] &\in \big[(1-o_n(1))P^{\zeta}[B_{n-q;q}^\zeta\in J], (1+o_n(1))P^{\zeta}[B_{n-q;q}^\zeta \in J] \big].
\end{align*}
\end{Kor}

\begin{proof}[Proof of Corollary~\ref{Kor:RWtoBM}]

Recall the shorthand notations for barrier probabilities introduced in \eqref{eq:pathprobsDef} and recall that we denote the time-reversed versions of these by $\widehat{p}_{\dots}^{\dots}(\dots)$, we drop the superscript $(n)$ indicating with respect to which endpoint we reverse time, since it doesn't change throughout the proof. We also recall the definition \eqref{eq:gknj} of $g_0^{(q)}$, also here we drop the superscript $(q)$, but emphasize now, that it is not equal to $n$, but rather to $q = \lfloor n/3\rfloor$ throughout the proof.

 We start  by proving an upper bound.
\paragraph{Applying Theorem~\ref{Sa:IndArgUpperBound} on $[0,q]$, $[n-q,q]$ and  a local limit theorem on $[q,n-q]$}
 Using the Markov property of random walks we have that
\begin{equation}
\begin{aligned}
p_{0,n}^h(y,J) &= \int_0^\infty p_{0,q}^h(y,z)p_{q,n}^h(z,J)\;\mathrm{d}z \\
&\le \int_0^\infty \int_0^\infty p_{0,q}^h(y,z) P^{\zeta}\left[z-w+b+\bar{H}_{n-q;q}-W_{n-q;q}^h \in J\right]\widehat{p}_{0,q}^{h}(b,w)\;\mathrm{d}z\mathrm{d}w,
\end{aligned} \label{eq:ThirdingUB1}
\end{equation}
where we have also dropped the barrier on $[q,n-q]$. We set $I_1 := [0,C_4 C_{\text{mid},n}n^{1/2}]$, $\Delta_n := \sqrt{n}\log(n)^{-2}$, $I_k := [C_4 C_{\text{mid},n}n^{1/2}+(k-2)\Delta_n,C_4 C_{\text{mid},n}n^{1/2}+(k-1)\Delta_n]$ for $k\ge 2$, where $C_4$ is some big enough constant. We fix $C_5\ge 1$ such that
\begin{equation}\label{eq:C2'Def}
P^{\zeta}[\bar{H}_q-W_{q}^h \ge C_5\sqrt{n\log(n)}]\vee P^{\zeta}[\bar{H}_{n;n-q}-W_{n;n-q}^h \ge C_5\sqrt{n\log(n)}] \le o_n\big(p_{0,n}^{B,h}(y,J) \big), 
\end{equation} which is possible using \cite[Theorem~15, p.~52]{petrov2012sums} and \cite[Lemma~11.6]{TimeInhKriech}.\footnote{We emphasize that we can't apply the result as written, since in \cite{TimeInhKriech} the variance of the Brownian motion isn't time dependent, but that the necessary changes to the argument are as in Lemma~\ref{Lem:UnifCrudeLB}.} We set $T := \lfloor C_5\log(n)^{2.5} \rfloor$. Continuing on from \eqref{eq:ThirdingUB1} we get that
\begin{equation}
\begin{aligned}
p_{0,n}^h(y,J)&\le \sum_{k,l=1}^\infty p_{0,q}^{B,h-g_0}(y,I_k)\big(\hskip-1.5pt\max_{(z,w)\in I_k\times I_l} \hskip-10ptP^{\zeta}[z-w+b+\bar{H}_{n-q;q}-W_{n-q;q}^h \in J]\big) \widehat{p}_{0,q}(b,I_l)\\
&\le \wt{O}_n(1) \sum_{k,l=1}^T \bigg(p_{0,q}^{B,h-g_0}(\wt{I}_k) \cdot \big(\max_{(z,w)\in I_k\times I_l} \hskip-5pt P^{\zeta}[z-w+b+\bar{H}_{n-q;q}-W_{n-q;q}^h \in J] \big)\\
&\hskip4cm \times\widehat{p}_{0,q}^{B,h,-g_0}(b,\wt{I}_l) \bigg)+o_n\big(p_{0,n}^{B,h}(y,J)\big),
\end{aligned}
\end{equation}
where in the second step we have used Theorem~\ref{Sa:IndArgUpperBound} (and its time reversed analog) and have set $\wt{I}_k := [\min(I_k), \max(I_k)+2\log(n-c_{k_0})^{1.5}]$.
Next, we want to replace the middle probability by its Gaussian analogue. Since $|z-w+b-W_{n-q;q}^h| \le  (2C_4+1)C_5C_{\text{mid},n}\sqrt{n\log(n)}$ and $a,b\le \log(n)^2$, $1\le b-a\le O_n(1)$ we can apply Corollary~\ref{Kor:GaussApprox} for this purpose, such that the last display implies
\begin{multline} \label{eq:p0nyJDropHigh}
p_{0,n}^h(y,J) \le\wt{O}_n(1) \sum_{k,l=1}^T \Bigg(p_{0,q}^{B,h-g_0}(y,\wt{I}_k) \\
\times \max_{(z,w)\in I_k\times I_l} P^{\zeta}[z-w+b+B_{n-q;q}^{\zeta}-W_{n-q;q}^h \in J]\cdot \widehat{p}_{0,q}^{B,h,-g_0}(b,\wt{I}_l) \Bigg) +o_n\left(p_{0,n}^{B,h}(y,J)\right).
\end{multline}
Next, we need to reintroduce the barrier on $[q,n-q]$.  However, for $z,w\in I_1$ (even in the homogeneous setting) reintroducing the barrier costs a factor which is polynomial in $n$ so that first we show
that we can drop the $k,l = 1$ summands.
\paragraph{Dropping the summands with $k=1$, $l=1$ in \eqref{eq:p0nyJDropHigh}}
We introduce the notation 
\[
I_{\text{mid},n} := [C_4 C_{\text{mid},n}n^{1/2}, (3C_{\text{mid,n}}+1)n^{1/2}]
\] and show that 
\begin{equation}\label{eq:MiddleNotLow}
\begin{aligned}
&\sum_{l=1}^T p_{0,q}^{B,h-g_0}(y,\wt{I}_1)\cdot \left(\max_{(z,w)\in I_1\times I_l} P^{\zeta}\left[z-w+b+B_{n-q;q}^{\zeta}-W_{n-q;q}^h \in J\right]\right) \widehat{p}_{0,q}^{B,h,-g_0}(b,\wt{I}_l)\\
&\qquad\le  \wt{O}_n(1) \sum_{l=2}^{T}p_{0,q}^{B,h-g_0}(y, I_{\text{mid},n}) \frac{1}{\sqrt{2\pi \sigma_{n-q;q}^2}} e^{-\frac{l^2 \log(n)^{-4}}{2}\frac{n}{\sigma_{n-q;q}^2}} \widehat{p}_{0,q}^{B,h,-g_0}(b,\wt{I}_l),\\
&\sum_{k=1}^T p_{0,q}^{B,h-g_0}(y,\wt{I}_k)\cdot \left(\max_{(z,w)\in I_k\times I_1} P^{\zeta}\left[z-w+b+B_{n-q;q}^{\zeta}-W_{n-q;q}^h \in J\right] \right) \widehat{p}_{0,q}^{B,h,-g_0}(b,\wt{I}_1)\\
&\qquad\le  \wt{O}_n(1)\sum_{k=2}^T p_{0,q}^{B,h-g_0}(y,\wt{I}_k)\cdot  \frac{1}{\sqrt{2\pi \sigma_{n-q;q}^2}} e^{-\frac{k^2\log(n)^{-4}}{2}\frac{n}{\sigma_{n-q;q}^2}} \widehat{p}_{0,q}^{B,h,-g_0}(b,I_{\text{mid},n}).
\end{aligned} 
\end{equation}
We only show the first inequality in \eqref{eq:MiddleNotLow}, the proof of the second one is done analogously. 
As a first step towards \eqref{eq:MiddleNotLow}, we upper bound $\max_{(z,w)\in I_1\times I_l} (\dots)$. Since $J$ is an interval of length $O_n(1)$ and $B_{n-q;q}^\zeta \sim \mathcal{N}(0,\sigma_{n-q;q}^2)$, $|W_{n-q;q}^h|\le C_{n,\text{mid}}\sqrt{n}$ and $|z-w|\ge l\sqrt{n}\log(n)^{-2}$ for $(z,w) \in I_1\times I_l$ we have 
\begin{equation} \label{eq:midUB}
\max_{(z,w)\in I_1\times I_l} P^{\zeta}[z-w+b+B_{n-q;q}^{\zeta}-W_{n-q;q}^h\in J ] \le \wt{O}_n(1)\cdot \frac{1}{\sqrt{2\pi \sigma_{n-q;q}^2}} \cdot e^{-\frac{l^2\log(n)^{-4}}{2}\frac{n}{\sigma_{n-q;q}^2}}.
\end{equation}
Next, we bound $p_{0,q}^{B,h-g_0}(y,\wt{I}_1)$. For this we recall the notation $p_{0,q}^{B,h-g_0}(y|z)$ from \eqref{eq:pathprobsDef}. We can use that $p_{0,q}^{B,h-g_0}(y|z)$ is monotone in $z$ to get that
\begin{equation}
\begin{aligned}
p_{0,q}^{B,h-g_0}(y,\wt{I}_1) &=  \int_{\wt{I}_1} g_{y-W_q^h+g_0(q), \sigma_q^2}(r) p_{0,q}^{B,h-g_0}(y|r) \;\mathrm{d}r\\
&\le \int_{\wt{I}_1} g_{y-W_q^h+g_0(q),\sigma_q^2}(r) \;\mathrm{d}r\cdot p_{0,q}^{B,h-g_0}(y|3C_{n,\text{mid}}\sqrt{n}+2\log(n-c_{k_0})^{1.5}).
\end{aligned}\label{eq:p0qUB}
\end{equation}
On the other hand, setting $\widetilde{I}_{\text{mid},n} := C_4C_{\text{mid},n}\sqrt{n}+[2\log(n-c_{k_0})^{1.5},\sqrt{n}]$,~we~see
\begin{equation}
\begin{aligned}
p_{0,q}^{B,h-g_0}(y,I_{\text{mid},n})\ge \int_{\widetilde{I}_{\text{mid},n}} g_{y-W_q^h+g_0(q),\sigma_q^2}(r)\;\mathrm{d}r\hskip-1.3pt\cdot\hskip-1.3pt p_{0,q}^{B,h-g_0}(y|3C_{n,\text{mid}}\sqrt{n}+2\log(n-c_{k_0})^{1.5})
\end{aligned}\label{eq:p0qLB}
\end{equation}
and a Gaussian calculation, keeping in mind \eqref{eq:Cmid}, yields
\[
\frac{\int_{\wt{I}_1} g_{y-W_q^h+g_0(q),\sigma_q^2}(r) \;\mathrm{d}r}{\int_{\widetilde{I}_{\text{mid},n}} g_{y-W_q^h+g_0(q),\sigma_q^2}(r)\;\mathrm{d}r } \le \wt{O}_n(1),
\] which together with \eqref{eq:p0qUB}, \eqref{eq:p0qLB} implies that
\begin{equation} \label{eq:p0qp0qtild}
\frac{p_{0,q}^{B,h-g_0}(y,\wt{I}_1)}{p_{0,q}^{B,h-g_0}(y,I_{\text{mid},n})} \le \wt{O}_n(1).
\end{equation} We combine \eqref{eq:midUB},~\eqref{eq:p0qp0qtild} to take care of the summands $l\in\{2,\dots, T\}$ in the first equation of \eqref{eq:MiddleNotLow}.  It remains to deal with $l = 1$. Similarly to the derivation of  \eqref{eq:p0qp0qtild} we can show that
\begin{equation}
\widehat{p}_{0,q}^{B,h,-g_0}(b,\wt{I}_1) \le \wt{O}_n(1)\hskip-4.5pt\sum_{l=2}^{\lceil \log(n)^2\rceil} \hskip-6.5pt\widehat{p}_{0,q}^{B,h,-g_0}(b,\wt{I}_l) \le \wt{O}_n(1) \hskip-4.5pt\sum_{l=2}^{\lceil \log(n)^2\rceil}  \hskip-6.5pte^{-\frac{l^2\log(n)^{-4}}{2}\frac{n}{\sigma_{n-q;q}^2}}\widehat{p}_{0,q}^{B,h,-g_0}(b,\wt{I}_l), \label{eq:Dumm}
\end{equation}
where the last inequality used that for $l\le \lceil\log(n)^2\rceil$ we have that $e^{-\frac{l^2\log(n)^{-4}}{2}} \ge e^{-1}$. 
Since $T\ge \lceil \log(n)^2\rceil$ together \eqref{eq:midUB},~\eqref{eq:p0qp0qtild} and \eqref{eq:Dumm} imply the first equation of \eqref{eq:MiddleNotLow}. As mentioned before, the second inequality in \eqref{eq:MiddleNotLow} can be proved analogously.

We plug \eqref{eq:MiddleNotLow} into \eqref{eq:p0nyJDropHigh} to see that
\begin{equation}\label{eq:p0nUBHalf}
\begin{aligned}
&p_{0,n}^h(y,J)\le o_n\left(p_{0,n}^{B,h}(y,J)\right)+\wt{O}_n(1)\Bigg(\sum_{k,l=2}^T p_{0,q}^{B,h-g_0}(y,\wt{I}_k)\widehat{p}_{0,q}^{B,h,-g_0}(b,\wt{I}_l)  \\
&\qquad\times \left(\max_{(z,w)\in I_k\times I_l} P^{\zeta}[z-w+b+B_{n-q;q}^\zeta-W_{n-q;q}^h \in J]\right)\Bigg)+\wt{O}_n(1)\Bigg(\sum_{l=2}^T  \frac{e^{-\frac{l^2\log(n)^{-4}}{2}}}{\sqrt{2\pi \sigma_{n-q;q}^2}}\\
&\qquad\times \Big(p_{0,q}^{B,h-g_0}(y, I_{\text{mid},n})\widehat{p}_{0,q}^{B,h,-g_0}(b,\wt{I}_l)+p_{0,q}^{B,h-g_0}(y, \wt{I}_l)\widehat{p}_{0,q}^{B,h,-g_0}(b,I_{\text{mid},n})  \Big)\Bigg).
\end{aligned}
\end{equation}
\paragraph{Corresponding lower bound on $p_{0,n}^{B,h-g^\ast}(y,J)$}
Recall the definition of $g_0$ from \eqref{eq:gknj}. We set $g^\ast(j) := \mathbf{1}_{j\le q}g_0(j)+\mathbf{1}_{j\ge n-q}\widehat{g}_0(n-j)+\mathbf{1}_{q<j<n-q} g_0(q)$. We want to derive a lower bound on $p_{0,n}^{B,h-g^\ast}(y,J)$ that can be matched with \eqref{eq:p0nUBHalf}. We start by giving a lower bound, which can be compared to the first summand in \eqref{eq:p0nUBHalf}. Using the Markov property at times $q$, $n-q$ yields
\begin{equation}\label{eq:UBGaussLB}
\begin{aligned}
p_{0,n}^{B,h-g^\ast}(y,J) &\ge \sum_{k,l = 2}^{T+1} p_{0,q}^{B,h-g_0}(y,I_k)\widehat{p}_{0,q}^{B,h,-g_0}(b,I_l)\cdot \min_{(z,w)\in I_k\times I_l}p_{q,n-q}^{B,h}(z,J+w-b).
\end{aligned}
\end{equation}
To compare the last display to the first summand in \eqref{eq:p0nUBHalf} we need to replace $\wt{I}_k$ (respectively $\wt{I}_l$) by $I_k$ (respectively $I_l$), drop the barrier on $[q,n-q]$ and lower bound $\min_{(z,w)\in I_k\times I_l}$ by $\max_{(z,w)\in I_k\times I_l}$. 

\subparagraph{Dropping the barrier on $[q,n-q]$}
 Let $(B_{k;q}^{u,v})_{k\in \{q,\dots, n-q\}}$ be (under $P^\zeta$) a discrete Brownian bridge on $\{q,\dots, n-q\}$ from $u$ to $v$ with variance structure inherited from $(B_{k;q}^{\zeta})_{k\in \{q,\dots, n-q\}}$. We have
\begin{equation}\label{eq:MiddleRegBayes}
\begin{aligned}
&p_{q,n-q}^{B,h}(z,J+w-b)\cdot P^{\zeta}[z-w+b+B_{n-q;q}^{\zeta}-W_{n-q;q}^h \in J ]^{-1} \\
&\qquad\qquad\ge   \min_{v \in J+w-b+W_{n-q;q}^h} P^{\zeta}[\forall_{k\in \{q+1,\dots, n-q\}}\, B_{k;q}^{z,v}-W_{k;q}^h\ge 0 ]\\
&\qquad\qquad\stackrel{\eqref{eq:Cmid}}{\ge} \min_{v \in J+w-b+W_{n-q;q}^h} P^{\zeta}[\forall_{k\in \{q+1,\dots, n-q\}}\, B_{k;q}^{z,v}\ge C_{\text{mid},n}\sqrt{n} ].
\end{aligned}
\end{equation} 
The second probability is monotone in $z,v$. For $z,w \in I_k$, $k\ge 2$, we have 
\[
\min_{v\in J+w-b+W_{n-q;q}^h} v \ge 2C_{\text{mid},n}\sqrt{n}\quad\text{and}\quad z\ge 2C_{\text{mid},n}\sqrt{n},\]
where we used \eqref{eq:Cmid} as well. Thus, we get

\begin{multline}\label{eq:BBridgeMin}
\min_{v \in J+w-b+W_{n-q;q}^h} P^{\zeta}[\forall_{k\in \{q+1,\dots, n-q\}}\, B_{k;q}^{z,v}\ge C_{\text{mid},n}\sqrt{n} ] \\
\ge P^{\zeta}\big[\min_{k\in \{q+1,\dots, n-q\}} B_{k;q}^{2C_{\text{mid},n}\sqrt{n},2C_{\text{mid},n}\sqrt{n}} \ge C_{\text{mid},n}\sqrt{n}\big] \ge \left(\wt{O}_n(1)\right)^{-1},\end{multline}
where the last step uses standard estimates on the minimum of Brownian bridges. Plugging \eqref{eq:MiddleRegBayes}, \eqref{eq:BBridgeMin} into \eqref{eq:UBGaussLB} yields that
\begin{multline} \label{eq:p0nJGaussLBAfterBarrMidRem}
p_{0,n}^{B,h-g^\ast}(y,J)\ge \big(\wt{O}_n(1)\big)^{-1}\bigg(\sum_{k,l=2}^{T+1} p_{0,q}^{B,h-g_0}(y,I_k)\widehat{p}_{0,q}^{B,h,-g_0}(b,I_l) \\
\times \min_{(z,w)\in I_k\times I_l} P^{\zeta}[ z-w+b+B_{n-q;q}^\zeta-W_{n-q;q}^h\in J]\bigg).
\end{multline}

\subparagraph{Replace $\min_{(z,w)\in I_k\times I_l}$ by $\max_{(z,w) \in I_k\times I_l}$}
Since for $k,l\in \{2,\dots, T\}$ we have 
\begin{align*}
|z-w+b-W_{n-q;q}^h| &\le (2C_4+1) C_5 C_{\text{mid},n}\sqrt{n\log(n)},\\
\max_{(z,w)\in I_k\times I_l} |z-w| - \min_{(z,w)\in I_k\times I_l} |z-w| &\le 2\sqrt{n}\log(n)^{-2}
\end{align*} we get that
\begin{equation} \label{eq:MiddleRegionMinMax}
\frac{\max_{(z,w)\in I_k\times I_l} P^{\zeta}[ z-w+b+B_{n-q;q}^\zeta-W_{n-q;q}^h\in J] }{\min_{(z,w)\in I_k\times I_l}P^{\zeta}[ z-w+b+B_{n-q;q}^\zeta-W_{n-q;q}^h\in J]} \le \wt{O}_n(1),
\end{equation}
which together with \eqref{eq:p0nJGaussLBAfterBarrMidRem} implies
\begin{multline} \label{eq:GaussFU}
p_{0,n}^{B,h-g^\ast}(y,J)\ge \big(\wt{O}_n(1)\big)^{-1}\bigg(\sum_{k,l=2}^{T+1} p_{0,q}^{B,h-g_0}(y,I_k)\widehat{p}_{0,q}^{B,h,-g_0}(b,I_l)\\
\times\max_{(z,w)\in I_k\times I_l} P^{\zeta}[z-w+b+B_{n-q;q}^{\zeta}-W_{n-q;q}^h\in J]\bigg).
\end{multline}
\begin{subparagraph}{Replacing $\wt{I}_k$ by $I_k$}
Finally, we need to exchange $\wt{I}_k$ by $I_k$ in \eqref{eq:p0nUBHalf}. As below \eqref{eq:UBl1I1k} we can see that 
\begin{equation} \label{eq:BlowUpAs+1}
\frac{\widehat{p}_{0,q}^{B,h,-g_0}(b,\wt{I}_l)}{\widehat{p}_{0,q}^{B,h,-g_0}(b,I_{l+1})}\frac{p_{0,q}^{B,h-g_0}(y,\wt{I}_k)}{p_{0,q}^{B,h-g_0}(y,I_{k+1})} \le \wt{O}_n(1).
\end{equation}
Furthermore, since $\max_{(z,w)\in I_k\times I_l}(\dots)$ only depends on $\{x \in \R : x = z-w, z\in I_k, w\in I_l\}$, which for $k,l\ge 2$ is invariant under increasing both $k$ and $l$ by 1 we have that
\begin{multline} \label{eq:MaxIncrBothEq}
\max_{(z,w)\in I_k\times I_l} P^{\zeta}[z-w+b+B_{n-q;q}^{\zeta}-W_{n-q;q}^h\in J] \\
= \max_{(z,w)\in I_{k+1}\times I_{l+1}} P^{\zeta}[z-w+b+B_{n-q;q}^{\zeta}-W_{n-q;q}^h\in J].
\end{multline}
Combining \eqref{eq:BlowUpAs+1}, \eqref{eq:MaxIncrBothEq} yields that
\begin{equation} \label{eq:Untildify}
\begin{aligned}
&\sum_{k,l=2}^T p_{0,q}^{B,h-g_0}(y,\wt{I}_k) \cdot \big(\max_{(z,w)\in I_k\times I_l} P^{\zeta}[z-w+b+\bar{H}_{n-q;q}-W_{n-q;q}^h \in J]\big)\cdot \widehat{p}_{0,q}^{B,h,-g_0}(b,\wt{I}_l) \\
&\le \sum_{k,l = 2}^{T+1} p_{0,q}^{B,h-g_0}(y,I_k)\big(\max_{(z,w)\in I_k\times I_l} P^{\zeta}[z-w+b+\bar{H}_{n-q;q}-W_{n-q;q}^h \in J]\big)\cdot \widehat{p}_{0,q}^{B,h,-g_0}(b,I_l) .
\end{aligned}
\end{equation}
\end{subparagraph}
Now, we combine \eqref{eq:p0nUBHalf}, \eqref{eq:GaussFU} and \eqref{eq:Untildify} to get
\begin{equation}
\begin{aligned}
p_{0,n}^h(y,J)&\le \wt{O}_n(1) p_{0,n}^{B,h-g^\ast}(y,J) +o_n\left(p_{0,n}^{B,h}(y,J)\right)+\wt{O}_n(1)\Bigg(\sum_{l=2}^T  \frac{e^{-\frac{l^2\log(n)^{-4}}{2}}}{\sqrt{2\pi \sigma_{n-q;q}^2}\frac{n}{\sigma_{n-q;q}^2}}\\
&\quad\times \big(p_{0,q}^{B,h-g_0}(y, I_{\text{mid},n})\widehat{p}_{0,q}^{B,h,-g_0}(b,\wt{I}_l)+p_{0,q}^{B,h-g_0}(y, \wt{I}_l)\widehat{p}_{0,q}^{B,h,-g_0}(b,I_{\text{mid},n})  \big) \Bigg).
\end{aligned}\label{eq:p0nUBThreeQuart}
\end{equation}
Next, we deal with the final summand in \eqref{eq:p0nUBThreeQuart}.We mirror the computation leading to \eqref{eq:p0nJGaussLBAfterBarrMidRem} to see that
\begin{align*}
&p_{0,n}^{B,h-g^\ast}(y,J)\ge \big(\wt{O}_n(1)\big)^{-1}\Bigg(\sum_{l=2}^{T+1} p_{0,q}^{B,h-g_0}(y,I_{\text{mid}})\widehat{p}_{0,q}^{B,h,-g_0}(b,I_l)\\
&\hskip5cm\times \min_{(z,w)\in I_{\text{mid}}\times I_l} P^{\zeta}[z-w+b+B_{n-q;q}^{\zeta}-W_{n-q;q}\in J ]\Bigg)\\
&\ge \left(\wt{O}_n(1)\right)^{-1}\sum_{l=2}^{T+1} p_{0,q}^{B,h-g_0}(y,I_{\text{mid}}) \widehat{p}_{0,q}^{B,h,-g_0}(y,I_l) \cdot \frac{1}{\sqrt{2\pi \sigma_{n-q;q}^2}} e^{-\frac{l^2 \log(n)^{-4}}{2}\frac{n}{\sigma_{n-q;q}^2}}.
\end{align*}
Analogously, we get that
\[
p_{0,n}^{B,h-g^\ast}(y,J) \ge\left(\wt{O}_n(1)\right)^{-1}\sum_{l=2}^{T+1} p_{0,q}^{B,h-g_0}(y,I_l)\widehat{p}_{0,q}^{B,h,-g_0}(y,I_{\text{mid}})\frac{e^{-\frac{l^2\log(n)^{-4}}{2}\frac{n}{\sigma_{n-q;q}^2}}}{\sqrt{2\pi \sigma_{n-q;q}^2}}.
\]
Combining the last two displays with \eqref{eq:BlowUpAs+1},~\eqref{eq:p0nUBThreeQuart}, while keeping in mind $e^{\frac{n}{\sigma_{n-q;q}^2}\frac{(-l^2+(l+1)^2)\log(n)^{-4}}{2}} \le \wt{O}_n(1)$ for $l\le T$, yields that 
\begin{equation}
p_{0,n}^h(y,J)\le \wt{O}_n(1) p_{0,n}^{B,h-g^\ast}(y,J)+o_n\big(p_{0,n}^{B,h}(y,J)\big). \label{eq:p0nyJUBLast}
\end{equation}
\paragraph{Removing the $g^\ast$}
To finish the upper bound we use the  results from \cite{TimeInhKriech} to remove $g^\ast$, which in particular allows us to absorb the second summand in \eqref{eq:p0nyJUBLast} into the first. We recall that by definition of $g_0$, $\widehat{g}_0$ in \eqref{eq:gknj},~\eqref{eq:gknjReversed} the function $g^\ast$ is weighted by the local variance of $B^\zeta$. Thus, a similar calculation to \cite[Section 9]{TimeInhKriech} yields
\[
p_{0,n}^{B,h-g^\ast}(y,J)\le \wt{O}_n(1) p_{0,n}^{B,h}(y,J),
\]
plugging this into \eqref{eq:p0nyJUBLast} yields 
\[
p_{0,n}^h(y,J) \le \wt{O}_n(1)p_{0,n}^{B,h}(y,J),
\]
which is the upper bound in \eqref{eq:RWtoBM}.

The lower bound can be proved analogously, switching the role of $p_{0,n}^h(y,J)$ and $p_{0,n}^{B,h-g^\ast}(y,J)$ in the above and using Theorem~\ref{Sa:IndArgLowerBound} instead of Theorem~\ref{Sa:IndArgUpperBound}. There are a few minor changes, like needing to use $\widehat{I}_k := [\min(I_k),\max(I_k)-2\log(n-c_{k_0})^{1.5}]$ instead of $\wt{I}_k$ and replacing \eqref{eq:BlowUpAs+1} by 
\begin{equation} \label{eq:ShrinkDownAs+1}
\frac{\widehat{p}_{0,q}^{B,h,g_0}(b,I_l)}{\widehat{p}_{0,q}^{B,h,g_0}(b,\widehat{I}_{l+1})}\frac{p_{0,q}^{B,h+g_0}(y,I_k)}{p_{0,q}^{B,h+g_0}(y,\widehat{I}_{k+1})} \le \wt{O}_n(1)
\end{equation}
for $k,l \ge 2$. The only major change is that in the step of dropping the barrier, i.e.\@ the analogue of \eqref{eq:MiddleRegBayes} the probability that a random walk, conditioned to end at a certain point, stays above a barrier is not monotone in the start- and endpoint. Because of this we need to establish 
\begin{equation}
\begin{aligned} 
&\widetilde{O}_n(1)\sum_{k,l = 2}^{T+1} p_{0,q}^{B,h+g_0}(y,\widehat{I}_k)\widehat{p}_{0,q}^{B,h,g_0}(b,\widehat{I}_k)\cdot \min_{(z,w)\in I_k\times I_l} p_{q,n-q}^h(z, J-b+w)\\
&\ge\sum_{k,l = 2}^{T+1} p_{0,q}^{B,h+g_0}(y,\widehat{I}_k)\widehat{p}_{0,q}^{B,h,g_0}(b,\widehat{I}_k)\min_{(z,w)\in I_k\times I_l} P^{\zeta}[z-w+b+B_{n-q;q}^{\zeta}-W_{n-q;q}^h\in J ], \label{eq:p0nyJLB}
\end{aligned}
\end{equation}
which is done in the next lemma.
\end{proof}

We use the following Lemma in the above.
\begin{Lem}\label{Lem:WallRemRW}
Keep the notations and assumptions from (the proof of) Corollary~\ref{Kor:RWtoBM}. We have for all $C_4 C_{\text{mid},n}n^{1/2}\le z,w \le C_4 C_{\text{mid},n}C_5\sqrt{n\log(n)}$
\begin{equation} \label{eq:pqnqBridgeRem}
p_{q,n-q}^h(z,J-b+w) \ge \left(\wt{O}_n(1)\right)^{-1} P^{\zeta}\left[z-w+b+B_{n-q;q}^{\zeta}-W_{n-q;q}^h \in J \right].
\end{equation}
\end{Lem}
\begin{proof}
By definition of $C_{\text{mid},n}$ we have $|W_{n-q;q}^h| \le C_{\text{mid},n}\sqrt{n}$. Thus, we have
\begin{align*}
p_{0,n}^h(z,J-b+w) &\ge P^{\zeta}\left[ \forall_{k\in\{q+1,\dots n-q\}}\, z+\sum_{j=q+1}^k \bar{\tau}_j\ge C_{\text{mid},n}\sqrt{n}, z-w+b+\bar{H}_{n-q;q}-W_{n-q;q}^h \in J\right].
\end{align*}
Next, we want to use the reflection principle. However, the increments of our random walk are not symmetric and because of this we set
\[
\bar{H}_{n-q;q}^{(l)}  := \sum_{k=q+1}^{l} \bar{\tau}_j-\sum_{k = l+1}^{n-q} \bar{\tau}_j,
\]
to be the at time $l$ reflected version of $\bar{H}_{n-q;q}^{(l)}$. We also recall that there is a constant $C_{\text{breach}}>0$ such that $\min_j \bar{\tau}_j\ge -C_{\text{breach}}$, since $\bar{\tau}_j = \tau_j-E^{\zeta}[\tau_j]$ and $\tau_j\ge 0$ as a hitting time, while $\min_{j\in\N} -E^{\zeta}[\tau_j]$ is $\mathbb{P}$-a.s.\@ bounded by \cite[Lemma~A.1.]{CernyDrewitz}. Thus, we get
\begin{align*}
p_{q,n-q}^h(z,J-b+w)&\ge P^{\zeta}\left[ z-w+b+\sum_{j=q+1}^{n-q}\bar{\tau}_j-W_{n-q;q}^h \in J\right]-\\
&\qquad-\max_{\genfrac{}{}{0pt}{}{l\in \{q+1,\dots, n-q\}}{d\in [C_{\text{mid},n}\sqrt{n}-C_{\text{breach}},C_{\text{mid},n}\sqrt{n}]}} P^{\zeta}\left[ z+w-b+\bar{H}_{n-q;q}^{(l)}+W_{n-q;q}^h-2d \in -J\right]\\
&= \min_{\genfrac{}{}{0pt}{}{l\in \{q+1,\dots, n-q\}}{d\in[C_{\text{mid},n}\sqrt{n}-C_{\text{breach}},C_{\text{mid,n}}\sqrt{n}]}}\Bigg( P^{\zeta}\left[z-w+b+B_{n-q;q}^\zeta-W_{n-q;q}^h \in J\right]-\\
&\hskip6cm-P^{\zeta}\left[z+w-b-2d+B_{n-q;q}^\zeta+W_{n-q;q}^h \in -J\right]\Bigg).
\end{align*}
Now we set $I_{\text{breach}} := [C_{\text{mid},n}\sqrt{n}-C_{\text{breach}},C_{\text{mid},n}\sqrt{n}]$ and use our conditions on $z,w$ and $|W_{n-q;q}^h| \le C_{\text{mid},n}\sqrt{n}$ and Corollary~\ref{Kor:GaussApprox} to get
\begin{align*}
p_{q,n-q}^h(z,J-b+w)&\ge \min_{\genfrac{}{}{0pt}{}{l\in \{q,\dots, n-q\}}{d\in I_{\text{breach}}}}\Bigg(  (1-o_n(1))P^{\zeta}\left[z-w+b+B_{n-q;q}^{\zeta}-W_{n-q;q}^h \in J\right]\\
&\qquad-(1+o_n(1))P^{\zeta}\left[z+w-b-2d+B_{n-q;q}^\zeta+W_{n-q;q}^h \in -J\right] \Bigg).\\
\end{align*}
From here a Gaussian calculation (after potentially increasing $C_5$) and our conditions on $z,w$ together with \eqref{eq:Cmid} yield that
\[
p_{q,n-q}^h(z,J-b+w)\ge (\wt{O}_n(1))^{-1}P^{\zeta}\left[z-w+b+B_{n-q;q}^\zeta-W_{n-q;q}^h \in J\right].\qedhere
\]

\end{proof}

\subsubsection{Concluding Lemmata~\ref{Lem:FMomSimplBarrComp},~\ref{Lem:FMomHardBarrComp},~\ref{Lem:SecMomBarrCompLB},~\ref{Lem:SecMomBarrCompUB},~\ref{Lem:SecMomBarrCompComplicated}}
In this section, we prove Lemmata~\ref{Lem:FMomSimplBarrComp},~\ref{Lem:SecMomBarrCompLB},~\ref{Lem:SecMomBarrCompUB} by first applying Corollary~\ref{Kor:RWtoBM} to move to Brownian motion barrier events and then applying the results from \cite{TimeInhKriech}. 
\begin{proof}[Proof of Lemmata~\ref{Lem:FMomSimplBarrComp},~\ref{Lem:SecMomBarrCompLB},~\ref{Lem:SecMomBarrCompUB}]
Recall the definitions \eqref{Def:Knzeta},~\eqref{Def:Wnzeta}, \eqref{Def:bprobUB} and \eqref{Def:bprobLB}. It is straightforward to check that both $\frownacc{m}_{n,\frownacc{h}}(\cdot)$, $\smileacc{m}_{n,\smileacc{h}}$ fulfill \eqref{Def:Ch} and $(W_k)_{k\in\N}$ fulfills \eqref{eq:Clog}. By applying Corollary~\ref{Kor:RWtoBM} we get that 
\begin{align*}
\barrierprobUB &\le \wt{O}_n(1){P^{\zeta}[\mathcal{B}_{[n],W-\smileacc{m}_{n,\smileacc{h}}}^{y,J_x}(B_{\cdot}^{\zeta}) ]}, \\
\barrierprobLB &\in \left[\big(\wt{O}_n(1)\big)^{-1}{P^{\zeta}[\mathcal{B}_{[n],W-{\frownacc{m}_{n,\frownacc{h}}}}^{y,J_{y_0}}(B_{\cdot}^{\zeta}) ]}, \wt{O}_n(1){P^{\zeta}[\mathcal{B}_{[n],W-{\frownacc{m}_{n,\frownacc{h}}}}^{y,J_{y_0}}(B_{\cdot}^{\zeta}) ]}\right].
\end{align*}
Now we can conclude Lemmata~\ref{Lem:FMomSimplBarrComp},~\ref{Lem:SecMomBarrCompLB},~\ref{Lem:SecMomBarrCompUB} by applying the results of Sections 8 and 9 in  \cite{TimeInhKriech},  where we also used that by \eqref{Def:pn} the object we want to compare to is already defined in terms of a Brownian motion. We also point out that both $\smileacc{h}$, $\frownacc{h}$ and the linear drift towards $\log(p_n)$ have been defined depending on $\chi_k^2 = \Var_{P^{\zeta}}[\bar{\tau}_k]$ in a way accounting for the local variance of the Brownian motion such that the arguments in \cite{TimeInhKriech} go through with only minor changes, we refer to the proof of Lemma~\ref{Lem:UnifCrudeLB} to see an example of how defining the functions in a way depending on the local variance enables one to adapt arguments from \cite{TimeInhKriech}. 
\end{proof}

The proof of both Lemmata~\ref{Lem:FMomHardBarrComp} and~\ref{Lem:SecMomBarrCompComplicated} runs along similar lines as the proof of Corollary~\ref{Kor:RWtoBM}, but since we are happy with polynomial (in $k$) lower bounds on $q_{k,\text{end}}(x)$ respectively $q_{k,\text{start}}$ we can allow bigger errors. Because of this, we only sketch the proofs.

\begin{proof}[Proof sketch for Lemmata~\ref{Lem:FMomHardBarrComp},~\ref{Lem:SecMomBarrCompComplicated}]
By \eqref{eq:qkstart} we have
\[
q_{k,\text{start}} = \inf_{x\in [1,2]} p_{0,k}^{\frownacc{m}_{n,\frownacc{h}}}(y_0, J_x).
\]
We set $q_k:= \lfloor k/3\rfloor$ and $h = (h_n)_{n\in\N}$ with $h_n = \frownacc{m}_{n,\frownacc{h}}$. As in the proof of Corollary~\ref{Kor:RWtoBM} we can use the Markov property at times $q$, $k-q$ to get that for $C_6\ge 1$, $I := [C_6\sqrt{k\log(k)}, 2C_6\sqrt{k\log(k)}]$
\[
q_{k,\text{start}}\ge   \inf_{x\in [1,2]} \left(p_{0,q}^{h}(y_0,I)\widehat{p}_{0,q}^{h,(k)}(x,I)\cdot \min_{z,w\in I} p_{q,n-q}^{h}(z,J-b+w)\right).
\]
Analogous to Lemma~\ref{Lem:UnifCrudeLB} we can show that there is a tight in $n$ sequence $\bar{\gamma}_n(C_6)\le \wt{O}_n(1)$ such that 
\[
\max\left\{p_{0,q}^{h}(y_0,I)^{-1}, \left(\inf_{x\in[1,2]} \widehat{p}_{0,q}^{h,(k)}(x,I) \right)^{-1}\right\} \le \wt{O}_{\text{unif}, k\le \log(n)^7}^n(k^{\bar{\gamma}_n}).
\]
Furthermore, similarly to Lemma~\ref{Lem:WallRemRW} we can show that (potentially after increasing $C_6$), there is a $\bar{\gamma}_n'\le \wt{O}_n(1)$ such that
\[
\left(\min_{z,w\in I} p_{q,n-q}^{h}(z,J-b+w)\right)\le \wt{O}_{\text{unif},k\le \log(n)^7}^n(k^{\bar{\gamma}_n'}).
\]
We emphasize that for the last display we don't use the analogue of $C_{\text{mid},n}$ from \eqref{eq:Cmid} to control $W_{l;q}$ for $l\in \{q+1,\dots, k-q\}$, but instead use that by \eqref{eq:Clog} and the niceness of $\frownacc{m}_{n,\frownacc{h}}$ we have
\[
\frac{|W_{k-q;q}^h|}{\sqrt{k\sqrt{\log(k)}}} \le (C_{\log}+C_{\frownacc{m}_{n,\frownacc{h}}})
\]
for all $k\in\N$. The last three displays and $\sup_n \sum_{k=0}^{\lfloor \log(n)^7\rfloor} e^{-\vartheta^\ast \frownacc{h}_n(k)}k^C<\infty$ for all $C>0$ then yields Lemma~\ref{Lem:SecMomBarrCompComplicated}.

The proof of Lemma~\ref{Lem:FMomHardBarrComp} goes the same way, using that by definition \eqref{eq:qkend}, we have
\begin{align*}
q_{k,\text{end}}(x) &= P^{\zeta}\left[\mathcal{B}_{\{n-k-1,\dots, n\}, -(\smileacc{m}_{n,\smileacc{h}}(\cdot)-\smileacc{m}_{n,\smileacc{h}}(n-k-1))}^{x,J_{y_0}}\big(H_\cdot-H_{n-k-1}-(\vartheta^\ast)^{-1}(K_\cdot^\zeta-K_{n-k-1}^\zeta) \big)\right]\\
&\ge P^{\zeta}\left[\mathcal{B}_{\{0,\dots, k+1\}, -(\smileacc{m}_{n,\smileacc{h}}(n-\cdot)-\smileacc{m}_{n,\smileacc{h}}(n-k-1)}^{y_0-1,J_{x}}\big(H_{n-\cdot}-H_{n-k-1})-(\vartheta^\ast)^{-1}(K_{n-\cdot}^\zeta-K_{n-k-1}^\zeta) \big) \right]\\
&=  \widehat{p}_{0,k+1}^{\smileacc{m}_{n,\smileacc{h}},(n)}(y_0-1,J_{x}),
\end{align*}
where in the second step we reversed time and in the last step we used Definition \eqref{eq:pt1t2TimeReversed}.
\end{proof}
\addcontentsline{toc}{section}{Bibliography}%
\printbibliography

\section*{Symbol Index}
\addcontentsline{toc}{section}{Symbol Index}%
\setlength\emergencystretch{\hsize}\hbadness=10000
\subsection*{Model and general objects}\strut
\vspace*{-.5cm}
\begin{longtable}{p{2.5cm}p{0.6\linewidth}p{0.1275\linewidth}}
Symbol & Description & Reference \\ \hline
$(B_s^{\zeta})_{s\ge 0}$ & Gaussian process with mean zero and local variance $\chi_k^2$ on $[k-1,k]$ & Below \eqref{eq:TiltPara}\\
$\mathcal{B}_{I,f}^{y,J}(Z_{\cdot})$& barrier event & \eqref{eq:GenBarrierEvent}\\
$c_{\chi}$, $C_{\chi}$ & constants such that for all $k\in\N$ $\chi_k^2\in [c_{\chi},C_{\chi}]$ $\mathbb{P}$-a.s. & Below \eqref{eq:TiltPara}\\ 
$\mathrm{ei}$ & essential infimum of $\xi(x)$ & \eqref{eq:EssinfEsssup}\\
$\mathrm{es}$& essential supremum of $\xi(x)$ & \eqref{eq:EssinfEsssup}\\
$\overline{\eta}$ & tilting parameter for which we have $\mathbb{E}[E^{\zeta,\overline{\eta}}[H_k]] = k/v_0$ & \eqref{eq:TiltPara}\\
$g_{\mu,\sigma^2}$ & Gaussian density with mean $\mu$, variance $\sigma^2$ & Section~\ref{Sec:Notation}\\
$\mathbf{H}_n$ & first time any particle hits $n$ & Section~\ref{Sec:Definitions and main result}\\
$H_k(X)$ & time at which $(X_s)_{s\ge 0}$ hits $k$ &  Above \eqref{eq:tilt} \\
$J_x$ & the interval $[x-1,x]$ & Section~\ref{Sec:Notation}\\
$K_n^\zeta$ & ``linear'' part of the centering around which $\mathbf{H}_n$ is tight & \eqref{Def:Knzeta}\\
$L_k^\zeta(\eta)$ &a logarithmic moment generating function  & \eqref{eq:LogMomMean}\\
$L_n^\zeta(\eta)$ & averaged logarithmic moment generating functions for $k\le n$ & \eqref{eq:LogMomMean}\\
$L(\eta)$ & expected logarithmic moment generating function & \eqref{eq:LogMomMean}\\
$L_{t_1,t_2}^{x_1,x_2}$ & piecewise linear interpolation between $x_1$ and $x_2$ on time-frame $t_1$, $t_2$ adapted to the local variance of $B_s^\zeta$ &  
\eqref{eq:AdaptedInter}\\
$\lambda$ & Lyapunov exponent & \eqref{eq:Lyap} \\ 
$M_t$ & maximal position at time $t$ & Section~\ref{Sec:Definitions and main result}\\
$m_n$ & centering making $\mathbf{H}_n-m_n$ tight & \eqref{Def:mn}\\
$\mathbb{N}_{\ge x}$ & natural numbers bigger than $x$ & Section~\ref{Sec:Notation}\\
$N(t)$ & particles alive at time $t$ & Section~\ref{Sec:Definitions and main result}\\
$\wt{O}_n(f(n))$ & functions being bounded by $C_n\cdot f(n)$, with $(C_n)_{n\in\N}$ tight w.r.t.\@ $\mathbb{P}$ & Section~\ref{Sec:Notation}\\
$\wt{O}_{\text{unif},g}^n(f(k,n))$& functions being bounded by $C_n\cdot f(n,k)$ uniformly in $k\le g$ with $(C_n)_{n\in\N}$ tight w.r.t.@ $\mathbb{P}$ & Section~\ref{Sec:Notation}\\
$p_n$ & probability that a Brownian motion stays above $W_\cdot$, ends near it & \eqref{Def:pn}\\
$P^{\zeta,\eta}$ & tilted random walk measure & \eqref{eq:tilt}\\
$P^{\zeta}$ & the measure $P^{\zeta,\overline{\eta}}$ for which $\mathbb{E}[E^{\zeta,\overline{\eta}}[H_k(X)]] = k/v_0$ &  Below \eqref{eq:TiltPara}\\
$\mathbb{P}$ & probability law of the environment $(\xi(x))_{x\in\Z}$ & Section~\ref{Sec:Definitions and main result} \\ 
$P_x^\xi$ & quenched law of our branching random walk given an environment $\xi$ with initial particle at $x$ & Section~\ref{Sec:Definitions and main result} \\
$P_x$ & annealed law of our branching random walk with first particle at~$x$ & Section~\ref{Sec:Definitions and main result}\\
$\sigma_k^2$ & variance of $H_k$ w.r.t.\@ $P^{\zeta}$ & Below \eqref{eq:TiltPara}\\
$\tau_k(X)$ & time $(X_s)_{s\ge 0}$ takes to hit $k$ starting from $k-1$ & Above \eqref{eq:tilt} \\
$\vartheta^\ast$ & $\mathrm{es}-\overline{\eta}$ & Below \eqref{eq:TiltPara}\\
$v_c$ & maximal $v$ so that $\lambda$ is linear on $[0,v]$ & Below \eqref{eq:Lyap} \\
$v_0$ & asymptotic velocity of $M_t$, has $\lambda(v_0) = 0$ & Below \eqref{eq:Lyap}\\
$(W_n)_{n\in\N}$ & fluctuation around the mean of the ``linear'' part of the centering making $\mathbf{H}_n$ tight & \eqref{Def:Wnzeta} \\ 
$(X_s)_{s\ge 0}$ & time continuous simple random walk with jump rate 1 and starting at $x$ w.r.t.\@ $P_x^\xi$ & Above \eqref{eq:tilt}\\ 
$X_{t_2;t_1}$ & increment of $X$ between times $t_1$, $t_2$, also used for other processes & Section~\ref{Sec:Notation}\\
$(\xi(x))_{x\in \Z}$ & random environment, $\xi(x)$ is the branching rate of particles at $x$ & Section~\ref{Sec:Definitions and main result}\\
$\chi_k^2$ & variance of $\tau_k$ w.r.t.\@ $P^{\zeta}$ & Below \eqref{eq:TiltPara} \\
$y_0$ & a fixed constant integer bigger than $e+1$ & Above~\eqref{Def:Wnzeta} \\
$(\zeta(x))_{x\in\Z}$& random environment shifted downwards by its essential supremum&  Above \eqref{eq:tilt}\\
$(\zeta^{(k)})_{x\in\Z}$ & random environment shifted in location by $k$ & Above \eqref{eq:qnsx}
\end{longtable}
\vspace*{-.65cm}
{\renewcommand{\arraystretch}{1.1}}
\subsection*{First moment calculation}\strut
\vspace*{-.5cm}
\begin{longtable}{p{2.5cm}p{0.6\linewidth}p{0.1275\linewidth}}
Symbol & Description & Reference \\ \hline
$A_n(y,x)$ & event that at least one particle $Y\in N(t_{n;y}^x)$ hits $n$ and has hitting times $H_k(Y)$, which stay above $t_{n;y}(k)$ & \eqref{eq:Anyx}\\
$\smileacc{h}_n(k)$ & curve of size roughly $-((n-k)\wedge k)^{1/6}$ added to the barrier in the first moment calculation & \eqref{eq:IncrPartBan},\eqref{eq:IncrPartEnd}, \eqref{eq:IncrPartMid} \\
$\smileacc{m}_{n,\smileacc{h}}(k)$ & $\smileacc{h}_n$ shifted by roughly $-(k/n)\log(p_n)$ & \eqref{Def:bprobUB}\\ 
$\smileacc{p}_n^x(y)$ & probability under $P^\zeta$ that $(H_k(X))_{k\le n}$ stays above $\smileacc{t}_{n;y}^0$ and has $H_n(X)\in [\smileacc{t}_{n;y}^x-1,\smileacc{t}_{n;y}^x]$ & \eqref{Def:bprobUB}\\
$q_{k,\text{end}}(x)$ & probability that $(x+H_l(X)-H_{n-k-1}(X))_{l\in \{n-k-1,\dots, n\}}$ stays above the barrier $\smileacc{t}_{n;0}(l)-\smileacc{t}_{n;0}(n-k-1)$ and has $H_n \le t_{n;x}^{y_0}$ & \eqref{eq:qkend}\\
$\smileacc{t}_{n;y}^x(k)$ & barrier used in the first moment shifted by $x-y$ & \eqref{Def:Deftnyx} \\
$\tau_n(y)$ & minimal $k\in[n]$ such that any particle of our BRW reaches $k$ before time $\smileacc{t}_{n;y}^{-1}$, or $n$ if no such time exists & \eqref{eq:tauny}\\
\end{longtable}
\vspace*{-.65cm}
\subsection*{Second moment calculation}\strut
\vspace*{-0.5cm}
\begin{longtable}{p{2cm}p{0.6\linewidth}p{0.17\linewidth}}
Symbol & Description & Reference \\ \hline
$A_n^y(s)$ & event that $(H_k(X))_{k\le k_s}$ stays above $t_{n;y}(k)$ and $H_{k_s}(X)\ge s$ & \eqref{eq:Anys}, Figure~\ref{fig:SecMom}\\
$\frownacc{h}_n(k)$ & curve of size roughly $((n-k)\wedge k)^{1/6}$ added to the barrier in the second moment calculation& \eqref{eq:SecMomhn1}, \eqref{eq:SecMomhn2}\\
$I_{k,n}$ & $[\frownacc{t}_{n;y}(k), \frownacc{t}_{n;y}(k+1)]$ & Below~\eqref{eq:qnsx}\\
$k_s$ & smallest $k$ such that $\frownacc{t}_{n;y}(k)>s$, set to be $n$ if no such $k$ exists& Definition~\ref{Def:TwoPartOnePartDefs}\\
$\frownacc{m}_{n,\frownacc{h}}(k)$ & $\frownacc{h}_n$ shifted by roughly $-(k/n)\log(p_n)$ & \eqref{Def:bprobLB}\\
$\frownacc{p}_{n}^{y_0}(y)$ & probability under $P^\zeta$ that $(H_k(Y))_{k\le n}$ stays above $\frownacc{t}_{n;y}^0(k)$ and $H_n(Y)\in[\frownacc{t}_{n;y}^{y_0}-1,\frownacc{t}_{n;y}^{y_0}]$ & \eqref{Def:tnyxk}\\
$q_{k,\text{start}}$ & probability under $P^\zeta$ that $(H_l(X))_{l\in\{0,\dots, k\}}$ stays above $\frownacc{t}_{n;y_0}(l)$ and has $H_k(X)\in [\frownacc{t}_{n;y_0}^x(k)-1, \frownacc{t}_{n;y_0}^x(k)]$ & \eqref{eq:qkstart}\\
$q_n(s,x)$ & probability that $(H_k((x+X_r)_{r\ge0}))_{k\in\{k_s,\dots, n\}}$ stays above $(\frownacc{t}_{n;y+s}(k))_{k\in\{k_s,\dots, n\}}$ and $H_n((x+X_r)_{r\ge0})\in [\frownacc{t}_{n;y+s}^{y_0}-1, \frownacc{t}_{n;y+s}^{y_0}]$ & \eqref{eq:qnsx}, Figure~\ref{fig:SecMom}\\
$q_n^\ast(k)$ &  & \eqref{eq:qnastk}\\
$\frownacc{t}_{n;y}^x(k)$ & barrier used in the second moment shifted by $x-y$ & \eqref{Def:tnyxk}\\
$Z_n(y)$ & number of particles $Y\in N(\frownacc{t}_{n;y}^{y_0})$ for which $(H_k(Y))_{k\le n}$ stays above $\frownacc{t}_{n;y}^0(k)$ and $H_n(Y)\in[\frownacc{t}_{n;y}^{y_0}-1, \frownacc{t}_{n;y}^{y_0}]$ &\eqref{Def:Zny}
\end{longtable}
\vspace*{-.65cm}
\subsection*{Barrier computations}\strut
\vspace*{-.5cm}
\begin{longtable}{p{0.2275\linewidth}p{0.5\linewidth}p{0.15\linewidth}}
Symbol & Description & Reference \\ \hline
 $A_0,\dots, A_3$, $B_0(u),\dots, B_2(u)$, $B_3$ & two partitions of $[0,\infty)$& \eqref{eq:CutPartition}\\
  $(c_k)_{k\in\N}$ & dyadic scale & Definition~\ref{Def:KMTBasicDef}\\
$C_h$ &  constant related to the (local) growth of $(h_n)_{n\in\N}$, finite iff $(h_n)_{n\in\N}$ is nice & \eqref{Def:Ch}\\
$C_{\log}$ &constant related to the growth of $(W_k)_{k\in\N}$ on different time-scales & \eqref{eq:Clog}\\
 $C_2(W,h,n)$, $C(\lambda)$ & constants related to a change of starting points in $p_{c_k,n}^{B,h\pm g_k}$ & Lemma~\ref{Lem:StartMove}\\
 $C_{\text{mid},n}(W^h)$ & constant related to the growth of $W_{k;n/3}$ on $\{n/3,\dots, 2n/3\}$ & \eqref{eq:Cmid}\\
 $g_k$, $g_k^{(n)}$ & a helper-function growing logarithmically & \eqref{eq:gknj}\\
 $\gamma_n$, $\widehat{\gamma}_n$, $\wt{\gamma}_n$ & tight sequences of exponents for polynomial lower bounds on $p_{c_k,c_{k+1}}^h(y,I)$, $p_{0,k}^h(y,I)$, $p_{c_k,l}^{B,h}(y_k, I_l)$ & Lemma~\ref{Lem:UnifCrudeLB}\\
 $k_0(n)$ & maximal $k$ for which $c_k\le \min\{n^{1/6},n/2\}$ & Definition~\ref{Def:KMTBasicDef}\\
$\bar{H}_k(X)$ & centered version of $H_k(X)$ & Definition~\ref{Def:KMTBasicDef}\\
$p_{t_1,t_2}^h(y,I)$ & probability w.r.t.\@ $P^{\zeta}$ that on $\{t_1,\dots, t_2\}$ $y+\bar{H}_{k;t_1}$ stays above $W_{k;t_1}^h$ and ends in $I+W_{t_2;t_1}^h$ & \eqref{eq:pathprobsDef} \\
$p_{t_1,t_2}^{B,h}(y,I)$ & probability w.r.t.\@ $P^{\zeta}$ that on $\{t_1,\dots, t_2\}$ $y+B_{k;t_1}^\zeta$ stays above $W_{k;t_1}^h$ and ends in $I+W_{t_2;t_1}^h$ & \eqref{eq:pathprobsDef}\\
$p_{t_1,t_2}^{B,h}(y|x)$ & probability w.r.t.\@ $P^{\zeta}$ that a discrete Brownian bridge (inheriting its variance structure from $B^{\zeta}$) from $y$ to $x$ stays above $W_{k;t_1}^h$ & \eqref{eq:pathprobsDef}\\
$p_{t_1,t_2}^h(y,z)$, $p_{t_1,t_2}^{B,h}(y,z)$ & density with respect to the Lebesgue measure corresponding to $p_{t_1,t_2}^h(y,\cdot)$, $p_{t_1,t_2}^{B,h}(y,\cdot)$ & Below \eqref{eq:pathprobsDef}\\
$\widehat{p}_{t_1,t_2}^{h,(n)}(y,I)$, $\widehat{p}_{t_1,t_2}^{B,h,(n)}(y,I)$ & time reversed versions of $p_{t_1,t_2}^h(y,I)$ and $p_{t_1,t_2}^{B,h}(y,I)$ & \eqref{eq:pt1t2TimeReversed}\\
$\bar{\tau}_k(X)$ & centered version of $\tau_k(X)$ & Definition~\ref{Def:KMTBasicDef}\\ 
$W_k^h$ & $W_k+h_n(k)$ & Definition~\ref{Def:KMTBasicDef}
\end{longtable}
\end{document}